\documentclass [12pt] {report}
\usepackage[dvips]{graphicx}
\usepackage{epsfig, latexsym,graphicx}
\usepackage{amssymb}
\newcommand {\D}[2] {\displaystyle\frac{\partial{#1}}{\partial{#2}}}

\newcommand {\al} {\alpha}

\newcommand {\ga} {\gamma}
\newcommand {\la} {\lambda}

\newcommand {\si} {\sigma}

\newcommand {\de} {\delta}
\newcommand {\prtl} {\partial}
\newcommand {\fr} {\displaystyle\frac}
\newcommand {\wt} {\widetilde}
\newcommand {\be} {\begin{equation}}
\newcommand {\ee} {\end{equation}}
\newcommand {\ba} {\begin{array}}
\newcommand {\ea} {\end{array}}
\newcommand {\bp} {\begin{picture}}
\newcommand {\ep} {\end{picture}}
\newcommand {\bc} {\begin{center}}
\newcommand {\ec} {\end{center}}
\newcommand {\bt} {\begin{tabular}}
\newcommand {\et} {\end{tabular}}
\newcommand {\lf} {\left}
\newcommand {\rg} {\right}

\newcommand {\cF} {{\cal F}}

\newcommand {\cR} {{\cal R}}
\newcommand {\cS} {{\cal S}}

\newcommand {\ses} {\medskip}

\newcommand {\e} {\mathop{\rm e}\nolimits}

\newcommand {\bibit} {\bibitem}
\newcommand {\nin} {\noindent}

\newcommand {\Rho} {\mbox{\large$\rho$}}

\newcommand {\ka} {\varkappa}

\newcommand {\Cos} {\mathop{\rm Cos}\nolimits}

\newcommand {\Sin} {\mathop{\rm Sin}\nolimits}

\newcommand {\Ka} {\mbox{\Large$\varkappa$}}

\usepackage{epsfig, latexsym,graphicx}
 \usepackage{amsmath}

\newcommand {\cD} {{\cal D}}

\def\2#1#2#3{{#1}_{#2}\hspace{0pt}^{#3}}
\def\3#1#2#3#4{{#1}_{#2}\hspace{0pt}^{#3}\hspace{0pt}_{#4}}
\newcounter{sctn}
\def\sec#1.#2\par{\setcounter{sctn}{#1}\setcounter{equation}{0}
                  \noindent{\bf\boldmath#1.#2}\bigskip\par}

\begin {document}

\begin {titlepage}

\vspace{0.1in}

\begin{center}

{\Large \bf  Finsleroid  gives rise to the angle-preserving connection}

\end{center}

\vspace{0.3in}

\begin{center}

\vspace{.15in} {\large G.S. Asanov\\} \vspace{.25in}
{\it Division of Theoretical Physics, Moscow State University\\
119992 Moscow, Russia\\
{\rm (}e-mail: asanov@newmail.ru{\rm )}} \vspace{.05in}

\end{center}

\begin{abstract}

\ses

\ses

The Finslerian unit ball is called the {\it Finsleroid} if the covering indicatrix
is a space of constant curvature.
We prove that Finsler spaces with such indicatrices possess the remarkable property
that the tangent spaces are conformally flat with the conformal factor of the power
dependence on the Finsler metric function.
It is amazing but the fact that in such spaces the notion of the two-vector angle
defined by the geodesic arc on the indicatrix
can readily be induced from the Riemannian space obtained upon the conformal transformation,
which opens up the straightforward way to induce also the connection coefficients and
the concomitant curvature tensor.
Thus, we are successfully inducing the Levi-Civita connection
 from the Riemannian space into the
Finsleroid space, obtaining the isometric connection.
 The resultant connection coefficients are not symmetric.
However, the metricity condition
 holds fine, that is, the produced covariant derivative
of the Finsleroid metric tensor vanishes identically.
The particular case underlined by the axial Finsleroid of the
 ${\mathbf\cF\cF^{PD}_{g}}$-type is explicitly evaluated in detail.

\medskip


\end{abstract}

\end{titlepage}

\vskip 1cm

\ses

\ses

\setcounter{sctn}{1} \setcounter{equation}{0}

\nin
  {\bf 1. Introduction and motivation}

\ses

\ses

What matter makes the Levi-Civita connection canonical in the Riemannian geometry?
The angle-parallelism property is the matter, namely  that the angle of two vectors
is left unchanged under the  parallel transportation of the vectors with that connection.
Can the property  be put ahead  also in the Finsler geometry?

In any tangent space $T_xM$
of a given   Finsler space  ${\cal F}^N$
introduced on
 an $N$-dimensional differentiable manifold $M$ of points $x\in M$,
 the angle
$ \al_{\text{Finslerian}}(x,y_1,y_2)  $
between two  tangent vectors $y_1,y_2\in T_xM$ can  be defined by means of
length
of the respective geodesic arc ${\cal A}(x,l_1,l_2)$ which
joins the intersection points of
the directions of vectors $y_1,y_2\in T_xM$ with the indicatrix
        ${\cal I}_x\subset T_xM$ supported by the point $x\in M$
and belongs to the indicatrix;
$l_1=y_1/F(x,y_1)$
and
$l_2=y_2/F(x,y_2)$
are respective unit vectors, and $F$ denotes the Finsler metric function of the space
${\cal F}^N$.
With such an angle, it can be hoped to prolong   the Levi-Civita connection
 from the Riemannian geometry to the Finsler geometry
 to arrive at
 the Finsler notion of the parallel connection with just the same meaning as the connection
notion
 that faces us
in the Riemannian geometry.
Namely, the parallel connection operates such that
under the parallel transportation generated by the connection
the vectors $y_1,y_2\in T_xM$ are transported in the parallel way and the
length of the arc ${\cal A}(x,l_1,l_2)$
does not change.
 However, the hard analytical difficulties can  arise that the
dependence of $ \al_{\text{Finslerian}}(x,y_1,y_2)  $
on $y_1,y_2$    cannot be written in an explicit form, except for (rare?)
 particular cases
of the  space  ${\cal F}^N$.
The desired prolongation  may exist but  not be explicit!
Various interesting investigations of  the  connection and angle
can be found  in the Finsler geometry literature
(see, for example, [1-9]).

Below,
we come to the Finsleroid-produced  parallel transports
which do not change neither the norms of vectors
nor the
Finsleroid-produced angles between vectors. Therefore,
  they are isometric:
  the produced displacements of the tangent spaces
along  curves of the underlined manifold $M$
are isometries.  In particular, they keep the  Finsleroid indicatrices into the  Finsleroid
indicatrices.
At the same time, in contrast to the Levi-Civita connection of the Riemannian geometry,
the Finsleroid-produced  parallel transports
are not linear in general, namely, the connection coefficients
processing  that transports depend
on tangent vectors $y\in M$ in a nonlinear way, except for possible particular cases.
The entailed connection
coefficients
 are {\it not} constructive
  from the Finslerian metric tensor and the first derivatives
of the tensor
--- they are obtained
 from the parallel transportation of the two-vector angle.

{%\pgbrk}

We  shall construct the angle with the help of the indicatrix arcs
as follows.

The embedded position  of the indicatrix
${\cal I}_x\subset T_xM$ in the tangent Riemannian space $\{T_xM, g^{\{F\}}(x,y)\}$
(where $g^{\{F\}}(x,y)$ denotes the Finslerian metric tensor with $x$ considered fixed
and $y$ used as being the variable) induces the Riemannian metric on the indicatrix
through  the well-known method  (see, e.g., Section 5.8 in [1])
and in this sense makes the indicatrix
 a Riemannian space.
Let $ U_x $ be a simply connected and geodesically complete region on the indicatrix
${\cal I}_x$ supported by a point $x\in M$.
Any point pair $u_1,u_2\in U_x$ can be joined by the respective
arc ${\cal A}(x,l_1,l_2) \subset {\cal I}_x$
of the  Riemannian  geodesic  line
drawn on $U_x$.
By identifying the  length of the arc  with the angle notion we arrive at
the {\bf geodesic-arc angle}
$ \al_{\{x\}}(y_1,y_2)$,
where $y_1,y_2\in T_xM$  are  two vectors issuing from the origin $0\in T_xM$ and
possessing the property that their direction rays $0y_1$ and $0y_2$
intersect
 the indicatrix at the  point pair $u_1,u_2\in U_x$.
Denoting
the respective Riemannian length of the geodesic arc ${\cal A}(x,l_1,l_2)$
by $s$,
we obtain
\be
\al_{\{x\}}(y_1,y_2)=s.
\ee

 With this angle notion, we can naturally introduce the scalar  product:
$$
<y_1,y_2>_{\{x\}} = F(x,y_1)F(x,y_2) \al_{\{x\}}(y_1,y_2).
$$

{%\pgbrk}

We may use   geodesic-arc  angle to construct
the {\bf  geodesic-arc sector}
$ {\cal S}_{\{x\}}(l_1,l_2)\subset T_xM$,
where $l_1,l_2\in T_xM$  are two unit vectors  issued from the origin $0\in T_xM$
(and pointed to the indicatrix). The sector said is the surface  swept by the unit vector
$l\in T_xM$ when its end runs along the geodesic arc from the point
$u_1$ to the point $u_2$.

Let us call the arc-curve joining the end points $u_1$ and $u_2$ of the sector
the {\it arc-side} of the sector.
When a point moves along  an arc-side parameterized by the geodesic length $s$,
 we obviously have
\be
(ds)^2=g_{mn}dl^mdl^n,
\ee
which can be written as
\be
g_{mn}\fr{dl^m}{ds}  \fr{dl^n}{ds}=1.
\ee

  Each
sector $ {\cal S}_{\{x\}}(l_1,l_2)$
can naturally be coordinatized by the pair
$z^1,z^2$ with
\be
z^1=F, \quad z^2=s,
\ee
 where $F$ is the value of the Finslerian metric function,
and
$s$ is the length of the geodesic arc  which goes from the ray direction of the left vector $l_1$
to the consideration point.
Let us consider the embedding $y=y(z^1,z^2)$ of the sector
$ {\cal S}_{\{x\}}(l_1,l_2)$
in the tangent Riemannian space
$\{T_xM, g^{\{F\}}(x,y)\}$.
Projecting the  tensor $g^{\{F\}}$
on the sector gives rise to the
{\it intrinsic metric tensor}, to be denoted by
$i_{\cal S}$.
The tensor $i_{\cal S}$ has the components:
$\{i_{11}, i_{12}=i_{21}, i_{22}\}$.
In terms of a local coordinate system $\{x^m\}$, we have
$y^m=y^m(z^1,z^2)$ and
\be
i_{11} =g_{mn}y^m_1y^n_1, \qquad
i_{12} =g_{mn}y^m_1y^n_2, \qquad
i_{22} =g_{mn}y^m_2y^n_2,
\ee
where $y^m_1=\partial y^m/\partial z^1$
and
$y^m_2=\partial y^m/\partial z^2$, and $g_{mn}$ are the Finslerian metric tensor components.
The equality   $y^m=Fl^m$  just entails that
$y^m_1=l^m$ and, therefore,
$i_{11} =1.$
Also,
$
l_my^m_2=0
$
(because  $l_m$ is the gradient vector $\partial F/\partial y^m$ and $F$ is independent of $s$),
which makes us conclude that $i_{12}=0$.
Finally, noting that in the case of an arc intersecting the ray at the distance $F$
from the origin  $0\in T_xM$ the equality  (1.3) must be modified to read
$F^2g_{mn}(dy^m/{ds})(dy^n/{ds})=1$,
we arrive at
$g_{mn}y^m_2y^n_2=1$.

Thus we have observed that the coordinates (1.4)
introduce the {\it orthogonal coordinate system} on the    geodesic-arc  sectors
and the
{\it intrinsic metric tensor of the   geodesic-arc  sector
is  Euclidean}:
\be
i_{11} =1, \qquad   i_{12} =0, \qquad
i_{22} =1.
\ee
In this sense, {\it each   geodesic-arc  sector is a Euclidean space.}

This observation just entails that in any Finsler space
{\it the area
$ ||{\cal S}_{\{x\}}(l_1,l_2)||$
of the   geodesic-arc  sector }
is presented by the  formula
\be
\bigl|\bigr.\bigl|\bigr. {\cal S}_{\{x\}}(l_1,l_2)\bigr|\bigr.\bigl|\bigr.=\fr12
\lf(\al_{\{x\}}(y_1,y_2)\rg)^2
\ee
which is {\it faithfully valid} in  all the Riemannian as well as Finslerian spaces.

{%\pgbrk}

Our underlined  idea  is to use the   geodesic-arc angle to
generate the notion of the
 angle-preserving  connection.
We shall  always assume that the connection is metrical.
Each  respective Finsleroid-produced  parallel transport
along a curve of the underlined manifold $M$    is an isometry for
the involved angles and, therefore, for   the involved
geodesic-arc sectors.

{%\pgbrk}

We shall develop  the idea by specifying the Finsler space as follows.

If an $N$-dimensional  Finsler space
${\cal F}^N$
is such that  the indicatrices of the space  possess the property
of constant curvature (with respect to the Riemannian metric induced from
the tangent Riemannian space
$\{T_xM, g^{\{F\}}(x,y)\}$), we call
${\cal F}^N$
the {\it Finsleroid-Finsler space},
 denote  the fundamental metric function of the space
by
\be
K=K(x,y),
\ee
and apply the following definitions.

\ses

\ses

 {\large  Definition}.  Within  any tangent space $T_xM$, the  metric function $K(x,y)$
 produces the {\it  Finsleroid}
\be
 \cF_{\{x\}}:=\{y\in    \cF_{\{x\}} : ~ y\in T_xM , K(x,y)\le 1\}.
\ee

\ses

 \ses

 {\large  Definition}. The {\it  Finsleroid Indicatrix}
\be
 {\cal  IF}_{\{x\}} \subset T_xM
 \ee
  is the boundary $\partial  \cF_{\{x\}}$
 of the  Finsleroid, that is,
\be
{\cal  IF}_{\{x\}} :=\{y\in  {\cal  IF}_{\{x\}} : ~ y\in T_xM, K(x,y)=1\}.
\ee

\ses

 \ses

The strong convexity and the positive homogeneity are assumed.

Let us also assume that  the  manifold $M$ can be endowed with a Riemannian metric tensor,
$a_{mn}(x)$  in terms of local coordinates $\{x^m\}$.
Then
in addition to the Finsleroid indicatrix  (1.11) we can bring to consideration
  the sphere obtainable from
 the {\it associated Riemannian space}
\be
{\cal R}^N:~=\{M,a_{mn}(x)\},
\ee
in accordance with

 \ses

 \ses

 {\large  Definition}. The {\it   sphere}
 \be
   \cS_{\{x\}} \subset T_xM
\ee
is defined by
\be
 \cS_{\{x\}} :=\{y\in   \cS_{\{x\}} : ~ y\in T_xM, \, a_{mn}(x)y^my^n=1\}.
\ee

\ses

 \ses

Our input stipulation is that the Finsler space
${\cal F}^N$
be  conformally   isomorphic to the
Riemannian space:
\be
{\cal F}^N=\Ka\cdot {\cal R}^N
\ee
in accordance with see (2.1)-(2.3).
 We shall see that such Finsler spaces are  the Finsleroid-type spaces.

{%\pgbrk}

 In Section 2  the basic theorems are sketched to enlighten the numerous
 remarkable properties which are
shown by  such Finsler spaces.
It proves that the
 $\Ka$-transformation (1.15)
 must involve the conformal factor which
is
of a power dependence on  the Finsler metric function,
as shown by (2.5).
In such Finsler spaces ${\cal F}^N$,
we have universally
the ratio
\ses\\
\be
{\fr{ \text{Area of the Finsleroid   geodesic-arc  sector}}
{ \text{Area of the Riemannian   geodesic-arc  sector}}}\Bigl|_{x\in M}\Bigr.
= \fr1{h^2(x)}
\ee
(at any admissible data of the leg-vectors $y_1,y_2$ of the sector),
where the denominator relates to the space
${\cal R}^N$
and $h^2(x)$ is the value of curvature of the indicatrix supported by  point $x$ (see (2.6)).

\ses

By inducing the angle $\al_{\text{Finslerian}}=\Ka\cdot \al_{\text{Riemannian}}
$
we obtain the remarkable equality
\be
\al_{\text{Finsleroid space}}
 =  \fr1h
 \al_{\text{Riemannian space}}
\ee
(see (2.6) and (2.7)).

At each point $x\in M$,
the ratio
$
{AREA_{\text{Finsleroid Indicatrix}}}
/{VOLUME_{\text{Finsleroid }}}
$
proves to be of  the universal value in each   dimension  $N $
(is independent of $h$), so that the ratio is exactly the same as it holds in the Riemannian limit
(that is, when $h=1$).
This property  is lucidly described by the formulas (3.15)-(3.20).
\ses

In Sections 3 and 4  we explain how the sought  connection
coefficients and the curvature tensor are obtainable
on the basis of  the
transformation (1.15),
provided $h=const$.
The connection coefficients involve  the $\Ka$-transition functions,
  {\it not} being
obtainable from the Finslerian metric tensor and the first derivatives
of the tensor in any algebraic manner.
The connection coefficients are not symmetric.
They are  non-linear in general regarding the $y$-dependence.
However, the {\it metricity},
that is, the condition that the covariant derivative
of the Finslerian metric tensor
be the nought (see (2.12)), can well be fulfilled.
The method is to postulate the {\it transitivity  of the covariant derivative}
under the $\Ka$-transformation. Analytically, the transitivity reads (4.9),
 and is fulfilled when the vanishing
(4.4)-(4.5) is postulated.
The angle (1.17) is preserved under  respective covariant displacements:
\be
d\al_{\text{Finsleroid space}}=0, \quad \text{if} \quad h=const,
\ee
where
 $d=dx^id_i$ and   $d_i\al_{\text{Finsleroid}}$
symbolizes  the left-hand part of (3.25).

In this way  we are quite able to
 successfully construct  in the Finsleroid  space ${\cal F}^N$
the
{\it angle-preserving connection}, to be denoted by  ${\cal FC}$,
by  adhering    faithfully at the method
\be
\Bigl\{\text{The Finsleroid-space connection} ~ {\cal FC}
\Bigr\}=\Ka\cdot
\Bigl\{\text{The Levi-Civita connection}  ~ {\cal LC}\Bigr\},
\ee
where ${\cal LC}$ is the canonical connection in the associated Riemannian space (1.13)
(the  connection coefficients of the
${\cal LC}$
are   the Christoffel symbols $a^m{}_{ij}=a^m{}_{ij}(x)$
constructed from the Riemannian metric tensor $a_{mn}(x)$, in accordance with
 (3.22)).

{%\pgbrk}

In Section 5,
we fix  a tangent space
and
consider the tensor $k_{ij}$ obtainable from the tensor
 $g_{ij}$ by performing the conformal transformation.
 We observe that
 the property of vanishing   the  curvature tensor produced by  $k_{ij}$
 is arisen upon fulfilling a simple ODE, which can explicitly
 be solved to establish  the  theorem 2.2.

 In Section 6,
we confine the consideration to
the   ${\mathbf\cF\cF^{PD}_{g} }$-space,
in which
the Finsleroid is of the axial type;
 $g$ denotes
the characteristic parameter; the upperscript  $\{PD\}$ means {\it positive-definite.}
The connection coefficients  are found explicitly.
They depend on vectors $y$ in a non-linear way in  dimensions $N\ge3$.
In the dimension $N=2$, however,
the connection is linear.
The structure of the appeared curvature tensor
$\Rho_k{}^n{}_{ij}$ has been elucidated, resulting in the explicit representation
(6.77)-(6.78).
The square of the tensor is given by the simple formula
(6.79).
\ses

Since the formula $h=\sqrt{1-(1/4)g^2}$  is  applicable in the
 ${\mathbf\cF\cF^{PD}_{g}}$-Finsleroid space,
 from the universal law (1.16)   we may conclude that
\ses\\
\be
{\fr{ \text{Area of the  ${\mathbf\cF\cF^{PD}_{g}}$-Finsleroid   geodesic-arc  sector}}
{ \text{Area of the Riemannian   geodesic-arc
 sector}}}\Bigl|_{x\in M}\Bigr.
>1 \qquad \text{when} \quad g\ne 0.
\ee

\ses

\ses

The angle measured by the lengths of geodesic arcs on the indicatrix
was found for   the ${\mathbf\cF\cF^{PD}_{g}}$-space
in the  work [7,8].
The underlying idea was to derive the angular measure from
 the solutions to the  respective geodesic equations.
The solutions have been derived in simple explicit forms.
This angle given by the formula (6.30)  can also be
obtained from the relationship (1.17),
because  in  the ${\mathbf\cF\cF^{PD}_{g}}$-space
we are able to propose the explicit knowledge
of the respective
$\Ka$-transformation,
 which is given by the formula
(6.26).
 The involved preferred vector field $b^i(x)$
  as well as the opposed vector  prove to be  the proper elements
of the   $\Ka$-transformation (see (6.28)).

{%\pgbrk}

Our  evaluations will everywhere  be of  local nature. However,
there   exists  a simple possibility to elucidate the global structure
of the ${\mathbf\cF\cF^{PD}_{g} }$-Finsleroid indicatrix.
Indeed,
in the ${\mathbf\cF\cF^{PD}_{g} }$-Finsleroid space
 the desired $\Ka$-transformation (1.15)
can be explicitly  given by means of the substitution $\zeta^i=\zeta^i (x,y)$
indicated in (6.26).
Inserting these $\zeta^i$
in the
associated Riemannian metric
$
S(x,\zeta)=\sqrt{a_{mn}(x)\zeta^m\zeta^n}
$
entails the equality (3.7), thereby producing the
 metric function
$K(x,y)$
of the ${\mathbf\cF\cF^{PD}_{g} }$-Finsleroid space.
Let us consider  the  sphere $  \cS^{\{\zeta\}}_{\{x\}} \subset T_xM$
in terms of the variables $\zeta$:
$
 \cS^{\{\zeta\}}  _{\{x\}} :=
 \{\zeta \in   \cS^{\{\zeta\}}_{\{x\}} : ~ \zeta\in T_xM,\,S(x,\zeta)=1\}.
$
The equality (3.7) manifests that
 the transformation
(6.26) maps regions of
the ${\mathbf\cF\cF^{PD}_{g} }$-Finsleroid indicatrix
 in regions of  the sphere
$  \cS^{\{\zeta\}}_{\{x\}} $.
Also,
the direction of the Finsleroid-axis vector $b^i$
as well as the direction of the opposed vector
$-b^i$
are left invariant under the transformation (6.26) (see (6.28) and (6.29)).
Denoting  by
 $\zeta^{\{\text{North}\}}$,
 {\it resp.} by
$\zeta^{\{\text{Sourth}\}}$, the point  which is obtained in intersection
of the direction of $b^i$,
{\it resp.} by $-b^i$,
with the  sphere,
we can  consider the pointed spheres
\be
 S^{\{+\}}_{\{x\}} = S^{\{\zeta\}}_{\{x\}}\setminus  \zeta^{\{\text{Sourth}\}},
 \qquad
 S^{\{-\}}_{\{x\}} = S^{\{\zeta\}}_{\{x\}}\setminus  \zeta^{\{\text{North}\}},
\ee
in which
the south pole,
 {\it resp.} the north pole,
is deleted.
The regions
$
  C^{\{+\}}_{g;\,\{x\}} = {\Ka}^{-1}  S^{\{+\}}_{\{x\}}
  $
  and
$  C^{\{-\}}_{g;\,\{x\}} = {\Ka}^{-1} S^{\{-\}}_{\{x\}}
$
may be used to yield  two {\it covering charts}
 for the indicatrix; they can be characterized by the angle ranges indicated in (6.22).
Similarly to the Riemannian case proper,
we need two charts to cover the indicatrix.
In its sense and role, the south chart  $ C^{\{-\}}_{g;\,\{x\}} $
is entirely similar to the ordinary Euclidean chart
obtainable by means of the so-called stereographic projection.
Then it can readily be seen that
 the substitution $\zeta^i=\zeta^i (x,y)$
indicated in (6.26) acts
{\it diffeomorphically on each the  chart}:
\be
  C^{\{+\}}_{g;\,\{x\}}  \stackrel{\Ka}{\Longleftrightarrow} S^{\{+\}}_{\{x\}},
\qquad
  C^{\{-\}}_{g;\,\{x\}}  \stackrel{\Ka}{\Longleftrightarrow} S^{\{-\}}_{\{x\}}.
\ee
Thus,  for the
 ${\mathbf\cF\cF^{PD}_{g} }$-Finsleroid space
  the   $\Ka$-transformation and, therefore,
the representations of the connection coefficients
 and the curvature tensor
obtained in Section 6,
are  meaningful globally   regarding the $y$-dependence.

We are also entitled to say that {\it the  ${\mathbf\cF\cF^{PD}_{g} }$-Finsleroid indicatrix
 is globally isometric to the Euclidean sphere of the radius $r=1/h \equiv 1/
 \sqrt{1-(1/4)g^2}$.}

Below when mentioning
the  ${\mathbf\cF\cF^{PD}_{g} }$-Finsleroid space,
we shall imply that we work on the upper regions
$  C^{\{+\}}_{g;\,\{x\}}$, unless otherwise stated explicitly.
The development of  extensions to
regions of $  C^{\{-\}}_{g;\,\{x\}}$ is, of course, a straightforward task.

{%\pgbrk}

In the Euclidean and Riemannian geometries,
an important role is played by the spherical coordinates.
Their use enables one to conveniently represent vectors,
evaluate squares and volumes,
study curvature of surfaces,
in many cases simplify consideration
and solve rigorously equations,
and also introduce and use
various trigonometric functions.
In the context of the
 ${\mathbf\cF\cF^{PD}_{g}}$-space
 theory, such coordinates
 can readily be arrived at.
 In the three-dimensional case, $N=3$,
they are given by (7.2), entailing the convenient representation (7.11)
for vectors as well as the generalized trigonometric functions indicated in (7.12).
The respective
squared linear element  $ds^2$ has been explicitly evaluated
to read
(7.14)
which lucidly manifests
the conformal-flat nature of tangent spaces as well as the validity
 of the key formulas (2.5)-(2.7).
With the help of these coordinates,
the equations for
the   arc
${\cal A}(x,l_1,l_2)$ can explicitly be integrated
in the  convenient  form  (7.26)-(7.27),

\ses

{%\pgbrk}

In Appendix A,
the basic representations of  objects of the
 ${\mathbf\cF\cF^{PD}_{g}}$-space
 are summarized.
In Appendix B,  the involved connection coefficients are evaluated.
In Appendix C,  many steps of calculation of the curvature tensor are presented.
In Appendix D, we evaluate the coefficients which enter the transformation of the curvature
tensor into the Riemannian space.
In the last Appendix E,
we show the explicit representation of components for the metric tensor of
the  ${\mathbf\cF\cF^{PD}_{g}}$-space in the fixed tangent space.

\ses

We are interested mainly in spaces of the dimension $N\ge3$.
The two-dimensional case has been studied in the preceding work [9].

{%\pgbrk}

\ses

\setcounter{sctn}{2} \setcounter{equation}{0}

\nin
  {\bf 2. Synopsis of main assertions}

\ses
\ses

We start with  the following idea of specifying the notion of a Finsler space.

\ses

\ses

INPUT STIPULATION.
A Finsler space ${\cal F}^N$ is  {\it conformally   isomorphic} to
the Riemannian space
${\cal R}^N$:
\be
{\cal F}^N=\Ka\cdot {\cal R}^N: ~ ~ ~ \{g_{mn}(x,y)\}=\Ka\cdot \{t_{mn}(x,y)\}
~ ~ \text{with}
~ ~ t_{mn}(x,y) = k^2(x,y)a_{mn}(x),
\ee
where  it is assumed that the applied $\Ka$-transformation
 does not influence any point $x\in M$ of the base manifold  $M$.
It is also natural to require that
the $\Ka$-transformation sends  unit vectors to unit vectors:
\be
{\cal  IF}_{\{x\}}= \Ka \cdot    {\cal S}_{\{x\}}.
\ee
It looks interesting to specify the conformal multiplier to be an algebraic function
of the fundamental Finsler metric function $K=K(x,y)$ used in ${\cal F}^N$,
so that
\be
 k = \ka(x,K).
\ee
The smoothness of class $C^5$ regarding the $y$-dependence,
and of class $C^4$ regarding the $x$-dependence,
is necessary to require from the $\Ka$-transformations.
\ses

\ses

Under the above stipulation,
the tangent spaces to the Finsler space ${\cal F}^N$ are conformally flat.
At the same time, the   space  ${\cal F}^N$
is not conformal to any  Riemannian  space,
unless $\partial k/\partial K=0$.
In terms of local coordinates, the stipulation (2.1) is described by the formulas (3.1)-(3.5).

In Finsler spaces ${\cal F}^N$
fulfilling the conditions (2.1)-(2.3)
 we can measure the angle by the conventional value
$\al_{\text{Riemannian}}$ used in  the Riemannian space ${\cal R}^N$
and  obtain in  ${\cal F}^N$
the  induced angle
\be
\al_{\text{Finslerian}}=\Ka\cdot \al_{\text{Riemannian}}.
\ee

{%\pgbrk}

Attentive calculations performed in Section 5 result in the following

\ses

\ses

{\bf Theorem 2.1.}
{
\it
The claimed conditions}
 (2.1)-(2.3)
 {\it
 are realized if and only if
 the function $\ka$ is taken to be  }
\be
\ka=\fr1hK^{1-h} \quad \text {with} \quad h=h(x),
\ee
{\it where $h$ is a positive scalar.   }

\ses

\ses

The theorem is compatible with  the homogeneity,
namely we adopt

\ses

\ses

{\bf Homogeneity  condition.}
{\it
Action of the $\Ka$-transformation } (2.1) {\it on tangent vectors
 possesses the property  of the positive
homogeneity of degree
} $h$.

\ses

\ses

The formulas (3.3) and (3.4) proposed in  Section 3 yield
the analytical representation to the last condition.

Also,  the following theorem  can  be obtained  (see (5.6)).

\ses

\ses

{\bf Theorem 2.2.}
{\it
Under the conditions formulated in the preceding theorem,
 the indicatrix supported by a point $x\in M$ is of the constant curvature $h^2(x)$,
  such that    }
\be
 \cR_{\text{Finsleroid Indicatrix} }={\cal C}(x) ~~ \text{with} ~ ~ {\cal C}(x)= h^2(x).
\ee

\ses

\ses

The space of a  constant curvature ${\cal C}$ is
 realized on the sphere of the radius $r=1/\sqrt{\cal C}$.
Therefore, since we adhere at measuring the angle $\al$ by  the geodesic arc-length
on the indicatrix,
from the above formulas (2.1)-(2.6) we are entitled to conclude the following.

\ses

\ses

{\bf Theorem 2.3.}
{\it
The simple property   }
\be
\al_{\text{Finsleroid space}}
 = \fr1{\sqrt{\cal C}}
 \al_{\text{Riemannian space}}
\ee
{\it is valid for the angle}.

\ses

\ses

To surely recognize the validity of this theorem, it is sufficient to take a glance on the
equality (3.9) which represents infinitesimally the squared length of the geodesic arc
on the indicatrix.

Inverting the theorem 2.2 can be justified by taking into account the
derived formulas
(5.4)-(5.8), namely  the following theorem is fulfilled.

\ses

\ses

{\bf Theorem 2.4.}
{\it If the indicatrix supported by a point $x\in M$ is a space of  constant curvature,
then
the conformal property }
(2.1)
{\it holds at the point} $x$.

\ses

\ses

{%\pgbrk}

The sought Finsleroid connection
\be
{\cal FC}=\{N^k_j,D^k{}_{ij}\}
\ee
involves the coefficients
$N^k_j=N^k_j(x,y)$
and
$D^k{}_{ih}=D^k{}_{ih}(x,y)$
which are  required to construct the operator
\be
d_i=
\D{}{x^i}+N^k_i   \D{}{y^k}
\ee
and  the  covariant derivative $\cD_i$ which action
 in the Finsleroid  space    ${\cal F}^N$
is exemplified in  the conventional way:
\be
\cD_iw^n{}_m=   d_iw^n{}_m  +   D^n{}_{ih}w^h{}_m   -  D^h{}_{im}w^n{}_h,
\ee
where
$w^n{}_m=w^n{}_m(x,y)$ is an arbitrary differentiable (1,1)-type tensor.
In Section 4, we  subject the connection
${\cal FC}$
to the condition that the covariant derivative
$\cD_i$ obeys the transitivity rule (4.9)
and that the $h$ entering  (2.5) is independent of the points $x\in M$, so that $h=const$.
These conditions result in the following assertion.

\ses

\ses

{\bf Theorem 2.5.}
{\it When  $h=const$, the  vanishing set
\be
\cD_iK=0, \quad \cD_iy^j=0,  \quad \cD_iy_j=0 ,
\ee
and  the {\it metricity }
\be
\cD_ig_{jn}=0
\ee
hold, when
the relations
\be
D^k{}_{in} =-\D{N^k_i}{y^n}, \qquad N^j_i=-D^j{}_{ik}y^k
\ee
are used and
the coefficients
$N^k_i$
are constructed in accordance with the explicit formulas
}
 (3.30)-(3.32).
{\it
 The angle-preserving   property } (1.18)
{\it is entailed.
}

\ses

\ses

The curvature tensor
can be explicated from the commutator of the covariant derivative
(2.10),
according to  (4.12)-(4.18).

{%\pgbrk}

All the above formulas are  explicitly (and brightly!)
realized in the
 ${\mathbf\cF\cF^{PD}_{g}}$-Finsleroid space
  (Section 6), in which
   the     metric function  can conveniently be   introduced  by the representation
\be
K(x,y)=\sqrt{B(x,y)}\, e^{- (1/2) g(x)\chi(x,y)},
\ee
where the formulas  (6.1)-(6.6) are of value.
The scalar  $\chi$  thus appeared possesses the lucid geometrical meaning
of the azimuthal angle measured from the direction assigned by the input vector $b^i(x)$
(see (6.14)).

 The metric function  of the Finsleroid type has been first appeared in the paper [10],
at  but the Minkowskian level.
Namely, the consideration in [10]
was  subjected to the following assumptions. ($Z_1$) The sought metric function
$ \breve K(y)$ is positively homogeneous and smooth locally
of at lest class $C^4$.
($Z_2$) The indicatrix of $ \breve K(y)$ is a surface of revolution,
say around the direction of the $N$-th component $y^N$ of the tangent vector $y$,
in which  case it is convenient to introduce the representation $ \breve K(y)=y^NV(w)$,
 where the generating metric function
$V$ depends on a single argument $w$.
($Z_3$) The induced Riemannian curvature on the indicatrix
is of  a constant-curvature type.
Treating the condition to be a differential equation to find the function $V$, we can arrive after straightforward
calculations (which are not short) to the ODE (which is non-linear and of the second order) which governs the $V$.
It is a bit surprising but the fact that the ODE can explicitly be resolved at a local level to specify the function $V$.
($Z_4$) The obtained function $V=V(w)$ should obey the requirement that the entailed Finslerian metric tensor
be  positive definite.
The final  condition is ($Z_5$): The indicatrix is closed and  regular.
The function $ \breve K(y)$ obtainable in this way, after fulfilling all the conditions
$(Z_1)$-$(Z_5)$,
 is just the Minkowskian version of   the ${\mathbf\cF\cF^{PD}_{g}}$-Finsleroid
  metric function $K$
given by the formulas (6.1)-(6.5).

\ses

Thus, from the standpoint of the indicatrix geometry,
 the ${\mathbf\cF\cF^{PD}_{g}}$-space
 occupies  a {\it unique position} in the class of the Finsler spaces
 specified by the condition
that the Finsler metric function
$K=K(x,y)$ be of the functional dependence
\be
K =\Phi \Bigl(g(x), b_i(x), a_{ij}(x),y\Bigr),
\ee
where $g(x)$ is a scalar, $b_i(x)$ is a covariant vector field, and
$a_{ij}(x)$ is a Riemannian metric tensor.
No Finsler metric function  allowing for the representation (2.15)
can meet the requirements that the entailed indicatrix is of constant curvature
and the smoothness class $C^k$ is attained at $k\ge3$ regarding the global $y$-dependence.
Admitting the class $C^2$ results uniquely in
the ${\mathbf\cF\cF^{PD}_{g}}$-Finsleroid metric function $K$
given by the formulas (6.1)-(6.6).
This  function $K$
when considered on the $b$-slit tangent bundle
${\cal T}_bM~:= TM \setminus0\setminus b\setminus -b$
 is smooth of the class $C^{\infty}$
regarding the global $y$-dependence.

Owing to the theorems 2.1 and 2.2,
the ${\mathbf\cF\cF^{PD}_{g}}$-space
also occupies the  unique position when in the above chain
$(Z_1)$-$(Z_5)$  the  condition $(Z_3)$ is replaced by
the stipulation (2.1)-(2.3) with prescription of the dependence (2.15).

{%\pgbrk}

In  the ${\mathbf\cF\cF^{PD}_{g}}$-space,
 the method of Sections 3 and 4 proves to produce the explicit
and simple representations
for the respective connection coefficients and curvature tensor,
on assuming $h=const$ (which entails $g=const$ because of (6.5)).
The success has been predetermined by the possibility to write down the explicit coefficients
(6.26)
of the $\Ka$-transformation.
 The involved preferred vector field $b^i(x)$ proves to be  the proper element
of the  ${\mathbf\cF\cF^{PD}_{g}}$-space $\Ka$-transformation.
The  metricity (2.12)
holds fine.
The obtained formulas (6.30) and (6.60)-(6.62)  straightforwardly entail the vanishing (1.18).

The
${\mathbf\cF\cF^{PD}_{g}}$-space
connection coefficients (6.49) involve the fraction $1/q$,
where
$ q=\sqrt{r_{mn}y^my^n}$ with $r_{mn}=a_{mn}-b_mb_n.
$
Since
the input 1-form $b$  is of the unit norm $||b||=1$,
the scalar  $q$ is zero when $y=b$. Therefore, we may apply the coefficients on but
 the  $b$-slit  tangent bundle
$
{\cal T}_bM~:= TM \setminus0\setminus b\setminus -b
$
(obtained
by deleting out in $TM\setminus 0$
all the directions which point along, or oppose,
 the directions  given rise to by the  1-form $b$), on which
 the coefficients
are  smooth  of the class
  $C^{\infty}$ regarding the $y$-dependence.

As we are entitled to conclude from the right-hand part of the
representation (6.49) of the
${\mathbf\cF\cF^{PD}_{g}}$-space connection coefficients $D^k{}_{mn},$
the coefficients are not equal to the Riemannian Christoffel symbols
$a^k{}_{mn}$
of the space  ${\cal R}^N,$ unless we meet the vanishing $\nabla_nb_m=0,$
where
$\nabla$ stands for the Riemannian covariant derivative in  the space
 ${\cal R}^N$.
The last vanishing means geometrically that
the  vector field $b_m(x)$ is parallel in the space  ${\cal R}^N,$
in which case the coefficients
$D^k{}_{mn}$ are equivalent to
the coefficients
$a^k{}_{mn}.$
So, we are entitled to set forth the following assertion.

\ses

\ses

{\bf Theorem 2.6.}
{\it
When the involved vector field $b_m(x)$ is parallel  in the associated
Riemannian space ${\cal R}^N,$
the
${\mathbf\cF\cF^{PD}_{g}}$-space connection  obtained
is equivalent to the
 Levi-Civita connection in the
space ${\cal R}^N.$

}

\ses

\ses

Otherwise the connection coefficients
$D^k{}_{mn}$
are nonlinear (regarding the $y$-dependence) in any  dimension $N\ge3$.

In the dimension $N=2$ we always have  $\eta^k_m=0$
(see (6.50) and (6.51))
and, therefore,
the connection coefficients  $D^k{}_{mn}$ are independent of  $y$
(see (6.55)), which in turn entails the independence of the  tensor (6.75)
of $y$.
Thus we can formulate the following remarkable result.

\ses

\ses

{\bf Theorem 2.7.}
{\it
The
${\mathbf\cF\cF^{PD}_{g}}$-space connection  obtained
 is linear in the dimension $N=2$.
}

\ses

\ses

{%\pgbrk}

In the ${\mathbf\cF\cF^{PD}_{g}}$-space
of  the dimension
 $N=3$,
 the   arc
${\cal A}(x,l_1,l_2)$
can be described by means of dependence of the azimuthal angle
$\chi$ (see (6.4) and (6.14))
and the polar angle $\phi$
on
the arc-length parameter $s$ (defined by  (1.2)).
Due attentive consideration   performed in Section 7
 leads to the following assertion.

\ses

\ses

{\bf Theorem 2.8.}
{\it
In the ${\mathbf\cF\cF^{PD}_{g}}$-space
of  the dimension
 $N=3$,
 the  geodesic equation for the arc
${\cal A}(x,l_1,l_2)$
can  be completely integrated,
 yielding  the following explicit dependence:
\be
\chi(s)=\fr1h\arccos\lf(\sqrt{1-h^2\wt C^2}\cos\bigl(h(s-\wt s)\bigr)\rg)
\ee
and
}
\be
\phi(s)= \wt\phi  -\fr{\pi}2  +
\arctan\lf(\fr1{h\wt C}\tan\bigl(h(s-\wt s)\bigr)\rg),  ~~ \text{if} ~ \wt C\ne0;
\quad
\phi=\wt\phi,  ~~ \text{if} ~ \wt C = 0,
\ee
{\it
where
 $\wt C,\wt s,\wt\phi$ are  integration constants.
}

\ses

\ses

Using this dependence, we obtain the following theorem.

\ses

\ses

{\bf Theorem 2.9.}
{\it
The behaviour of
the unit vector $l^i$ along the  arc  ${\cal A}(x,l_1,l_2)$
of the
}
{\it
 $(N=3)$-dimensional space ${\mathbf\cF\cF^{PD}_{g}}$
is governed
by the  expansion
\be
l^i(s)=k_1(s)l^i_1+k_2(s)l^i_2 +k_3(s)b^i,
\ee
which in addition to
the pair $l^i_1,l^i_2$
 involves the Finsleroid-axis vector $b^i$.
}

\ses

\ses

The coefficients
$k_1(s),k_2(s),k_3(s)$
are given explicitly by means of the formulas
(7.43)-(7.46).

 {%\pgbrk}

\ses

\ses

\setcounter{sctn}{3} \setcounter{equation}{0}

\nin
  {\bf 3. Preliminary observations}

\ses

\ses

Let  the desired $\Ka$-transformation
(2.1)  of a Finsleroid space
 ${\cal F}^N$
 be realized over the tangent vectors by means of
a convenient diffeomorphic transformation
\be
y=\Ka\cdot\zeta: ~ ~ ~  y^i=y^i(x,\zeta).
\ee
Denote the inverse  by
\be
\zeta=\Ka^{-1}\cdot y: ~ ~  ~  \zeta^i=\zeta^i (x,y).
\ee
In (3.1), as well as in (3.2), it is implied that $y\in T_xM$ and $\zeta\in T_xM$
with the same point $x\in M$ of support.
The homogeneity  condition formulated below theorem 2.1 takes on the explicit form
\be
\zeta^i(x, \ga y) =
\ga^h\zeta^i(x,y), \qquad  \ga>0,\,\forall y,
\ee
which entails the identity
\be
y^n\zeta^i_n = h\zeta^i,
\ee
where
$\zeta^i_n =\partial{\zeta^i}/\partial{y^n}.
$
The transformation (2.1) can be written in the tensorial form
\be
g_{mn}(x,y)= \ka^2 \zeta^i_m(x,y) \zeta^j_n (x,y) a_{ij}(x).
\ee
From (3.4) and (3.5) it
just ensues that
the Finsleroid  metric function
 $K(x,y)=\sqrt{g_{mn}y^my^n}$ and
the Riemannian  metric function
$
S(x,\zeta)=\sqrt{a_{mn}(x)\zeta^m\zeta^n}
$
are connected by means of the relation
\be
K=h\varkappa    S.
\ee
Owing to
$\ka=(1/h)K^{1-h}$
(see (2.5)),
from  (3.6) we can
 obtain the remarkable equality
\be
\lf(K(x,y)\rg)^{h(x)}=S(x,\zeta).
\ee
The indicatrix  property (2.2) is a direct implication of the formulas (3.7)
and
\be
l=\Ka\cdot L: ~ ~ ~  l^i=y^i(x,L);
~ \qquad
L=\Ka^{-1}\cdot l: ~ ~  ~  L^i=\zeta^i (x,l)
\ee
(see (3.1) and (3.2)),
where
$l^i=y^i/K(x,y)$ and $L^i=\zeta^i/(S(x,\zeta)$ are the respective unit vectors
which possess the properties $K(x,l)=1$ and $S(x,L)=1$.
From (2.5) and (3.5) it follows that
\be
g_{mn}(x,l)dl^mdl^n=\fr1{h^2}a_{ij}(x)dL^i dL^j.
\ee
No support vector enters here the right-hand part.

{%\pgbrk}

Any two nonzero tangent vectors
  $y_1,y_2\in T_xM$ in  a fixed tangent space  $T_xM$
form the {\it Finsleroid-space angle}
  \be
 \al_{\{x\}}(y_1,y_2) =
 \fr1{h(x)}
 \arccos\la,
  \ee
 where the scalar
\be
\la =
  \fr{a_{mn}(x)\zeta_1^m  \zeta_2^n}
 { \sqrt{S^2(x, \zeta_1)}\,\sqrt{S^2(x, \zeta _2)} },
\qquad \text{with} \quad  \zeta_1^m = \zeta^m (x,y_1)  \quad   \text{and} \quad
 \zeta_2^m = \zeta^m (x,y_2)    ,
   \ee
is of the entire Riemannian  meaning in the
space ${\cal R}^N$.
These representations (3.10) and (3.11) realize the claimed relation (2.7),
with making the choice
$\sqrt{\cal C}=h$ in accordance with (2.6).

{%\pgbrk}

Let
the Finsleroid indicatrix ${\cal  IF}_{\{x\}}$
supported by   a fixed point $x\in M$
be parameterized by means of a convenient variable set $u^a$
(for instance,  we can take
$  u^a=\zeta^a/\zeta^N $     in regions  with $\zeta^N\ne0,$
or
$ u^a=y^a/y^N $, whenever  $y^N\ne0$.
The indices $a,b,c,d,e$ will be specified over
the range
$(1,...,N-1).$
Using the parametrical representation $l^i=l^i(u^a)$
of the indicatrix, where $l^i$ are unit vectors (possessing the property
$K(l)=1$),
 we can construct the induced metric tensor
\be
i_{ab}(u^c)=g_{mn}t^m_at^n_b
\ee
on the  indicatrix
by the help of the projection factors
$
t^m_a=\partial {l^m}/\partial {u^a}
$
(the method was described in detail in Section 5.8 of [1]).
Applying (3.5) yields the equality
\be
i_{ab}=\fr1{h^2}\wt i_{ab},
\ee
where
$
\wt i_{ab}(u^c)=a_{mn}\wt t^m_a \wt t^n_b \quad \text{with} \quad
\wt t^m_a= t^j_a\zeta^m_j
$
is the
Riemannian version of  the indicatrix induced metric tensor, obtainable when one puts $h=1$.
We have taken into account the fact that
the conformal factor $\ka$, having been proposed by (2.5), equals $1/h(x)$
 on the indicatrix.
From this standpoint, it is easy to make the search into the curvature of the indicatrix.
Indeed, since $h$ is independent of the vectors $y$,
 the associated Christoffel symbols
$$
i_a{}^c{}_b
=\fr12i^{ce}\lf(\D{i_{ea}}{u^b}+\D{i_{eb}}{u^a}-\D{i_{ab}}{u^e}\rg),
\qquad
 \wt i_a{}^c{}_b
=\fr12 \wt   i^{ce}\lf(\D{ \wt i_{ea}}{u^b}+\D{ \wt i_{eb}}{u^a}
-\D{ \wt i_{ab}}{u^e}\rg)
$$
are equivalent:
$
  i_a{}^c{}_b   =   \wt i_a{}^c{}_b .
$
Therefore,  the indicatrix curvature tensor
$$
I_a{}^c{}_{bd}=
\D{i_a{}^c{}_b}{u^d}-
\D{i_a{}^c{}_d}{u^b}
+i_a{}^e{}_b    i_e{}^c{}_d
-i_a{}^e{}_d    i_e{}^c{}_b
$$
\ses\\
 is identical to  the tensor
$\wt I_a{}^c{}_{bd}$
constructible by the same rule from the tensor
$ \wt i_{ab}$,
that is,
$I_a{}^c{}_{bd}
= \wt I_a{}^c{}_{bd}$.
Let us now consider the tensors
$
I_{acbd}  = i_{ce}I_a{}^e{}_{bd}
$ and
$\wt I_{acbd}
= \wt i_{ce}\wt I_a{}^e{}_{bd}.
$
Since
$
\wt I_{acbd}=  -
(\wt i_{cb}\wt i_{ad}-\wt i_{cd}\wt i_{ab})
$
is the ordinary case characteristic of  the  Riemannian  geometry
% upon  assuming $g=0$
(which  reflects the fact that  the    curvature of the unit sphere is equal to 1),
 we get
$
I_{acbd}=  -
 (\wt i_{cb}\wt i_{ad}-\wt i_{cd}\wt i_{ab})/h^2
$
and, then,  arrive at  the representation
\be
I_{acbd}=  -
  h^2
(i_{cb}i_{ad}-i_{cd}i_{ab})
\ee
which manifests that the indicatrix curvature is constant and equals
$h^2$
(in compliance with (2.6)).

{%\pgbrk}

Also,  from (3.13)
we have  $\det(i_{ab})=h^{-2(N-1)} \det(\wt i_{ab})$.
The area of the indicatrix is the volume $\int\sqrt{\det(i_{ab})}\,du^1...du^{N-1}$
of the internal indicatrix space,
 where integration is
performed over all the space.
Then we are entitled to conclude that
 in any dimension  $N $ and at each point $x\in M$
the ratio
\be
\fr{AREA_{\text{Finsleroid Indicatrix}}}
{AREA_{\text{Euclidean Unit Sphere}}}
=\fr1{h^{N-1}}
\ee
is valid.

The tensor $i_{ab}(u^c)$ is defined on the indicatrix.
We can, however, extend the meaning of the parameters $u^a$ by
homothety  to obtain the scalars $u^a(y)$ defined at any point of the Finsleroid,
using
the zero-degree homogeneity $u^a(ky)=u^a(y), ~ k>0,  \forall y. $
With their help we can obtain the tensor
 $i^*_{ab}(y)=i_{ab}(u(y))$ meaningful at any point of the Finsleroid
and construct  the extended tensor
$f_{AB}    = f_{AB} (u^0,u^a)$ with    $u^0=\ln K$ as follows:
\be
f_{AB} =\e^{2u^0}f^*_{AB},
\quad
\text{with}
\quad
f^*_{ab} =i^*_{ab}, \quad f^*_{a0}=0,   \quad f^*_{00}  = 1,
\ee
where the sets  $u^A=\{u^0,u^a\}$ and
$u^B=\{u^0,u^b\}$ have been used.
In terms of
the parametrization $y^m=y^m(u^A)$
thus arisen, it can readily be seen that
the tensor
$f_{AB}$
is the covariant transform of the Finslerian metric tensor
$g_{mn}$:
\be
f_{AB}=g_{mn}\D{y^m}{u^A}  \D{y^n}{u^B}.
\ee
So,
the volume $\int\sqrt{\det(g_{mn})}\,d^Ny$
of the Finsleroid when written in terms of the coordinates $u^A$
is given by the integral
$$
\int\sqrt{\det(f_{AB})}\,du^1...du^{N-1} du^0
=\int \e^{Nu^0} du^0\int\sqrt{\det(f^*_{ab})}\,du^1...du^{N-1},
$$
in which $u^0\in (-\infty,0)$,
so that the integral  is equal to the $1/N$ multiplied by the area of the indicatrix.
Whence in addition to the law (3.15) we have the ratio
\be
\fr{VOLUME_{\text{Finsleroid }}}
{VOLUME_{\text{Euclidean Unit Ball}}}
=\fr1{h^{N-1}}.
\ee

\ses

{%\pgbrk}

Therefore,
 at each point $x\in M$
the following law is valid:
\be
\fr{AREA_{\text{Finsleroid Indicatrix}}}
{VOLUME_{\text{Finsleroid }}}
=
                          \fr{AREA_{\text{Euclidean Unit Sphere}}}
{VOLUME_{\text{Euclidean Unit Ball }}},
\qquad \text{when} \quad N\ge3,
 %{%%%%%=\text{not the Universal Value}}
\ee
\ses
and
\be
\fr{LENGTH_{\text{Finsleroid Indicatrix}}}
{AREA_{\text{Finsleroid }}}
=
                          \fr{LENGTH_{\text{Euclidean Unit Circle}}}
{AREA_{\text{Euclidean Unit  Circle }}}
=\fr{\fr{2\pi}h}
{\fr{\pi}h}=2,
\qquad \text{when} \quad N =2.
\ee

{%\pgbrk}

With the help of the  derivative coefficients
\be
\zeta^i_n =\D{\zeta^i}{y^n}, \quad \zeta^i_{nk} =\D{\zeta^i_n}{y^k},
  \qquad y^i_m=\D{y^i}{\zeta^m},  \quad  y^i_{mh}=\D{y^i_m}{\zeta^h},
\ee
it is possible
 to develop a direct method to induce the
connection   in the Finsleroid
 space    ${\cal F}^N$ from the Riemannian space ${\cal R}^N$.
To this end we can naturally use in
${\cal R}^N$
the Levi-Civita connection
 \be
{\cal LC}=\{L^m_j,L^m{}_{ij}\}: \qquad    L^m_{j}=-L^m{}_{ij}\zeta^i,
   \quad
   L^m{}_{ij}=a^m{}_{ij},
\ee
 with $a^m{}_{ij}=a^m{}_{ij}(x)$ standing for the Christoffel symbols
constructed from the Riemannian metric tensor $a_{mn}(x)$.

First of all, we need the coefficients $N^k_i(x,y)$ to construct the
operator $d_i$ indicated in (2.9).
It proves fruitful to obtain the coefficients by means of the map
\be
\{N^k_i\}=\Ka\cdot \{L^k_i\}.
\ee
Namely,
starting with the fundamental property of
 the Levi-Civita connection
that the Riemannian
angle is preserving under the parallel displacements, which in terms of our notation
can be written as
\be
\Biggl(\fr{\partial}{\partial x^i}
+L_i^k{(x,\zeta_1)} \fr{\partial}{\partial \zeta_1^k}
+L_i^k{(x,\zeta_2)} \fr{\partial}{\partial \zeta_2^k}\Biggr)
{\al_{\text{Riemannian space}}(x,\zeta_1,\zeta_2)}
=0,
\ee
we want to have  the similar vanishing  in the  Finsleroid  space ${\cal F}^N$:
\be
\Biggl(\fr{\partial}{\partial x^i}
+N^k_i{(x,y_1)} \fr{\partial}{\partial y_1^k}
+N_i^k{(x,y_2)} \fr{\partial}{\partial y_2^k}\Biggr)
\al_{\{x\}}(y_1,y_2)
=0,
\ee
assuming also
that the vanishing
\be
\Biggl(\fr{\partial}{\partial x^i}
+N^k_i{(x,y)} \fr{\partial}{\partial y^k}
\Biggr) K(x,y)=0
\ee
arises after performing the $\Ka$-transformation of the Riemannian vanishing
\be
\Biggl(\fr{\partial}{\partial x^i}
+L^k_i{(x,\zeta)} \fr{\partial}{\partial \zeta^k}
\Biggr) S(x,\zeta)=0.
\ee

{%\pgbrk}

With
an arbitrary  differentiable scalar
$w(x,y)$,
we consider the $\Ka$-transform
\be
W(x,\zeta)=w(x,y), \quad \text{ which entails} \quad \D W{\zeta^n}= y^k_n \D w{y^k},
\ee
and  postulate  that the $\Ka$-transformation is {\it covariantly transitive},
so that
\be
\Biggl(\fr{\partial}{\partial x^i}
+N_i^k{(x,y)} \fr{\partial}  {\partial y^k}
\Biggr)
w(x,y)
=
\Biggl(\fr{\partial}{\partial x^i}
+L_i^k{(x,\zeta)} \fr{\partial}{\partial \zeta^k}
\Biggr)
W(x,\zeta).
\ee
Since the field $w$ is arbitrary, the last equality is fulfilled if and only if
\be
N^m_n=\D{y^m(x,\zeta)}{x^n}  +   y^m_iL^i_n.
\ee
This is the    representation  which is required  to
realize the map (3.23).

Since the equality (3.30)  can be written in the form
\be
\D{\zeta^i}{x^n}+N^m_n\zeta^i_m+a^i{}_{kn}\zeta^k=0,
\ee
we have
\be
N^m_n=-y^m_i \lf(\D{\zeta^i}{x^n}+a^i{}_{kn}\zeta^k\rg).
\ee

\ses

It can readily be noted that the transitivity property (3.29) can
straightforwardly
be extended to scalars
dependent on two vectors. Namely,
if
\be
W(x,\zeta_1,\zeta_2)=w(x,y_1,y_2),
\ee
then
\be
\Biggl(\fr{\partial}{\partial x^i}
+N_{1i}^k \fr{\partial}  {\partial y_1^k}
+N_{2i}^k \fr{\partial}  {\partial y_2^k}
\Biggr)
w(x,y_1,y_2)
=
\Biggl(\fr{\partial}{\partial x^i}
+ L_{1i}^k \fr{\partial}{\partial \zeta_1^k}
+ L_{2i}^k \fr{\partial}{\partial \zeta_2^k}
\Biggr)
W(x,\zeta_1,\zeta_2),
\ee
where $N_{1i}^k=N_i^k(x,y_1),\,N_{2i}^k=N_i^k(x,y_2),\,
L_{1i}^k =L_i^k{(x,\zeta_1)}, \,L_{2i}^k =L_i^k{(x,\zeta_2)} $.
The  equality (3.34) is verified by  using (3.31)  or (3.32).
When this implication is applied to the equality (3.10),
the Riemannian vanishing (3.24) just entails the Finsleroid-space vanishing (3.25),
whenever $h=const$.

{%\pgbrk}

Differentiating (3.5) with respect to $y^k$ yields the representation

\be
2C_{mnk}=
(1-h)\fr2{K}l_k
g_{mn}
+\ka^2 (\zeta^i_{mk}\zeta^j_n +\zeta^i_m\zeta^j_{nk}) a_{ij}
\ee
for the Cartan tensor.
Contracting this  by $y^n$ results in the equality

\be
\ka^2 h\zeta^i_{mk}\zeta^j  a_{ij}=
(1-h)(h_{km} -l_k  l_m ),
\ee
where the vanishing $C_{mnk}y^n=0$ and the homogeneity identity (3.4) have
 been taken into account.

{%\pgbrk}

From (3.11) it follows that
$$
\D{\la}{x^i}=
  \fr{a_{mn,i}\zeta_{1}^m  \zeta_2^n}
{ \sqrt{S^2(x, \zeta_1)}\,\sqrt{S^2(x, \zeta _2)} }+
\fr1{{ \sqrt{S^2(x, \zeta_1)}\,\sqrt{S^2(x, \zeta _2)} }}
  a_{mn}\lf(\D{\zeta_{1}^m}{x^i}  \zeta_2^n
  +
\zeta_{1}^m  \D{\zeta_2^n  }{x^i}
\rg)
$$

\ses

\be
-
\fr12\la
\lf[
\fr1{ S^2(x, \zeta_1) }
\lf(a_{mn,i}\zeta_{1}^m\zeta_1^n
+ 2 a_{mn} \D{\zeta_{1}^m}{x^i}  \zeta_1^n
 \rg)
    +
\fr1{ S^2(x, \zeta_2) }
\lf(a_{mn,i}\zeta_{2}^m\zeta_2^n
+ 2 a_{mn} \D{\zeta_{2}^m}{x^i}  \zeta_2^n
 \rg)
\rg],
\ee
\ses
where   $a_{mn,i}=\partial a_{mn}/\partial x^i$,
and
\be
\D{\la}{y^k_1}=     \fr{a_{mn}\zeta_{1k}^m  \zeta_2^n}
{ \sqrt{S^2(x, \zeta_1)}\,\sqrt{S^2(x, \zeta _2)} }
 -      \fr{a_{mn}\zeta_{1k}^m  \zeta_1^n}
{ S^2(x, \zeta_1)}\la,
\qquad
\D{\la}{y^k_2}=     \fr{a_{mn}\zeta_{2k}^m  \zeta_1^n}
{ \sqrt{S^2(x, \zeta_2)}\,\sqrt{S^2(x, \zeta _1)} }
 -      \fr{a_{mn}\zeta_{2k}^m  \zeta_2^n}
{ S^2(x, \zeta_2)}\la.
\ee
With  (3.32)
we find
\be
N^k_{1i}
\D{\la}{y^k_1}=
- \lf(\D{\zeta^m_1}{x^i}+a^m{}_{ti}\zeta^t_1\rg)
\lf[
  \fr{a_{mn}  \zeta_2^n}
{ \sqrt{S^2(x, \zeta_1)}\,\sqrt{S^2(x, \zeta _2)} }  -
   \fr{a_{mn}  \zeta_1^n}
{ S^2(x, \zeta_1) }\la
\rg]
\ee
and
\be
N^k_{2i}
\D{\la}{y^k_2}=
- \lf(\D{\zeta^m_2}{x^i}+a^m{}_{ti}\zeta^t_2\rg)
\lf[
  \fr{a_{mn}  \zeta_1^n}
{ \sqrt{S^2(x, \zeta_1)}\,\sqrt{S^2(x, \zeta _2)} }
 -
   \fr{a_{mn}  \zeta_2^n}
{ S^2(x, \zeta_2) } \la
\rg],
\ee
where the identity
$y^k_j\zeta^m_k=\de^m_j$
has been taken into account.

{%\pgbrk}

The formulas (3.37), (3.39), and (3.40) just entail the vanishing
\be
\D{\la}{x^i}
+N^k_{1i}
\D{\la}{y^k_1}
+
N^k_{2i}
\D{\la}{y^k_2}
=0.
\ee
Thus, from (3.10) we may conclude that whenever $h=const$ the angle preservation
(1.18) holds fine.

{%\pgbrk}

\bigskip

\setcounter{sctn}{4} \setcounter{equation}{0}

\nin
  {\bf 4. Entailed connection coefficients and curvature tensor}

\bigskip

Let us trace the validity of the theorem 2.5
and the involved formulas (2.11)-(2.13).
Since   $\cD_iK=d_iK$,  the  vanishing $d_iK=0$ indicated in (2.11) ensues from  (3.7)
and (3.27).
The second vanishing   in (2.11) is tantamount to
$
N^j_i=-D^j{}_{ik}y^k
$
(because of  $\partial y^j/\partial x^i=0$). The third equality entered  (2.11) reads
\be
\D{y_j}{x^i}+N_i^kg_{kj} - D^h{}_{ij}y_h=0.
\ee
Let us differentiate this equality with respect to $y^n$. We obtain
\be
d_ig_{jn}+\D{N_i^k}{y^n}g_{jk} - D^h{}_{ij}g_{hn}
-y_h\D{ D^h{}_{ij}}{y^n}
=0.
\ee
By making  the choice
$
D^k{}_{in} =-\partial {N^k_i}/\partial{y^n}
$
we obtain from (4.2)    the  metricity
$
\cD_ig_{jn}=0,
$
if
\be
y_h\D{ D^h{}_{ij}}{y^n}=0.
\ee

From  (3.30) and
$
D^k{}_{in} =-\partial {N^k_i}/\partial{y^n}
$
  it follows that
\be
\D{ y^n_k}{x^i}+L^t_i y^n_{kt}+D^n{}_{is}y^s_k-L^h{}_{ik}y^n_h=0.
\ee
Since
$ y^n_k\zeta^k_j=\de^n_j$,
the previous identity can be written as
\be
\D{\zeta_m^s}{x^i}+N^h_i\zeta^s_{mh}+L^s{}_{it}\zeta^t_m-D^h{}_{im}\zeta^s_h=0.
\ee

\ses

Can the last vanishing be materialized?

\ses

{%\pgbrk}

Let us realize the action of the $\Ka$-transformation (3.1)-(3.2)
on tensors  by the help of the transitivity  rule,
that is,
\be
\{w^n{}_m(x,y)\}  =\Ka\cdot \{W^n{}_m(x,\zeta)\}: ~ ~ ~  w^n{}_m=  y^n_h\zeta^j_mW^h{}_j,
\ee
and define the
   covariant derivative  $\nabla$ in ${\cal R}^N$
according to the  conventional  Riemannian rule:
\be
\nabla_iW^n{}_m=\D{W^n{}_m}{x^i}+L^k_i   \D{W^n{}_m}{\zeta^k}
+
L^n{}_{hi}W^h{}_m
-L^h{}_{mi}W^n{}_h
\ee
and
\be
\nabla_iS=0, \qquad \nabla_i\zeta^j=0, \qquad  \nabla_ia_{mn}=0.
\ee
Due to   (4.4) and (4.5), we have the transitivity property
\be
\cD_iw^n{}_m=y^n_h\zeta^j_m\nabla_iW^h{}_j.
\ee
Applying  the rule (4.9) to the transformation (3.5) of the metric tensor yields
\be
g_{mn}d_i\lf(\fr1{\ka^2}\rg)+
\fr1{\ka^2}
\cD_ig_{mn}=0.
\ee
Thus,
{\it  the metricity condition }
 $
\cD_ig_{jn}=0
$
  {\it       holds if and only if   }
$d_i\ka=0$.

{%\pgbrk}

Applying
 (2.5) and   the vanishing  $\cD_i K=0$  to (4.10)
makes us conclude that  the following assertion is valid.

\ses

\ses

{\bf Theorem 4.1.}
{\it
Under the  input stipulation} (2.1)-(2.3), {\it the covariant derivative $\cD_i$
obtained through the transitivity} (4.9) {\it fulfills
  the metricity condition
}
 $   \cD_ig_{jn}=0
 $
{\it   if and only if
}
\be
\D{h}{x^i}=0: ~~~ h=\text{const}.
\ee

\ses   \ses

The last condition entails the vanishing  (4.3);
in the ${\mathbf\cF\cF^{PD}_{g}}$-Finsleroid space
the validity of this implication can  explicitly be verified
with the help of the representation (6.55) derived in Section 6.

\ses

Commuting the covariant derivative (2.10) yields the equality
\be
\lf[D_iD_j-D_jD_i\rg] w^n{}_k=M^h{}_{ij}\D{w^n{}_k}{y^h}  -E_k{}^h{}_{ij}w^n{}_h
+E_h{}^n{}_{ij}w^h{}_k
\ee
with the  tensors
\be
 M^n{}_{ij}~:= d_iN^n_j-d_jN^n_i
\ee
and
\be
E_k{}^n{}_{ij}~: =
d_iD^n{}_{jk}-d_jD^n{}_{ik}+D^m{}_{jk}D^n{}_{im}-D^m{}_{ik}D^n{}_{jm}.
\ee
\ses

If the choice
$ D^k{}_{in} =-N^k{}_{in} $
is made (see (2.13)),
 the tensor (4.13) can be written in the form
%Owing to the equality $N^j_i=-D^j{}_{ik}y^k$,
% we can write the tensor (4.13) in the form
\be
 M^n{}_{ij}= \D{N^n_j}{x^i}-\D{N^n_i}{x^j}
-N_i^hD^n{}_{jh}+N_j^hD^n{}_{ih}.
\ee
By applying the commutation rule (4.12)
to the particular choices $\{K, y^n, y_k,g_{nk}\}$,
we obtain the identities
\be
y_n M^n{}_{ij}=0,
\qquad
y^kE_k{}^n{}_{ij}  = -   M^n{}_{ij},   \qquad
y_nE_k{}^n{}_{ij}  = M_{kij},
\ee
and
\be
E_{mnij}+E_{nmij}=2C_{mnh}
M^h{}_{ij} \quad \text{with} \quad C_{mnh}=\fr12\D{g_{mn}}{y^h}.
\ee
Differentiating (4.15) with respect to $y^k$ and using the equality
$
N^j_i=-D^j{}_{ik}y^k
$
(see (2.13))
 yield
\be
E_k{}^n{}_{ij}= -\D{M^n{}_{ij}}{y^k}.
\ee

The cyclic identity
\be
D_kM^n{}_{ij}+D_jM^n{}_{ki}+D_iM^n{}_{jk}=0
\ee
is valid,
where
\be
D_k M^n{}_{ij}=\D{ M^n{}_{ij}}{x^k}+N^m_k\D{ M^n{}_{ij}}{y^m}+D^n{}_{kt} M^t{}_{ij}
 - a^s{}_{ki}   M^n{}_{sj}  -  a^s{}_{kj}   M^n{}_{is}.
\ee

{%\pgbrk}

It proves pertinent to replace in the commutator (4.12)
the partial derivative
$\partial  w^n{}_k/\partial y^h$
by the definition
\be
\cS_h{w^n{}_k}=\D{w^n{}_k}{y^h}+C^n{}_{hk}w^h{}_k-C^m{}_{hk}w^n{}_m
\ee
which has the meaning of the covariant derivative in the tangent
Riemannian space
$\cR_{\{x\}}$.
With
the {\it curvature tensor}
\be
\Rho_k{}^n{}_{ij}=E_k{}^n{}_{ij}
-
M^h{}_{ij}C^n{}_{hk},
\ee
the commutator (4.12) takes on the form
\be
\lf(\cD_i\cD_j-\cD_j\cD_i\rg) w^n{}_k=
M^h{}_{ij}\cS_h w^n{}_k  -\Rho_k{}^h{}_{ij}w^n{}_h
+\Rho_h{}^n{}_{ij}w^h{}_k.
\ee
The skew-symmetry
\be
 \Rho_{mnij}=-\Rho_{nmij}
\ee
 holds (cf.  (4.17)).

{%\pgbrk}

\setcounter{sctn}{5} \setcounter{equation}{0}

\nin
  {\bf 5. Specifying the conformal multiplier}

\ses

\ses

In a fixed tangent space endowed with
 a Finslerian metric tensor $g_{ij}$
 produced by
 a Finslerian metric function $F$,
  we may consider the conformal transform
\be
k_{ij}={\e^{2{\psi}}}{g_{ij}}, \qquad {\psi}={\psi}(F),
\ee
 where
${\psi}= {\psi}(F)   $ is a smooth function,
and construct from $k_{ij}$ the
Christoffel symbols $k^n{}_{ij}$, which yields
$
k^n{}_{ij}=
{\psi}'l_i\de^n_j  +  {\psi}'l_j\de^n_i  -  {\psi}'l^ng_{ij}
+C^n{}_{ij},
$
where the prime means differentiation with respect to $F$
and
$C_{nij}=(1/2)\partial g_{ij}/\partial y^n$.
We  directly derive the equalities
$$
\D{k^n{}_{ij}}{y^m}-  \D{k^n{}_{im}}{y^j}  =
\lf({\psi}''l_il_m+  {\psi}'\fr1Fh_{im}\rg)
\de^n_j
 - \lf({\psi}''l^nl_m+{\psi}'\fr1Fh^n_m
\rg)  g_{ij}
-2 C^n{}_{lm}C^l{}_{ij}
-[jm]
$$
and
$$
k^t{}_{ij}  k^n{}_{tm}   -   k^t{}_{im}  k^n{}_{tj}
=
{\psi}'l_j
{\psi}'l_i\de^n_m
-  {\psi}'g_{ij}
{\psi}'h^n_m
+C^n{}_{tm}C^t{}_{ij}
-[jm],
$$
together with
$$
\D{k^n{}_{ij}}{y^m}-  \D{k^n{}_{im}}{y^j} + k^t{}_{ij}  k^n{}_{tm}   -   k^t{}_{im}  k^n{}_{tj}
+C^n{}_{tm}C^t{}_{ij}  -C^n{}_{tj}C^t{}_{im}
$$

\ses

$$
=
-\lf({\psi}''l_il_j+  {\psi}'\fr1Fh_{ij}\rg)
\de^n_m
 -\lf( {\psi}''l^nl_m+{\psi}'\fr1Fh^n_m\rg)  g_{ij}
+
{\psi}'l_j
{\psi}'l_i\de^n_m
-  {\psi}'g_{ij}{\psi}'h^n_m
-[jm].
$$
The notation $[jm]$ symbolizes the skew-symmetric terms.
We introduce the associated tensor
\be
L_i{}^n{}_{mj}=
\D{k^n{}_{ij}}{y^m}-  \D{k^n{}_{im}}{y^j} + k^t{}_{ij}  k^n{}_{tm}   -   k^t{}_{im}  k^n{}_{tj}
\ee
\ses
and   the indicatrix curvature tensor
\be
\hat R_i{}^n{}_{mj}=
C^n{}_{tm}C^t{}_{ij}  -C^n{}_{tj}C^t{}_{im},
\ee
\ses
obtaining
$$
L_i{}^n{}_{mj}
+
\hat R_i{}^n{}_{mj}
=
-({\psi}''-{\psi}'{\psi}')l_i(l_j h^n_m  -l_m h^n_j)
-
{\psi}''l^n(l_m h_{ij}-l_j h_{im})
$$

\ses

$$
+  {\psi}'\fr1F(h_{im}   \de^n_j-h_{ij}   \de^n_m)
+{\psi}'\fr1F(h^n_j g_{im}-h^n_m g_{ij})
+  {\psi}'{\psi}'(g_{im} h^n_j-g_{ij} h^n_m)
$$
\ses
and
$$
\e^{-2\psi} L_{inmj}   +   \hat R_{inmj}=
-({\psi}''-{\psi}'{\psi}')l_i(l_j h_{nm}  -l_m h_{nj})
-
{\psi}''l_n(l_m h_{ij}-l_j h_{im})
$$

\ses

$$
+  {\psi}'\fr1F(h_{im}   g_{nj}-h_{ij}  g_{nm})
+{\psi}'\fr1F(h_{nj} g_{im}-h_{nm} g_{ij})
+  {\psi}'{\psi}'(g_{im} h_{nj}-g_{ij} h_{nm}),
$$
\ses
where
$  L_{inmj} =k_{nk}L_i{}^k{}_{mj}$
and
$  \hat R_{inmj}  = g_{nk}\hat R_i{}^k{}_{mj}$.
Using the equality $h_{ij}=g_{ij}-l_il_j$ yields

{%\pgbrk}

$$
\e^{-2\psi} L_{inmj}   +   \hat R_{inmj}=
-\lf({\psi}''+{\psi}'\fr1F\rg)
\Bigl[l_i(l_j h_{nm}  -l_m h_{nj})    -      l_n(l_j h_{im}  -l_m h_{ij})\Bigr]
$$
\ses
\be
- \lf( {\psi}'\fr2F+  {\psi}'{\psi}'\rg)
(h_{ij} h_{nm}-h_{im} h_{nj}).
\ee

If  the curvature tensor $L_i{}^n{}_{mj}$
obtained  under the conformal transformation (5.1)
vanishes identically,
then
because of the known Finslerian
identities  $y^i C_{ijk}=0$ and  $y^i h_{ij}=0$
the equality (5.4) would entail the equation
$
{\psi}''F   +  {\psi}'   =0,
$
which solution is
\be
\e^{\psi}=c_1F^{c_2},
\ee
where $c_1>0$ and $c_2$ are integration constants.
The result (5.5)
permits writing  the conformal multiplier (2.3)
in the form  $\ka=(1/c_1)F^{1-h}$, where we have identified $c_2$ with $h-1$.
This entails that
the tensor $\zeta^i_m(x,y) \zeta^j_n (x,y) a_{ij}(x)$ appeared in the right-hand part of
the Finslerian metric tensor representation (3.5)
is positively homogeneous of degree $2h-2$. The last observation is in agreement with the
homogeneity condition (3.3), whence we have
 $c_1=h$.
Therefore,
the  theorem 2.1 of Section 2 is valid.

If we put
$L_i{}^n{}_{mj}=0$ and insert (2.5) into the right-hand part of (5.4),
 we obtain  for the indicatrix curvature tensor
(5.3)
the representation
\be
F^2\hat R_{inmj}
=
(1-h^2)(h_{ij} h_{nm}  - h_{im} h_{nj} )
\ee
which says us that the indicatrix is of  the constant curvature
 $h^2$.
Thus,
the  theorem 2.2 of Section 2 is fulfilled.

{%\pgbrk}

Inversely,
let  the indicatrix  be of  constant curvature, so that
\be
F^2\hat R_{inmj}
=
(1-{\cal C})(h_{ij} h_{nm}  - h_{im} h_{nj} ),
\ee
where ${\cal C}>0$ is the curvature value (and ${\cal C}$
 is independent of the tangent vectors $y$).
Inserting (5.7) in (5.4), performing the conformal transformation (5.1),
and applying  (5.5) with the choice of the exponent $c_2$ according to the condition
\be
(1+c_2)^2={\cal C},
\ee
from (5.4) we obtain the vanishing
$L_i{}^n{}_{mj}=0$, which means that the space is conformally flat.
Therefore,
the  theorem  2.4 of Section 2 is also  valid.

{%\pgbrk}

\ses

\ses

\setcounter{sctn}{6} \setcounter{equation}{0}

\nin
  {\bf 6. Performing  the choice of the ${\mathbf\cF\cF^{PD}_{g}}$-Finsleroid space}

\ses

\ses

Let us  assume that in addition to a
Riemannian  metric
 $\sqrt{a_{ij}(x)y^iy^j}$
the manifold $M$ admits a non-vanishing 1-form
$ b=b_i(x)y^i$ of the unit length:
\be
a_{ij}(x)b^i(x)b^j(x)=1,
\ee
where
$b^i(x)=a^{ij}(x)b_j(x).$
The tensor $a^{ij}(x)$ is reciprocal to $a_{ij}(x)$, so that
$a_{ij}a^{jn}=\de^n_i$, where $\de^n_i$ stands for the Kronecker symbol.
The Finsleroid  space is  specified in accordance with  the condition
that the   metric function
$K(x,y)$ is (2.14) with
\be
B=b^2+gbq+q^2\equiv A^2+h^2q^2 \quad {\text{with}} \quad A=b+\fr12gq,
\ee
where
\be
q=\sqrt{r_{mn}y^my^n} \quad \text{and} \quad     r_{mn}=a_{mn}-b_mb_n,
\ee
so that
\be
 a_{ij}(x)y^iy^j=b^2+q^2.
 \ee

The scalar $g(x)$ obtained through
\be
h(x)=\sqrt{1- \fr{g^2(x)}4}, \qquad \text{with} \quad -2<g(x)<2,
\ee
plays the role of the characteristic parameter.
The variable $\chi$ entering  the exponential representation
(2.14) of the ${\mathbf\cF\cF^{PD}_{g}}$-Finsleroid metric function
$K$
is given as it follows:
\be
\chi=\fr1h \Bigl(
-\arctan   \fr G2   +\arctan\fr{L}{hb}\Bigr),    ~  {\rm if}  ~ b\ge 0;
\quad
\chi=\fr1h \Bigl(
 \pi-\arctan
\fr G2
+\arctan\fr{L}{hb}\Bigr),
~   {\rm if}
~ b\le 0,
\ee
 with the function
$    L =q+(g/2) b $
fulfilling   the identity
\be
 L^2+h^2b^2=B.
 \ee
The definition range
$$
0\le\chi\le\fr1h\pi
$$
is of value to describe all the tangent space.
The normalization in (6.6)
is such that
\be
\chi\bigl|_{y=b}\bigr. =0.
\ee
The quantity (6.6) can conveniently be written as
\be
\chi  =  \fr1h  f
\ee
with
the function
\be
f=\arccos \fr{ A(x,y)}   {\sqrt{B(x,y)}}
\ee
ranging as follows:
\be
0\le f\le \pi.
\ee

{%\pgbrk}

\nin
The Finsleroid-axis vector $b^i$ relates to the value $f=0$, and
the opposed vector $-b^i$ relates to the value $f=\pi$:
\be
f=0 ~~  \sim ~~ y=b;  \qquad f=\pi ~~ \sim ~~ y=-b.
\ee
            It is  frequently
            convenient to  represent the function $K$ in the form
\be
K=\sqrt B\, J, \qquad \text{with} ~~ J=\e^{-\frac12 g \chi}.
\ee
The normalization is such that
\be
K(x,b(x))=1
\ee
(notice that
 $q=0$ at $y^i=b^i$).
The positive  (not absolute) homogeneity  holds:
 $K(x, \ga y)=\ga K(x,y)$ for any $\ga >0$   and all admissible $(x,y)$.

 Under these conditions, we call $K(x,y)$
the {\it  ${\mathbf\cF\cF^{PD}_{g}}$-Finsleroid   metric function},
obtaining  the ${\mathbf\cF\cF^{PD}_{g}}$-{\it Finsler space}
\be
{\mathbf\cF\cF^{PD}_{g}} :=\{M;\,a_{ij}(x);\,b_i(x);\,g(x);\,K(x,y)\}.
\ee

\ses

\ses

 {\large  Definition}.  Within  any tangent space $T_xM$, the  metric function $K(x,y)$
  produces the {\it    ${\mathbf\cF\cF^{PD}_{g} } $-Finsleroid}
 \be
 \cF\cF^{PD}_{g;\,\{x\}}:=\{y\in   \cF\cF^{PD}_{g; \, \{x\}}: y\in T_xM , K(x,y)\le 1\}.
  \ee

\ses

 \ses

 {\large  Definition}. The {\it    ${\mathbf\cF\cF^{PD}_{g} } $-Indicatrix}
 $ {\cal I}\cF^{PD}_{g; \, \{x\}} \subset T_xM$ is the boundary of the
    ${\mathbf\cF\cF^{PD}_{g} } $-Finsleroid, that is,
 \be
{\cal I}\cF^{PD}_{g\, \{x\}} :=\{y\in {\cal I}\cF^{PD}_{g\, \{x\}} : y\in T_xM, K(x,y)=1\}.
  \ee

\ses

 \ses

 {\large  Definition}. The scalar $g(x)$ is called
the {\it Finsleroid charge}.
The 1-form $b=b_i(x)y^i$ is called the  {\it Finsleroid--axis}  1-{\it form}.

\ses

\ses

The entailed components $y_i :=(1/2)\partial {K^2}/ \partial{y^i}$)
of the  covariant tangent vector $\hat y=\{y_i\}$
can be found in the simple form
\be
y_i=(u_i+gqb_i) J^2,
\ee
where $u_i=a_{ij}y^j$.

{%\pgbrk}

By making due inspection into the formulas (6.6)-(6.12)
it proves convenient
to  separate the tangent space
 $ T_xM$
into  the unification
\be
 T_xM=T^{\{+\}}_{g;\,\{x\}}  \cup   T^{\{-\}}_{g;\,\{x\}}
\ee
of the  regions
\be
T^{\{+\}}_{g;\,\{x\}} ~:=
\{y\in T^{\{+\}}_{g;\,\{x\}} ~: ~ ~  y\in T_xM , ~ ~ 0\le f<h\pi\}
\ee
\ses
and
\be
T^{\{-\}}_{g;\,\{x\}}M ~:=
\{y\in T^{\{-\}}_{g;\,\{x\}}~:~  ~ y\in T_xM, ~ ~  (1-h)\pi< f \le\pi\},
\ee
which depend on value of $g$.
We have
the range correspondence
\be
0\le f<h\pi \quad\sim  \quad \chi\in \lf[0,\pi\rg),
\qquad \text{and} \qquad
(1-h)\pi< f \le\pi  \quad \sim \quad  \chi\in \lf(\fr{1-h}h\pi,\fr1h\pi\rg].
\ee
The directions involving the Finsleroid axis vector $y=b$ belong to the region
(6.20), and the opposed cases belong to the region (6.21), that is,
\be
b(x)\in T^{\{+\}}_{g;\,\{x\}}
 \quad  \text{and} \quad
 -b(x)\in T^{\{-\}}_{g;\,\{x\}}.
\ee
The intersections
\be
  C^{\{+\}}_{g;\,\{x\}}=
  T^{\{+\}}_{g;\,\{x\}}  \cap
 {\cal I}\cF^{PD}_{g; \, \{x\}}
 \quad  \text{and} \quad
C^{\{-\}}_{g;\,\{x\}}=T^{\{-\}}_{g;\,\{x\}}  \cap
 {\cal I}\cF^{PD}_{g; \, \{x\}}
\ee
yield  two {\it covering charts}
 for the indicatrix
(6.17).

It will be noted that, in contrast to
$b^i$,
the opposed vector $-b^i$ is not unit. Therefore, we introduce the normalized vector
\be
b^{\{-\}i}=-b^i\e^{\frac12 g \pi}
\ee
which is unit: it can readily be seen that
$$
K\lf(x,b^{\{-\}}(x)\rg)=1.
$$

{%\pgbrk}

In this space the $\Ka$-transformation (3.2) can be realized in the explicit and simple form
\be
\zeta^i=\lf[hv^i  +   ( b+\fr12gq)b^i    \rg] \fr J{\varkappa h},
\ee
where $v^i=y^i-bb^i$ and  $\ka=(1/h)K^{1-h}$.
Obviously, the right-hand part in (6.26)
  possesses the  homogeneity properties (3.3)-(3.4)
of degree
 $h$.
Simple calculation shows that
\be
\det(\zeta^i_m)   =\lf(\fr J{\varkappa}\rg)^N .
\ee

Since $v^i=q=0$ at $y^i=\pm b^i$, from (6.26) it follows that
\be
\zeta^i(x,b)=b^i(x),  \qquad
\zeta^i(x,b^{\{-\}})=b^{\{-\}i}(x),
\ee
where the equality $K(x,b)=1$ (see (6.14)) has been taken into account.
Therefore,
{\it  the involved preferred vector field $b^i(x)$ as well as the opposed field
$b^{\{-\}i}(x)$
are the proper elements
of the
  ${\mathbf\cF\cF^{PD}_{g}}$-space
  $\Ka$-transformation
  }
  (6.26),
  {\it  that is,}
\be
\lf\{b^i(x)\rg\}\stackrel{\Ka }{\to} \lf\{b^i(x)\rg\},
\qquad
\lf\{b^{\{-\}i}(x)\rg\}\stackrel{\Ka }{\to} \lf\{b^{\{-\}i}(x)\rg\}.
\ee

\ses

\ses

{%\pgbrk}

When the substitution (6.26) is applied,
from (3.10)-(3.11) we obtain
 the ${\mathbf\cF\cF^{PD}_{g}}$-{\it angle}
  \be
 \al_{\{x\}}(y_1,y_2) =
 \fr1h\arccos\la
\qquad \text{with} \quad
\la =
\fr{ A(x,y_1)A(x,y_2)+ h^2v_{12}}
 { \sqrt{B(x,y_1)}\,\sqrt{B(x,y_2)} },
 \ee
where
$
 v_{12}=   r_{mn}(x) y_1^my^n_2.
$

If, fixing a point $x$,  we consider the angle
$ \al_{\{x\}}(y,b)$ formed by a vector $y\in T_xM$ with
 the input characteristic vector $b^i(x)$,
from (6.30) we get
the respective value to be
 \be
 \al_{\{x\}}(y,b) =    \fr1h\arccos \fr{ A(x,y)}
 {\sqrt{B(x,y)}}
 \equiv
 \chi
 \ee
(notice that    $q=0$ and $A=1$  whenever   $y=b$).
In terms of the variables $\zeta^n$, the last formula reads
\be
 \al_{\{x\}}(y,b)   = \fr1h    \arccos \fr{\zeta^nb_n(x)}    {\sqrt{ S^2(x,\zeta)}}.
 \ee

We have
\be
A=\sqrt B\,\cos(h\chi), \qquad    qh=\sqrt B\,\sin (h\chi),
\ee
so that the transformation (6.26)  can clearly  be written in terms of the angle
$\chi$:
\be
\zeta^i= \lf(\fr{v^i}q \sin(h\chi)  +   b^i\cos(h\chi) \rg ) K^h.
\ee

\ses

From (6.32) we can conclude
that
\be
 \al_{\{x\}}\lf(b^{\{-\}},b\rg)
   =    \fr1h \pi
\ee
which is {\it more than $\pi$ whenever }
$g\ne0$.

{%\pgbrk}

Let us verify that the transformation (6.26) obeys the input stipulation
(2.1).
Differentiating (6.26) leads to
\be
\zeta^m_n
=E^m_n
+\fr1N \zeta^mC_n
-\fr1{\varkappa}\varkappa_n\zeta^m
\quad  \text{with} \quad
E^m_n
= \Biggl[   h(\de^m_n-b_nb^m)+\lf(  b_n+{}  \fr1{2q}gv_n\rg)b^m\Biggr]J\fr1{\varkappa h},
\ee
\ses
where
\be
\varkappa_n=\partial \varkappa/\partial y^n   =(1-h)y_n \varkappa /K^2  ,
\ee
 and we have used the equality
$\partial \ln J/\partial {y^n} =  (1/N) C_n
$
(which is a direct implication of the formulas (A.6)
and
$ C_n=     \partial{\ln\bigl(\sqrt{\det(g_{ij})}\bigr)}/\partial{y^n}$).
It is useful to  take into account  that
$$
E^m_ny^n=\zeta^m , \qquad    E^m_nb^n=   b^m        J\fr1{\varkappa h}.
$$
Noting also the vanishing
$$
-
\fr{\zeta^i\zeta^j}{N^2}
C_mC^m
-\fr{(1-h)^2}{K^2}\zeta^i\zeta^j
+\fr {1-h}{K^2}
(\zeta^iy^m E_m^j+\zeta^jy^m E_m^i)
=0,
$$
where the equality $K^2C_mC^m=N^2g^2/4$ has been applied (see  (A.7) in Appendix A),
and using
the contravariant components
$g^{mn}$
written in   (A.12) of Appendix A,
we get
$$
g^{mn}\zeta_m^i\zeta_n^j
=
\fr1{N}(\zeta^iC^m E_m^j+\zeta^jC^m E_m^i)
$$

\ses

$$
+
\Biggl[
a^{mn}E_m^iE_n^j
+
 \fr gqbb^mb^nE_m^iE_n^j
-
 \fr gq  (b^my^n+y^mb^n)    E_m^iE_n^j
+
 \fr g{Bq}(b+gq)  y^m  y^n E_m^i  E_n^j
\Biggr]
\fr1{J^2}.
$$
Here,
$$
 a^{mn}E_m^iE_n^j=
h
\Biggl[
   h(a^{ij}-b^ib^j)   +    \fr1{2q}gv^i  b^j
\Biggr] \fr{J^2}{\ka^2h^2}
+
\Biggl[
1-\fr1{2q}gb\Biggr]b^ib^j
\fr{J^2}{\ka^2h^2}
+
b^i  \fr1{2q}g\zeta^j \fr J{\varkappa h}.
$$
We can write
$$
 h^2g^{mn}\zeta_m^i\zeta_n^j
=
-
\fr g{qK^2}(b+gq) \zeta^i\zeta^j
h^2
+
\fr g2\fr1{q}
(\zeta^i b^j+\zeta^jb^i)  J\fr1{\varkappa} h\fr B{K^2}
$$

\ses

\ses

$$
+h^2a^{ij}\fr1{\ka^2}
-
h^2b^ib^j\fr1{\ka^2}
+
\fr1{2q}gb^j hv^i\fr1{\ka^2}
+
\Biggl[
1-
\fr1{2q}gb\Biggr]b^ib^j
\fr1{\ka^2}
+
b^i  \fr1{2q}g\zeta^j   \fr1J\fr1{\varkappa} h
$$

\ses

$$
+
 \fr gqb b^ib^j\fr1{\ka^2}
-
 \fr gq  (b^i\zeta^j+\zeta^ib^j)    \fr1J\fr1{\varkappa} h
+
 \fr g{Bq}(b+gq) \zeta^i\zeta^j
\fr1{J^2}
h^2,
$$
\ses
or
$$
 h^2g^{mn}\zeta_m^i\zeta_n^j
=
-\fr g2\fr1{q}
\zeta^i b^j  J\fr1{\varkappa} h\fr B{K^2}
+h^2a^{ij}\fr1{\ka^2}
-
h^2b^ib^j\fr1{\ka^2}
+
\fr1{2q}gb^j hv^i\fr1{\ka^2}
+
\Biggl[
1+
\fr1{2q}gb\Biggr]b^ib^j
\fr1{\ka^2},
$$
\ses
so that
\be
 g^{mn}\zeta_m^i\zeta_n^j=\fr1{\varkappa^2} a^{ij}.
\ee
The metric tensor transformation (6.38) can be inverted to read
(3.5).
Thus the verification is complete.

{%\pgbrk}

Several interesting relations can be found.
First of all,
constructing the function $S^2=a_{ij}\zeta^i\zeta^j$ with the help of the choice (6.26),
 applying (6.2) and the identity $v^i b_i=0$,
 we obtain the useful
equality
\be
\lf( \fr J{\varkappa h}\rg)^2=\fr{S^2(x,\zeta)}{B(x,y)}.
\ee
Also, from (6.26) it follows that
\be
\zeta^ib_i=\lf(b+\fr12gq \rg) \fr J{\varkappa h},
\qquad
\zeta^i-  (\zeta^nb_n)b^i=
 \fr J{\varkappa h}
 hv^i,
 \ee
and
\be
\sqrt{r_{mn}\zeta^m\zeta^n} =hq
\fr J{\varkappa h},
\qquad
b=
\Bigl(
\zeta^nb_n   -     \fr1{2h}g\sqrt{r_{mn}\zeta^m\zeta^n}
\Bigr)
 \fr {h\varkappa }J.
\ee
From (6.39) and (6.13)
we can obtain the equality
\be
h\varkappa=
\lf(S^2(x,\zeta)\rg)^{(1-h)/(2h)}
\ee
which is equivalent to (3.6).
According to (6.13), (6.32), and (6.33), we have
\be
\fr1J=\e^{\frac12 g \chi}
 \qquad                      \text{with}           \quad
\chi
   = \fr1h    \arccos \fr{\zeta^nb_n(x)}    {\sqrt{ S^2(x,\zeta)}}.
 \ee

The indicated formulas allow us to write down the explicit form of the
inverse
to the transformation (6.26), namely we find
\be
y^i=y^i(x,\zeta) ~~ ~\text{with} ~~
y^i=bb^i+
\fr1h
\Bigl(\zeta^i- (\zeta^nb_n)b^i\Bigr)
 \fr {h\varkappa }J,
 \ee
where $b$ can be taken from (6.41).
It is possible to find straightforwardly the coefficients
$y^i_j = \partial {y^i}/\partial{\zeta^j}$,
obtaining
\ses\\
\be
y^i_j
=
\Biggl[
\Bigl(
b_j  -          \fr {gT}{2h}
r_{jk}\zeta^k
\Bigr)
b^i+
\fr1h
\Bigl(\de^i_j- b_jb^i\Bigr)
\Biggr]
 \fr {h\varkappa }J
+
\Biggl[
\fr{1-h}{hS^2}    \zeta_j
+
\fr {gT}{2h}
\lf(\fr{\zeta^kb_k}{S^2}
\zeta_j-b_j\rg)
\Biggr]
y^i,
\ee
where
$ \zeta_j=a_{jk}\zeta^k   $
and $T=1/\sqrt{r_{mn}\zeta^m\zeta^n}$, which entails
 the useful identities
\be
 y_iy^i_j=  \fr1h\fr{K^2}{S^2}  \zeta_j,
\qquad
\fr1h\fr{K^2}{S^2} \zeta_i\zeta^i_n = y_n.
\ee

{%\pgbrk}

In the rest of this section we assume that
\be
\D g{x^i}=0 ~ : ~ ~ ~  g=const \quad \text{and} \quad h=const.
\ee

\ses

The right-hand part in (6.26)  is such that
$$
\D{\zeta^i}{x^n}=
\D{\zeta^i}{b_j}\partial_nb_j
+
\D{\zeta^i}{a_{mj}}\partial_na_{mj}.
$$
Under these conditions,
straightforward calculations with the help of the representation (3.32)
result in
\be
N^k_n=
-
\Bigl((1-h)b  +\fr12gq    \Bigr)
\fr1h
a^{kj}  \nabla_nb_j
-
\lf(
\fr g{2q}v^k
-(1-h)b^k
\rg)
\fr1h
y^j\nabla_nb_j
-
  a^k{}_{nj}y^j,
\ee
 where
 $
 \nabla_nb_j=\partial b_j/\partial x^n-b_ka^k{}_{nj}.
 $
 Evaluating the coefficients
$
D^k{}_{nm} = - \partial {N^k_n}/\partial {y^m}
$
(see (2.13))
yields
\be
D^k{}_{nm}=
\Bigl((1-h)b_m  +\fr g{2q}v_m    \Bigr)
\fr1h
a^{kj}  \nabla_nb_j
+
      \fr g{2qh} \eta^k_m
y^j\nabla_nb_j
+
\lf(
\fr g{2q}v^k
\!-\!(1\!-\!h)b^k
\rg)
\fr1h
\nabla_nb_m
+  a^k{}_{nm},
\ee
where
\be
\eta^k_m=r^k_m-\fr1{q^2}v^kv_m\equiv a_{mn}\eta^{kn}, \qquad
\eta^{kn}=a^{kn}-\fr1{q^2}v^kv^n.
\ee
This tensor obeys the nullification
\be
y_k\eta^k_m=b_k\eta^k_m=0.
\ee
We can alternatively write (6.48) as follows:
\be
N^k_n=
\Biggl[
ba^{kj}-b^ky^j
-\fr1h\lf(b+\fr12gq\rg)\eta^{kj}
+\fr1h\lf(b^k-\fr1{q^2}
(b+gq)v^k\rg)  y^j
\Biggr]
\nabla_nb_j
-     a^k{}_{nj}y^j,
\ee
\ses
or
\be
N^k_n=\Biggl[\lf(b\!-\!\fr1h\lf(b+\fr12gq\rg)\rg)
\eta^{kj}
\!+\!
\lf(
\fr1{q^2} v^k \lf(b \!  -\!  \fr1h(b+gq)\rg)
\!+
\!\lf(\fr1h\!-\!1\rg)b^k
\rg)
y^j
\Biggr]
 \nabla_nb_j
-     a^k{}_{nj}y^j.
\ee
By  using
$y_k=(u_k+gqb_k) {K^2}/B$
(see  (A.3) in Appendix A), we obtain
\be
y_kN^k_n=-gqJ^2y^j   \nabla_nb_j
-
y_k  a^k{}_{nj}y^j.
\ee

{%\pgbrk}

The subsequent differentiation of the coefficients (6.49) results in
\be
 \D{D^k{}_{nm}}{y^i}
 =
\fr g{2qh}\eta_{mi}
\eta^{kj}  \nabla_nb_j
-       \fr g{2q^3h} (\eta^k_{m}v_i+  \eta^k_{i}v_m)
y^j\nabla_nb_j
         +
\fr g{2qh}(\eta^k_m   \nabla_nb_i+\eta^k_i\nabla_nb_m).
\ee
Owing to the identity   (6.51),
the vanishing
$
y_k \partial{D^k{}_{nm}}/\partial{y^i}=0
$
(see (4.3)) holds true.

{%\pgbrk}

The coefficients (6.48) show the properties
\be
u_kN^k_n=-\fr1h gqy^j \nabla_nb_j-
u_k  a^k{}_{nj}y^j,
\qquad
b_kN^k_n=\fr1h(1-h)y^j \nabla_nb_j
-
b_k  a^k{}_{nj}y^j,
\ee
and
\be
d_nb\equiv
\D{b}{x^n}+b_kN^k_n=\fr1hy^j \nabla_nb_j,
\qquad
d_nq\equiv
\D{q}{x^n}+\fr1q v_kN^k_n=-\fr1{hq}
(b+gq)
 y^j \nabla_nb_j,
\ee
together with
\be
d_n\lf(\fr qb\rg)
=
-\fr1{b^2qh}B
 y^j \nabla_nb_j,
 \quad
d_n B
=
-\fr g{qh}B
 y^j \nabla_nb_j,
   \quad
d_n\fr B{b^2}=
-\fr{2q+g b}
{b^3qh}B
 y^j \nabla_nb_j .
 \ee

\ses

With the formulas (6.48)-(6.58) it is possible to verify directly the validity of the desired
vanishing set
\be
\D K{x^n}+N^m_nl_m=0
\ee
(see (3.26)),
\ses
\be
\D {y_j}{x^n}+N^m_ng_{mj} - D^m{}_{nj}y_m=0
\ee
(see (4.1)),
and
\be
\D {g_{ij}}{x^n}+2N^m_nC_{mji} - D^m{}_{nj}g_{mi}-
D^m{}_{ni}g_{mj}
=0
\ee
(see (4.2)-(4.3)).

\ses

With the help of (6.54) and (6.58) the coefficients
(6.53) can be transformed to
\be
N^k_n
=-l^k\D K{x^n}
+
\Biggl[\lf(b-\fr1h\lf(b+\fr12gq\rg)\rg)
{\cal H}^{kj}\fr{K^2}B
+\lf(\fr1{hq}-\fr1q +\fr{gb}{B}\rg)Km^k
y^j
\Biggr]
 \nabla_nb_j
-    h^k_t a^t{}_{nj}y^j
\ee
(see Appendix B),
where $m^k=(2/Ng)A^k$ and ${\cal H}^{kj}=h^{kj}-m^km^j\equiv (K^2/B) \eta^{kj}$
(see (A.30)--(A.35) in Appendix A).

The entailed coefficients $
 N^k_{nm}=\partial N^k_{n}/\partial y^m$
are found as follows:
$$
 N^k_{nm}
=
-l^k
g\fr qK  \fr{K^2}B   \nabla_nb_m
+
\fr1K\lf[q-\fr1hq+\fr1{2h}g(b+gq)\rg]
{\cal H}^{kj}m_m\fr{K^2}B
 \nabla_nb_j
 $$

\ses

\ses

$$
+
\lf(b-\fr1h\lf(b+\fr12gq\rg)\rg)
\fr{K}B
\lf(
{\cal H}^{kj}l_m-l^k{\cal H}^j_m  -l^j{\cal H}^k_m
+\fr bq( m^j{\cal H}^k_m
+ m^k {\cal H}^j_m)
\rg)
 \nabla_nb_j
$$

\ses

$$
+
K
\fr b{q^2}
\lf(\fr1{h}-1 \rg)
 m_m
 m^k
l^j
 \nabla_nb_j
-
\fr Kq\lf(\fr1{h}-1 \rg)
m_ml^k
l^j
 \nabla_nb_j
$$

\ses

\ses

$$
-
K
\lf(\fr1{hq}(b+gq)-\fr bq\rg)
\fr 1{q}
{\cal H}^k_m
l^j
 \nabla_nb_j
+
\lf(\fr1{hq}-\fr1q +\fr{gb}{B}\rg)Km^k
 \nabla_nb_m
-
a^k{}_{nm}
$$
(see (B.5) in Appendix B).
They  fulfill the equality
 $ N^k_{nm}=- D^k{}_{nm}$ with $D^k{}_{nm}$ given by (6.49).

\ses

{%\pgbrk}

Moreover,  with the coefficients
$N^k_i$ given by (6.48)
we get straightforwardly
the vanishing
\be
\D{\la(x,y_1,y_2)}{x^i}
+ N^k_i(x,y_1)\D{\la(x,y_1,y_2)}{y^k_1}
+ N^k_i(x,y_2)\D{\la(x,y_1,y_2)}{y^k_2}=0, ~~ \text{when $ h=$ const},
\ee
where $\la$ is the scalar indicated in (6.30).
To verify the statement, it is worth
deriving  the equality
\be
\D{\la}{y^k_1}=
h^2\fr{
B_1v_{2k}+q^2_1b_kA_2
- b_1A_2
v_{1k}
-v_{12}
\lf(h^2v_{1k}+\lf(b_k+\fr12g\fr{1}{q_1}v_{1k}\rg)A_1
\rg)
}
 {B_1 \sqrt{B_1}\,\sqrt{B_2} }
\ee
with the counterpart
\be
\D{\la}{y^k_2}=
h^2\fr{
B_2v_{1k}+q^2_2b_kA_1
- b_2A_1
v_{2k}
-v_{12}
\lf(h^2v_{2k}+\lf(b_k+\fr12g\fr{1}{q_2}v_{2k}\rg)A_2
\rg)
}
 {B_2 \sqrt{B_2}\,\sqrt{B_1} },
\ee
where
$
A_1=A(x,y_1),\,   A_2=A(x,y_2),\,
B_1=B(x,y_1),\,   B_2=B(x,y_2),\,
q_1=q(x,y_1),\,q_2=q(x,y_2),\,
b_1=b(x,y_1),\,   b_2=b(x,y_2),$
together with
$
v_{1i}= r_{in}(x)y^n_1$ and $  v_{2i}= r_{in}(x)y^n_2.$
Plugging these derivatives in (6.63) results in the claimed vanishing
after attentive couplepage reductions.

\ses

It will be noted that
$$
b^k\D{\la}{y^k_1}=  h^2\fr{  q^2_1A_2  -v_{12}   A_1   }
 {B_1 \sqrt{B_1}\,\sqrt{B_2} },  \qquad
 b^k\D{\la}{y^k_2}=  h^2\fr{  q^2_2A_1  -v_{12}   A_2   }
 {B_2 \sqrt{B_1}\,\sqrt{B_2} }.
$$

{%\pgbrk}

From (6.30) we have also
$$
\D{\la}g=-\fr12\lf(\fr {b_1q_1}{{B_1}}+\fr{b_2q_2}{{B_2}}\rg)  \la
+\fr{q_1A_2+q_2A_1-gv_{12}}
 {2 \sqrt{B_1}\,\sqrt{B_2} },
$$
or
$$
\D{\la}g=
\fr{1}
 { 2\sqrt{B_1}\,\sqrt{B_2} }
 \lf[
\fr{q_1^2 A_2}{B_1}\si_1
+
\fr{q_2^2A_1}{B_2}\si_2
-v_{12}
 \lf(
\fr{ A_1}{B_1}\si_1
+
\fr{ A_2}{B_2}\si_2
 \rg)
 \rg],
 $$
where
\be
\si_1=\fr g2A_1+h^2q_1\equiv q_1+\fr g2b_1,  \qquad   \si_2=\fr g2A_2+h^2q_2\equiv q_2+\fr g2b_2.
\ee
There arises the equality
\be
\D{\la}g=\fr1{2h^2}\lf[
\si_1b^k  \D{\la}{y^k_1}  +   \si_2b^k  \D{\la}{y^k_2}
\rg].
\ee
Using  the formula (A.26)
of Appendix A,
we arrive at
\be
\D{\la}g=\fr1{h^2}\lf[
\mu_1C_1^k  \D{\la}{y^k_1}  +   \mu_2C_2^k  \D{\la}{y^k_2}
\rg],
\ee
where
\be
\mu_1=   \fr{q_1K_1^2}{NgB_1}   \si_1,  \qquad
\mu_2=   \fr{q_2K_2^2}{NgB_2}   \si_2.
\ee

\ses

{%\pgbrk}

 The associated Riemannian curvature tensor is constructed as follows:
\be
a_n{}^i{}_{km}=\D{a^i{}_{nm}}{x^k}-\D{a^i{}_{nk}}{x^m}+a^u{}_{nm}a^i{}_{uk}-a^u{}_{nk}a^i{}_{um}.
\ee
The evaluation of the  tensor (4.15) from the coefficients (6.48) gives us
\be
 M^n{}_{ij}=
\Biggl[
\Bigl((1-h)b  +\fr12gq    \Bigr)
a^{nt}
+
\lf(
\fr g{2q}v^n
-(1-h)b^n
\rg)  y^t
\Biggr]
\fr1h
b_la_t{}^l{}_{ij}
-
a_t{}^n{}_{ij}
 y^t,
\ee
or
\be
 M_{nij}= \Biggl[
  \Bigl((1-h)b  +\fr12gq    \Bigr)
\fr1h
b^la_{tlij}
+
y^la_{tlij}
\Biggr]
\fr {K^2}B{\cal H}_n^t
- \fr {1}{q}
Km_n
 \fr1h
y^tb^la_{tlij}
\ee
(see (A.45)),
which entails the equalities
\be
\fr1NC_nM^n{}_{ij}= -\fr g{2qh}
b_n a_t{}^n{}_{ij}y^t, \quad
A_{knm}M^m{}_{ij}=
-
K
 {\cal H}_{kn}
\fr g{2qh}
b_m a_t{}^m{}_{ij}y^t
+\fr1N (A_kM_{nij}+A_nM_{kij})
\ee
((A.8) has been used and
the tensor ${\cal H}_{kn}$
has been defined in (A.30)),
and
\be
M^n{}_{ij}=
-
y^n_t\zeta^h
a_h{}^t{}_{ij}.
\ee
The  tensor (4.14) is found to read
\be
E_k{}^n{}_{ij}=
-
\Biggl[
\Bigl((1-h)b_k  +\fr12\fr gqv_k   \Bigr)    a^{nt}
+   \fr g{2q} \eta^n_k  y^t
+
\lf(
\fr g{2q}v^n
-(1-h)b^n
\rg)
 \de^t_k
\Biggr]
\fr1h
b_ma_t{}^m{}_{ij}
+  a_k{}^n{}_{ij},
\ee
which entails
\be
E_k{}^n{}_{ij}=
y_h^n\zeta^h_{km}M^m{}_{ij}
+
y^n_ma_h{}^m{}_{ij} \zeta^h_k.
\ee

Obeying the identities (4.16)-(4.18) can straightforwardly be verified.

\ses

{%\pgbrk}

 The following explicit representation
for the curvature tensor
(4.22) can be proposed:
\be
\Rho_{knij}=-\fr1K(l_kM_{nij} - l_nM_{kij})
+(m_k{\cal H}_n^t-m_n{\cal H}_k^t)
P_{tij}
 +
{\cal H}_k^t {\cal H}_n^l  a_{tlij}
 \fr {K^2}B
 \ee
with
$$
P_{tij}=
\Biggl[
-\lf[hq^2+b\lf(b+\fr12gq\rg)\rg]
b^la_{tlij}
+
\lf(b+\fr12gq\rg)
y^l a_{tlij}
\Biggr]
\fr {K}{qB},
$$
\ses
or
\be
P_{tij}=
\Biggl[
-hq^2
b^la_{tlij}
+
\lf(b+\fr12gq\rg)
v^l a_{tlij}
\Biggr]
\fr {K}{qB}
\ee
(see (C.14)-(C.17) in Appendix C).

{%\pgbrk}

With $\Rho_{knij}$ we associate the tensor
$$
\Rho^{knij}=g^{pk}g^{qn}a^{mi}a^{nj}\Rho_{pqmn}.
$$
We can straightforwardly obtain  the contraction
\be
\Rho^{knij}\Rho_{knij}=
a^{knij}a_{knij}
+\fr 2{S^2}
\lf(\fr1{h^2}-1\rg)
\zeta^l
a_l{}^{nij}
\zeta^h
a_{hnij}
\ee
\ses
(see Appendix C).

{%\pgbrk}

We can also find  that
the tensor
\be
 M_{nij}= g_{nm}M^m{}_{ij}
 \ee
 possesses the simple representation
\be
\fr B{K^2} M_{nij}=
\Bigl((1-h)b  +\fr12gq    \Bigr)
\fr1h
b_la_n{}^l{}_{ij}
-
\lf(
\fr g{2q}v_n
+
(1-h)b_n
\rg)
\fr1h
 y^t
 b_la_t{}^l{}_{ij}
-
a_{tnij}
 y^t.
\ee
The identity
$
y^n M_{nij}=0
$
holds.
Squaring this tensor leads to
the quadratic expressions
\be
\fr B{K^2}M^{nij} M_{nij}
=
\Biggl(\fr1{h}
\Bigl((1-h)b  +\fr g2q    \Bigr)
b_ha^{nhij}
-
a_h{}^{nij}
 y^h
\Biggr)
\Biggl(\fr1{h}
\Bigl((1-h)b  +\fr g2q    \Bigr)
b_la_n{}^l{}_{ij}
-
a_{tnij}
 y^t
\Biggr)
\ee
and
\be
 M^{nij} M_{nij}=
\ka^2
\zeta^l
a_l{}^{nij}
\zeta^h
a_{hnij}
\ee
(as shown in Appendix A); here,
$\ka^2= {K^2}/h^2S^2$ in accordance with (3.6).

\ses

{%\pgbrk}

Using (6.76)  together with
\be
y_h^n\zeta^h_{km}M^m{}_{ij}=
C^n{}_{km}M^m{}_{ij}+(1-h)
\fr1{K^2}\lf(
 y^ng_{km}M^m{}_{ij}
-y_k
 M^n{}_{ij}
\rg)
\ee
(see (D.4) in Appendix D)
 reduces the curvature tensor
(4.22) to the sum
\be
\Rho_k{}^n{}_{ij}=
y^n_ma_h{}^m{}_{ij} \zeta^h_k
 +(1-h)
\fr1{K^2}\lf(
 y^ng_{km}M^m{}_{ij}
-y_k
 M^n{}_{ij}
\rg).
\ee

 Make the transform
$$
y^t_k\zeta^n_l\Rho_t{}^l{}_{ij}
=
a_k{}^n{}_{ij}
 +(1-h)
\fr1{K^2}\lf(
h\zeta^n
\fr{K^2}{h^2S^2}
a_{kr}
\zeta^r_lM^l{}_{ij}
-
 \fr1h\fr{K^2}{S^2}  \zeta_k
\zeta^n_l M^l{}_{ij}
\rg),
$$
where
(3.4)-(3.6) and (6.46) have been taken into account.
Using
(6.74) leads to
\be
y^t_k\zeta^n_l\Rho_t{}^l{}_{ij}
=
a_k{}^n{}_{ij}
 +(1-h)
\fr1{hS^2}
( \de^n_l\zeta_k-\zeta^na_{kl})
\zeta^h
a_h{}^l{}_{ij}.
\ee

\ses

\ses

Now we contract
\ses\\
$$
\Rho_l{}^t{}^{ij}
\Rho_t{}^l{}_{ij}
=
y^p_n\zeta^k_q\Rho_p{}^q{}^{ij}
y^t_k\zeta^n_l\Rho_t{}^l{}_{ij},
$$
\ses
so that
$$
\Rho_l{}^t{}^{ij}
\Rho_t{}^l{}_{ij}
=
y^p_n\zeta^k_q\Rho_p{}^q{}^{ij}
\lf[
a_k{}^n{}_{ij}
 +(1-h)
\fr1{hS^2}
( \de^n_l\zeta_k-\zeta^na_{kl})
\zeta^h
a_h{}^l{}_{ij}
\rg]
$$

\ses

\ses

$$
=
a_n{}^{kij}
\lf[
a_k{}^n{}_{ij}
 +(1-h)
\fr1{hS^2}
( \de^n_l\zeta_k-\zeta^na_{kl})
\zeta^h
a_h{}^l{}_{ij}
\rg]
$$

\ses

\ses

$$
 +(1-h)
\fr1{hS^2}
\zeta_n
\zeta^s
a_s{}^{kij}
\lf[
a_k{}^n{}_{ij}
 +(1-h)
\fr1{hS^2}
( \de^n_l\zeta_k-\zeta^na_{kl})
\zeta^h
a_h{}^l{}_{ij}
\rg]
$$

\ses

\ses

$$
-
(1-h)
\fr1{hS^2}
\zeta^k
\zeta^s
a_{nl}a_s{}^{lij}
\lf[
a_k{}^n{}_{ij}
 +(1-h)
\fr1{hS^2}
( \de^n_l\zeta_k-\zeta^na_{kl})
\zeta^h
a_h{}^l{}_{ij}
\rg],
$$
\ses
or
$$
\Rho_l{}^t{}^{ij}
\Rho_t{}^l{}_{ij}
=
a_n{}^{kij}
a_k{}^n{}_{ij}
 +(1-h)
\fr1{hS^2}
\zeta_n
\zeta^s
a_s{}^{kij}
\lf[
2
a_k{}^n{}_{ij}
-
(1-h)
\fr1{hS^2}
\zeta^na_{kl}
\zeta^h
a_h{}^l{}_{ij}
\rg]
$$

\ses

\ses

$$
-
(1-h)
\fr1{hS^2}
\zeta^k
\zeta^s
a_{np}a_s{}^{pij}
\lf[
2
a_k{}^n{}_{ij}
 +(1-h)
\fr1{hS^2}
 \de^n_l\zeta_k
\zeta^h
a_h{}^l{}_{ij}
\rg]
$$

\ses

\ses

\ses

$$
=
a_n{}^{kij}
a_k{}^n{}_{ij}
-
2(1-h)
\fr1{hS^2}
\zeta^s
a_s{}^{lij}
\lf[
2
\zeta^k
a_{klij}
 +
 (1-h)
\fr1{h}
\zeta^h
a_{hlij}
\rg].
$$
\ses
The result
$$
\Rho_l{}^t{}^{ij}
\Rho_t{}^l{}_{ij}
=
a_n{}^{kij}
a_k{}^n{}_{ij}
-
2(1-h^2)
\fr1{h^2S^2}
\zeta^s
a_s{}^{lij}
\zeta^k
a_{klij}
$$
\ses
is equivalent to (6.79).

{%\pgbrk}

Let us introduce the object
$
Y^n_k(x,y)=y^n_k(x,\zeta).
$
We get
$$
\D{Y^n_k}{x^i}=\D{y^n_k}{x^i}+y^n_{ks}\D{\zeta^s}{x^i}.
$$
Using here  (4.4) and (3.31) leads to
$$
\D{Y^n_k}{x^i}=
-L^t_i y^n_{kt}-
D^n{}_{is}y^s_k
+
L^h{}_{ik}y^n_h
-
y^n_{ks}
\lf(
N^r_i\zeta^s_r+a^s{}_{ui}\zeta^u
\rg).
$$
Therefore, from the representation
 $M^n{}_{ij}=-Y^n_t\zeta^ha_h{}^t{}_{ij}$
(see (6.74))
we find
the partial derivative
$$
\D{M^n{}_{ij}}{x^k}
+
D^n{}_{ks}M^s{}_{ij}
=
\lf[
L^h_i y^n_{th}
-
L^h{}_{kt}y^n_h
+
y^n_{ts}
\lf(
N^r_k\zeta^s_r+a^s{}_{uk}\zeta^u
\rg)
\rg]
\zeta^h
a_h{}^t{}_{ij}
$$

\ses

$$
+
Y^n_t
\lf(
N^u_k\zeta^h_u+a^h{}_{uk}\zeta^u
\rg)
a_h{}^t{}_{ij}
-
Y^n_t\zeta^h
\D{a_h{}^t{}_{ij}}{x^k}.
$$
The covariant derivative
 (4.20) can now be written in the form
$$
D_k M^n{}_{ij}
=
\lf[
L^u_i y^n_{tu}
-
L^u{}_{kt}y^n_u
+
y^n_{ts}
\lf(
N^r_k\zeta^s_r+a^s{}_{uk}\zeta^u
\rg)
\rg]
\zeta^h
a_h{}^t{}_{ij}
+
Y^n_t
\lf(
N^u_k\zeta^h_u+a^h{}_{uk}\zeta^u
\rg)
a_h{}^t{}_{ij}
$$

\ses

$$
-
Y^n_t\zeta^h
\D{a_h{}^t{}_{ij}}{x^k}
-
N^m_k
\lf[
y^n_{ts}\zeta^s_m\zeta^h
a_h{}^t{}_{ij}
+
Y^n_t\zeta^h_m
a_h{}^t{}_{ij}
\rg]
 - a^s{}_{ki}   M^n{}_{sj}  -  a^s{}_{kj}   M^n{}_{is}.
 $$
Cancelling here similar terms leaves us
with
$$
D_k M^n{}_{ij}
=
\lf(
L^u_i y^n_{tu}
-
L^u{}_{kt}y^n_u
+
y^n_{ts}
a^s{}_{uk}\zeta^u
\rg)
\zeta^h
a_h{}^t{}_{ij}
$$

\ses

$$
+
y^n_t
a^h{}_{uk}\zeta^u
a_h{}^t{}_{ij}
-
y^n_t\zeta^h
\D{a_h{}^t{}_{ij}}{x^k}
 - a^s{}_{ki}   M^n{}_{sj}  -  a^s{}_{kj}   M^n{}_{is}.
$$
Recollecting the equalities
$  L^m_{j}=-L^m{}_{ij}\zeta^i$ and $   L^m{}_{ij}=a^m{}_{ij}$
indicated in (3.22),
we obtain simply
\be
D_k M^n{}_{ij}
=
-
y^n_t\zeta^h
\nabla_k a_h{}^t{}_{ij},
\ee
\ses
where
\be
\nabla_k a_h{}^t{}_{ij}
=
\D{a_h{}^t{}_{ij}}{x^k}
+
 a^t{}_{ku}  a_h{}^u{}_{ij}
- a^u{}_{kh}  a_u{}^t{}_{ij}
- a^u{}_{ki}  a_h{}^t{}_{uj}
- a^u{}_{kj}  a_h{}^t{}_{iu}
\ee
is the Riemannian covariant derivative of the Riemannian curvature tensor.

The cyclic identity (4.19) proves to be a direct implication of the known Riemannian
identity
\be
\nabla_ka^n{}_{ij}+\nabla_ja^n{}_{ki}+\nabla_ia^n{}_{jk}=0.
\ee

{%\pgbrk}

Using the equality
\be
g_{nm}y^m_i=\ka^2\zeta^j_n a_{ij}
\ee
(ensued from (3.5)),
we can obtain the tensor (6.80) to read
\be
M_{nij}=-\ka^2\zeta^h\zeta^m_na_{hmij}
\ee
and
juxtapose to (6.85) the tensor
$\Rho_{knij}=
g_{mn}\Rho_k{}^m{}_{ij}$
which is
\be
\Rho_{knij}=T_{kn}{}^{hm}a_{hmij},
\ee
where
\be
T_{kn}{}^{hm}= \ka^2
\lf[
\fr12(\zeta^h_k\zeta^m_n-\zeta^m_k\zeta^h_n)
+(1-h)\fr1{K^2}(y_k\zeta^h\zeta^m_n-y_n\zeta^h\zeta^m_k)
\rg].
\ee

Since
\be
D_lT_{kn}{}^{hm}=0,
\ee
we have
\be
D_l\Rho_{knij}=T_{kn}{}^{hm}\nabla_la_{hmij},
\ee
together with the cyclic identity
\be
D_l\Rho_{knij}+D_j\Rho_{knli}  +D_i\Rho_{knjl} =0.
\ee

{%\pgbrk}

\ses

\ses

             \setcounter{sctn}{7}
\setcounter{equation}{0}

\nin
{\bf 7.
${\mathbf\cF\cF^{PD}_{g}}$-space coordinates and angles}

\ses

\ses

We now  fix the tangent space (in accordance with  Appendix E)
and choose  the three-dimensional case
\be
N=3, \qquad  R^p=\{R^1,R^2,R^3\}.
\ee
It is convenient to relabel the coordinates $R^p$ as follows:
\be
R^1=x, \quad  R^2=y,  \quad R^3=z.
\ee
We get
\be
q=\sqrt{x^2+y^2},  \qquad
B= x^2 +y^2 + z^2    +    gz q.
\ee
In terms of such coordinates, the metric tensor components
$g_{pq}$ can be  obtained
 from  the list  (E.6)-(E.7).
The result reads
\be
  g_{11}=\Bigl(   1 - \fr {g}{Bq}zx^2   \Bigr)  J^2,
    \qquad
g_{22}=\Bigl(   1 - \fr {g}{Bq} zy^2   \Bigr)  J^2,
 \qquad
g_{33}=\lf(  1+\fr {gq} {B} (z+gq) \rg)  J^2,
\ee

\ses

\be
g_{12}= -\fr {g}{Bq}zxy    J^2,
 \qquad
g_{13}=    \fr {gqx}{B}    J^2,
 \qquad
g_{23}=\fr {gqy}{B}   J^2.
  \ee
From the formulas (E.8)-(E.9) it follows that

\be
g^{11}=\Bigl(   1+\fr {g}{Bq}(z+gq) x^2   \Bigr)  \fr1{J^2},
   \quad
g^{22}= \Bigl(   1+  \fr {g}{Bq}(z+gq) y^2   \Bigr)   \fr1{J^2},
\quad
g^{33}=\Bigl(  1  -\fr {gq}{B}z
\Bigr)  \fr1{J^2},
\ee

\ses

\be
g^{12}=
\fr {g}{Bq} (z+gq) xy      \fr1{J^2},
 \qquad
g^{13}=    - \fr {gq}{B}x  \fr1{J^2},
 \qquad
g^{23}= - \fr {gq}{B}y  \fr1{J^2}.
\ee

\ses

The
${\mathbf\cF\cF^{PD}_{g}}$-{\it space coordinates}
$\{z^p\}$
are given by
\be
z^1=K,  \qquad z^2=\phi, \qquad z^3=\chi \equiv \fr1h f,
\ee
where $K$ is the Finsleroid metric function (6.13),
 $\phi$ is the {\it polar angle} in the $R^1\times R^2$-plane,
and
$\chi$ plays the role of the    Finsleroid  {\it azimuthal angle}
 measured from the direction
of the input vector $b^i$ (see (6.32)).
The indices $p,q,...$ will be specified over the range 1,2,3.
For the vector $\{R^p\}$ we
 construct
 the representation
\be
R^p=R^p(g;z^q)
\ee
which
  possesses the {\it invariance property}
\be
K\bigl(g;R^p(g;z^q)\bigr)= z^1.
\ee

{%\pgbrk}

The  representations
\ses
\be
R^1=K \Sin\chi \cos\phi,
\qquad
R^2=K \Sin\chi \sin\phi,
\qquad
R^3=K \Cos\chi,
\ee
with
\be
\Sin\chi=\fr1{Jh} \sin f, \qquad
\Cos\chi=
\fr1J\lf(\cos f-\fr G2\sin f\rg),
\ee
can readily be arrived at,
entailing
\be
  q=K\Sin\chi, \qquad        b+\fr12gq =\fr KJ\cos f,  \qquad b=K \Cos\chi.
\ee
The arisen functions
$\Sin\chi $ and $\Cos\chi $
can be interpreted as
the required extensions of  the trigonometric functions
to the ${\mathbf\cF\cF^{PD}_{g}}$-space.

The squared  linear element  $ds^2= g_{rs}(R)dR^rdR^s$
is found to be of the diagonal form
\be
(ds)^2=(dz^1)^2 +
(z^1)^2\lf[(d\chi)^2+
\fr1{h^2}\sin^2 (h\chi)   (d\phi)^2
\rg].
\ee

On the other hand,
when the components (7.11)  are inserted in   the ${\mathbf\cF\cF^{PD}_{g}}$-angle (6.30),
 the following result is obtained:
\be
 \al_{\{x\}}(y_1,y_2)    = \fr1h\arccos \tau_{12},
 \qquad \text {with}
\quad
\tau_{12}=
\cos(f_2-f_1)-(1-\cos(\phi_2-\phi_1)\bigr)\sin f_1\sin f_2.
\ee
This $\tau_{12}$
 does not involve any support vector $y$.
By developing here the infinitesimal version,
 putting $\chi_1=\chi$, $ \chi_2=\chi+d\chi$,
 $\phi_1=\phi$, $ \phi_2=\phi+d\phi$,
we come to the infinitesimal angle
$d\al=d\al_{\{x\}}$
which square reads
\be
(d\al)^2=
(d\chi)^2+
\fr1{h^2}\sin^2 (h\chi)\,   (d\phi)^2.
\ee
By comparing (7.14) with (7.16) we conclude that
\be
(ds)^2=(dz^1)^2+(z^1)^2(d\al)^2.
\ee
This formula
is remarkable because showing us frankly that
 $d\al$ {\it is the infinitesimal arc-length on the indicatrix}
 (keeping in mind that  $z^1=1$  holds along the indicatrix).

The obtained metric
(7.14)
is of the  conformally flat type
\be
(ds)^2= \varkappa^2 (ds)^2{\Bigl|\Bigr.}_{\text{Euclidean}}
\qquad
\text{with} \quad  \varkappa=\fr1{h}K^{1-h}.
\ee

        {%\pgbrk}

To verify this assertion,
it is appropriate to
use the substitution
$z^1=\e^{\si}$  in (7.14),
which  yields
\be
(ds)^2=    e^{2\si}
\Bigl[
(d\chi)^2+
\fr1{h^2}\sin^2 (h\chi)   (d\phi)^2
+
(d\si)^2
\Bigr].
\ee
\ses
Making  the coordinate transformation
\be
\rho=\e^{h\si}\sin(h\chi),
\qquad
\tau=\e^{h\si}\cos(h\chi)
\ee
 leads to (7.18) with
$$
\Bigl((ds)^2\Bigr){\Bigl|\Bigr.}_{\text{Euclidean}}=
(d\rho)^2+
\rho^2\sin^2 (h\chi)   (d\phi)^2
+(d\tau)^2
$$
and
$$
\ka=\fr1{h}(\rho^2+\tau^2)^{(1-h)/2h}.
$$

The observations can be summarized by  formulating the following

\ses

\ses

{\bf Proposition.}
{\it
Given  an ${\mathbf\cF\cF^{PD}_{g}}$-space of the dimension $N=3$.
In terms of the Finsleroid coordinates \rm(7.2),
{\it
which directly extend the spherical coordinates
applied conventionally in the tangent spaces to the three-dimensional Riemannian space,
the induced metric on the indicatrix is of the diagonal representation
}
 (7.14).
{\it
The metric is of the conformally flat type as shown by
}
(7.18).

\ses

\ses

{%\pgbrk}

According to
         (7.17),
         the length element on the indicatrix is given by the representation
\be
ds{\Bigl|\Bigr.}_{\text{${\mathbf\cF\cF^{PD}_{g}}$-indicatrix}}=\sqrt{d\chi^2+\fr1{h^2}\sin^2(h \chi)d\phi^2},
\ee
which  can be used to
 find the geodesics which are the solutions of the Euler-Lagrange equation
 written by the help of the Lagrangian
$L=\sqrt{(d\chi/dt)^2+(1/{h^2})\sin^2(h \chi)(d\phi/dt)^2},
$
where $t$ is an appropriate parameter.
Since $\phi$ is a cyclic coordinate, we have
\be
\fr1{h^2}\sin^2(h \chi)\phi'=\wt C,
\ee
where $\wt C$ is a constant,
thereafter  we get
\be
\chi''=h^3 \wt C^2 \fr{ \cos(h\chi)}{\sin^3(h\chi)}
\ee
and
\be
\chi'=\sqrt{1-h^2 \wt C^2 \fr1{\sin^2(h\chi) }}.
\ee
The prime $'$ means differentiation with respect to the parameter $s$ defined by (7.21).
It follows that
\be
\lf(\fr1h  \sin f\chi'\rg)'
=
\cos f\chi'\chi'
+
\fr1h  \sin f\chi''
=\cos f.
\ee

The  equation (7.24)
can readily be integrated,
 yielding  the explicit dependence
\be
\chi(s)=\fr1h\arccos\lf(\sqrt{1-h^2\wt C^2}\cos\bigl(h(s-\wt s)\bigr)\rg),
\ee
where $\wt s$ is an integration constant.
So,
\be
\cos(h\chi)=\sqrt{1-h^2\wt C^2}\cos\bigl(h(s-\wt s)\bigr),
\quad
\sin(h\chi)=
\sqrt{
1-\lf(1-h^2\wt C^2\rg)\cos^2\bigl(h(s-\wt s)\bigr)
},
\ee
and from (7.22) we get
\be
\phi'= h\wt C\fr{h}{1-\lf(1-h^2\wt C^2\rg)\cos^2\bigl(h(s-\wt s)\bigr)}.
\ee
Integrating yields  explicitly
\be
\phi(s)=\arctan\lf(\fr1{h\wt C}\tan\bigl(h(s-\wt s)\bigr)\rg),  ~~ \text{if} ~ \wt C\ne0;
\quad
\phi=\fr{\pi}2,  ~~ \text{if} ~ \wt C = 0,
\ee
from which we have
\be
\cos\phi=\fr{h\wt C} {\sqrt{1-\lf(1-h^2\wt C^2\rg)\cos^2\bigl(h(s-\wt s)\bigr)}},
\quad
\sin\phi=\fr{\sqrt{1-h^2\wt C^2}\,\sin \bigl(h(s-\wt s)\bigr)}
 {\sqrt{1-\lf(1-h^2\wt C^2\rg)\cos^2\bigl(h(s-\wt s)\bigr)}}.
\ee

{%\pgbrk}

With these representations, we are able to obtain from the formulas (7.11)-(7.13)
 the explicit behavior of the unit vector components
$l^1=R^1/K$,
$l^2=R^2/K$,
and
$l^3=R^3/K$
along the geodesic arc. The result reads
\be
 \e^{-\frac12 g \chi(s)}
 l^1(s)=  \wt C \sin (hs),
\ee
\ses
\be
 \e^{-\frac12 g \chi(s)}
 l^2(s)=  \fr1h \sqrt{1-h^2\wt C^2}\,\sin (hs),
\ee
and
\be
\e^{-\frac12 g \chi(s)}
l^3(s)=
\sqrt{1-h^2\wt C^2}\cos(hs)
-\fr G2
\sqrt{
1-\lf(1-h^2\wt C^2\rg)\cos^2(hs)
},
\ee
where $\chi(s)$ is the function (7.26) and we have put $\wt s=0$.

Given a geodesic-arc ${\cal A}(x,l_1,l_2)$. Let the left-side vector $l_1$
correspond to  $s=s_1$,
and the right-side vector $l_2$  relate to a value $s_2>s_1$.
From (7.31)-(7.33) we obtain
\be
 \e^{-\frac12 g \chi_1}
 l^1_1=  \wt C \sin (hs_1),
\ee
\ses
\be
 \e^{-\frac12 g \chi_1}
 l^2_1=  \fr1h \sqrt{1-h^2\wt C^2}\,\sin (hs_1),
\ee
\ses
\be
\e^{-\frac12 g \chi_1}
l^3_1=
\sqrt{1-h^2\wt C^2}\cos(hs_1)
-\fr G2
\sqrt{
1-\lf(1-h^2\wt C^2\rg)\cos^2(hs_1)
},
\ee
where
\be
\chi_1=\fr1h\arccos\lf(\sqrt{1-h^2\wt C^2}\cos(hs_1)\rg),
\ee
\ses
and
\be
 \e^{-\frac12 g \chi_2}
 l^1_2=  \wt C \sin (hs_2),
\ee
\ses
\be
 \e^{-\frac12 g \chi_2}
 l^2_2=  \fr1h \sqrt{1-h^2\wt C^2}\,\sin (hs_2),
\ee
\ses
\be
\e^{-\frac12 g \chi_2}
l^3_2=
\sqrt{1-h^2\wt C^2}\cos(hs_2)
-\fr G2
\sqrt{
1-\lf(1-h^2\wt C^2\rg)\cos^2(hs_2)
},
\ee
where
\be
\chi_2=\fr1h\arccos\lf(\sqrt{1-h^2\wt C^2}\cos(hs_2)\rg),\ee

{%\pgbrk}

With the last formulas, the representations (7.31)-(7.33) can be written as follows:
\be
 l^1(s)=  k_1  l^1_1 +k_2  l^1_2,
\qquad
 l^2(s)= k_1  l^2_1+   k_2  l^2_2,
\qquad
 l^3(s)= k_1  l^3_1 + k_2  l^3_2  +k_3 ,
\ee
\ses
where
\be
k_1=
\lf[
\fr{\sin (h(s-s_1))}
{\sin(h(s_2-s_1))}
\fr{\sin (hs_2)}
{\sin(hs_1)}
+
\fr{\sin (h(s_2+s_1))}
{\sin(h(s_2-s_1))  \, \sin(hs_1)}
\rg]
\e^{\frac12 g (\chi(s)-\chi_1)},
\ee
\ses
\be
k_2=
\lf[
\fr{\sin (h(s_2-s))}
{\sin(h(s_2-s_1))}
\fr{\sin (hs_1)}
{\sin(hs_2)}
-
\fr{\sin (h(s_2+s_1))}
{\sin(h(s_2-s_1))  \, \sin(hs_2)}
\rg]
\e^{\frac12 g (\chi(s)-\chi_2)},
\ee
\ses
\be
k_3=\fr g{2h} Y,
\ee
\ses
and
$$
Y=
-
\sqrt{
1-\lf(1-h^2\wt C^2\rg)\cos^2(hs)
}
\,\e^{\frac12 g \chi(s)}
+k_1
\sqrt{
1-\lf(1-h^2\wt C^2\rg)\cos^2(hs_1)
}
\,\e^{\frac12 g \chi_1}
$$

\ses

\be
+k_2
\sqrt{
1-\lf(1-h^2\wt C^2\rg)\cos^2(hs_2)
}
\,\e^{\frac12 g \chi_2}.
\ee
Thus we have arrived at the vector representation
\be
l^p=k_1l_1^p+k_2l^p_2+k_3\de^p_3.
\ee
With respect to  arbitrary local coordinates $x^i$,
we eventually obtain
the expansion
\be
l^i(s)=k_1(s)l^i_1+k_2(s)l^i_2 +k_3(s)b^i,
\ee
which does involve the vector $b^i$ in addition to
$l^i_1,l^i_2.$

The function (7.46) possesses the property
\be
\wt C=0 \quad \Rightarrow \quad Y=0
\ee
(at any value of $g$).

{%\pgbrk}

\ses

\ses

\setcounter{equation}{0}

{\nin \bf  Appendix A:   Involved  ${\mathbf\cF\cF^{PD}_{g} } $-notions   }

\ses

\ses

By $K$ we denote the metric function obtainable from the formulas (2.14) and (6.1)--(6.6).

\ses

\ses

 {\large  Definition}.  Within  any tangent space $T_xM$, the   function $K(x,y)$
  produces the {\it    ${\mathbf\cF\cF^{PD}_{g} } $-Finsleroid}
 \be
 \cF\cF^{PD}_{g;\,\{x\}}:=\{y\in   \cF\cF^{PD}_{g; \, \{x\}}: y\in T_xM , K(x,y)\le 1\}.
  \ee

\ses

 \ses

 {\large  Definition}. The {\it    ${\mathbf\cF\cF^{PD}_{g} } $-Indicatrix}
 $ {\cal I}\cF^{PD}_{g; \, \{x\}} \subset T_xM$ is the boundary of the
    ${\mathbf\cF\cF^{PD}_{g} } $-Finsleroid, that is,
 \be
{\cal I}\cF^{PD}_{g\, \{x\}} :=\{y\in {\cal I}\cF^{PD}_{g\, \{x\}} : y\in T_xM, K(x,y)=1\}.
  \ee

\ses

 \ses

 {\large  Definition}. The scalar $g(x)$ is called
the {\it Finsleroid charge}.
The 1-form $b=b_i(x)y^i$ is called the  {\it Finsleroid--axis}  1-{\it form}.

\ses

\ses

We can   explicitly extract from the function $K$ the  distinguished Finslerian tensors,
 and first of all
the covariant tangent vector $\hat y=\{y_i\}$ from $y_i :=(1/2)\partial {K^2}/ \partial{y^i}$,
obtaining
\be
y_i=(u_i+gqb_i) \fr{K^2}B,
\ee
where $u_i=a_{ij}y^j$.
After that, we can find
the  Finslerian metric tensor $\{g_{ij}\}$
together with the contravariant tensor $\{g^{ij}\}$ defined by the reciprocity conditions
$g_{ij}g^{jk}=\de^k_i$, and the  angular metric tensor
$\{h_{ij}\}$, by making  use of the following conventional  Finslerian  rules in succession:
$$
g_{ij} :
=
\fr12\,
\fr{\prtl^2K^2}{\prtl y^i\prtl y^j}
=\fr{\prtl y_i}{\prtl y^j}, \qquad
h_{ij} := g_{ij}-y_iy_j\fr1{K^2},
$$
thereafter   the Cartan tensor
\be
 A_{ijk}~ := \fr K2\D{g_{ij}}{y^k}
 \ee
and the contraction
\be
 A_i~:=g^{jk}A_{ijk} =  K\D{\ln\bigl(\sqrt{\det(g_{mn})}\bigr)}{y^i}
\ee
can readily be evaluated.

It can straightforwardly be verified that
\be
\det(g_{ij})=\biggl(\fr{K^2}B\biggr)^N\det(a_{ij})>0.
\ee
Contracting the components $A_i$ and $A^i$ yields the formula
\be
A^iA_i=     \fr{N^2g^2}4
\ee

{%\pgbrk}

\nin
and
 evaluating  the Cartan tensor results in the lucid representation
\be
A_{ijk}= \fr1N
 \lf[
A_ih_{jk}  +A_jh_{ik}  +A_kh_{ij}
-\fr 4{N^2g^2}
A_iA_jA_k
\rg].
\ee

\ses

If we insert (A.8)  into the indicatrix curvature tensor (5.3), we obtain
the representation   (5.6)
which manifests  that
 in the ${\mathbf\cF\cF^{PD}_{g}}$-space the indicatrix is of  constant positive
 curvature
(in compliance  with (2.6)).

We use
 the Riemannian covariant derivative
\be
\nabla_ib_j~:=\partial_ib_j-b_ka^k{}_{ij},
\ee
where
\be
a^k{}_{ij}~:=\fr12a^{kn}(\prtl_ja_{ni}+\prtl_ia_{nj}-\prtl_na_{ji})
\ee
are the  associated Riemannian
Christoffel symbols.

The associated Riemannian metric tensor $a_{ij}$ has the meaning
$$
 a_{ij}=g_{ij}\bigl|_{g=0}\bigr. .
$$

The following explicit representation is obtained:
\be
g_{ij}=
\biggl[a_{ij}
+\fr g{B}\Bigl((gq^2-\fr{bS^2}q)b_ib_j-\fr bqu_iu_j+
\fr{ S^2}q(b_iu_j+b_ju_i)\Bigr)\biggr]\fr{K^2}B.
\ee
The reciprocal components $(g^{ij})=(g_{ij})^{-1}$ read
\be
g^{ij}=
\biggl[a^{ij}+\fr gq(bb^ib^j-b^iy^j-b^jy^i)+\fr g{Bq}(b+gq)y^iy^j
\biggr]\fr B{K^2}.
\ee

In many cases it is convenient to use the variables
\be
v^i~:=y^i-bb^i, \qquad v_m~:=u_m-bb_m=r_{mn}y^n\equiv r_{mn}v^n\equiv a_{mn}v^n,
\ee
where
 $r_{mn}=a_{mn}-b_mb_n.$
 Notice that
\be
r^i{}_n~:=a^{im}r_{mn}=\de^i{}_n-b^ib_n=\D{v^i}{y^n},
\ee
\ses
\be
v_ib^i=v^ib_i=0, \qquad r_{ij}b^j=r^i{}_jb^j=b_ir^i{}_j=0,\qquad u_iv^i=v_iy^i=q^2,
\ee
\ses
\be
q=\sqrt{r_{ij}v^iv^j},
\ee
and
\be
\D b{y^i}=b_i, \qquad \D q{y^i}=\fr{v_i}q, \qquad \D{(b/q)}{y^i}=\fr{2B}{NKgq^2}A_i.
\ee

{%\pgbrk}

In terms of the variables (A.13) we obtain the representations
\be
y_i=\Bigl(v_i+(b+gq)b_i\Bigr)\fr{K^2}B,
\ee
\ses
\be
g_{ij}=
\biggl[a_{ij}
+\fr g{B}\Bigl (q(b+gq)b_ib_j+q(b_iv_j+b_jv_i)-b\fr{v_iv_j}q\Bigr)\biggr]\fr{K^2}B,
\ee
and
\be
g^{ij}=
\biggl[a^{ij}+\fr gB\Bigl(-bqb^ib^j-q(b^iv^j+b^jv^i)+(b+gq)\fr{v^iv^j}q\Bigr)
\biggr]\fr B{K^2}
\ee
which are alternative to (A.11)--(A.12).

We have
\be
y_ib^i=(b+gq)\fr{K^2}B, \qquad
g_{ij}b^j=\bigl(b_i+gq\fr{y_i}{K^2}\bigr)\fr{K^2}B,
\ee
\ses
\be
g_{ij}v^j=\bigl(S^2v_i+gq^3b_i\bigr)\fr{K^2}{B^2},
\ee
\ses
\be
g^{ij}a_{ij}=\fr{NB+gq^2}{K^2},
\qquad
h_{ij}b^j=\bigl(b_i-b\fr{y_i}{K^2}\bigr)\fr{K^2}B.
\ee

By the help of the formulas (A.5) and (A.12) we find
\be
A_i=\fr {NK}2g\fr1{q}(b_i-\fr b{K^2}y_i),
\ee
or
\be
A_i=\fr {NK}2g\fr1{qB}(q^2b_i- bv_i),
\ee
and
\be
A^i=\fr N2g\fr 1{qK}
\Bigl[Bb^i-(b+gq)y^i\Bigr],
\ee
or
\be
A^i=\fr N2g\fr 1{qK}
\Bigl[q^2b^i-(b+gq)v^i\Bigr],
\ee
together with
\be
A_ib^i=\fr N2gq\fr KB, \qquad    A^ib_i=\fr N2gq\fr 1K.
\ee
These formulas are convenient to verify the contraction (A.7) and
the algebraic structure (A.8).

Since
$$
\fr{v^iv^j}q\to 0 \quad {\text{when}}\quad v^i\to 0
$$
(notice (A.16)) the components  $g_{ij}$ and $g^{ij}$ given by (A.11) and (A.12) are smooth on all the slit tangent bundle.
However, the  components of the Cartan tensor are singular at $v^i=0$, as this is apparent from the above
formulas (A.24)--(A.28) in which the {\it pole singularity} takes place at $q=0$.
 Therefore,  {\it
on the slit tangent bundle the $\cF\cF^{PD}_g $--space is smooth
 of the class $C^2$ and not of the class $C^3$}.

{%\pgbrk}

Also,
\be
A_{ij}~:= K\partial A_i/\partial y^j+l_iA_j=
%\ee and get \be
-\fr N2\fr{gb}{q}{\cal H}_{ij}+\fr{2}NA_iA_j
\ee
with
the tensor
\be
{\cal H}_{ij}=h_{ij}-\fr{A_iA_j}{A_nA^n}.
\ee

It can readily be verified that
\be
g^{ij}{\cal H}_{ij}=N-2,
\ee
\ses
\be
 g^{mn}{\cal H}_{im}{\cal
H}_{jn}={\cal H}_{ij},
\ee
\ses
\be
{\cal H}_{ij}=
\Bigl( r_{ij}-\fr1{q^2}v_iv_j\Bigr)\fr {K^2}{B}
\equiv
\eta_{ij}\fr {K^2}{B}
\quad
\text{and} \quad
{\cal H}_i{}^j~:=g^{jn}{\cal H}_{ni}=
r_i{}^j-\fr1{q^2}v_iv^j
\equiv \eta_i{}^j.
\ee
The last tensor fulfills obviously the identities
\be
{\cal H}_{ij}y^j=0,\qquad
{\cal H}_{ij}b^j=0,
\ee
which in turn entails
\be
{\cal H}_{ij}A^j=0
\ee
because $A^i$ are linear combinations of $y^i$ and $b^i$
(see (A.26)).
We also have
\be
K\Bigl(\D{{\cal H}_{ij}}{y^k}
-\D{{\cal H}_{kj}}{y^i}
\Bigr)
=l_k{\cal H}_{ij}-l_i{\cal H}_{kj}
-\fr1{A_nA^n}
\fr{Ngb}{2q}
(A_k{\cal H}_{ij}-A_i{\cal H}_{kj}
)
\ee
together with
\be
\D{(q^2{\cal H}_i^j)}{y^k}
-\D{(q^2{\cal H}_k^j)}{y^i}
=3\Bigl[
(\de_i{}^j-b_ib^j)(a_{kn}y^n-bb_k)
-
(\de_k{}^j-b_kb^j)(a_{in}y^n-bb_i)
\Bigr]
=3({\cal H}_i^jv_k-{\cal H}_k^jv_i).
\ee

The  structure (A.8) of the $\cF\cF^{PD}_g $--space Cartan tensor is such that
\be
A_k\3Aikj=
\fr1N(A_iA_j+h_{ij}A_kA^k)=\fr1N(2A_iA_j+{\cal H}_{ij}A_kA^k),
\ee
so that the tensor
\be
\tau_{ij}~:
=
A_{ij}-A_k\3Aikj
=-\fr N4\fr{g(2b+gq)}{q}{\cal H}_{ij}
\ee
obeys the identities
\be
 \tau_{ij}b^j=b_j\tau_i{}^j=\tau_{ij}A^j=0.
\ee
The tensor
\be
\tau_{ijmn}~:=K\partial A_{jmn}/\partial y^i-A_{ij}{}^hA_{hmn}-A_{im}{}^hA_{hjn}
-A_{in}{}^hA_{hjm}
+l_jA_{imn}+l_mA_{ijn}+l_nA_{ijm}
\ee
can be expressed as follows:
\be
\tau_{ijmn}=-\fr{g(2b+gq)}{4q}({\cal H}_{ij}{\cal H}_{mn}+{\cal H}_{im}{\cal H}_{jn}+{\cal H}_{in}{\cal H}_{jm}),
\ee
showing the total symmetry in all four indices and the properties
$$
\tau_{ij}= g^{mn}\tau_{ijmn}
$$
and
\be
\qquad y^i\tau_{ijmn}=0,\qquad A^i\tau_{ijmn}=0,\qquad b^i\tau_{ijmn}=0.
\ee

{%\pgbrk}

Evaluations frequently involve the vector
$m_i=(2/Ng)A_i$ which possesses the properties
$$
g^{ij}m_im_j=1, \qquad y^im_i=0.
$$
From  (A.24) it follows that
\be
m_i=K\fr1{q}(b_i-\fr b{K^2}y_i).
\ee
The equality
\be
K\D{m_i}{y^n}
=-m_nl_i
+
gm_nm_i
-\fr b{q}{\cal H}_{in}
\ee
holds.
The contravariant components $m^i$ can be taken from
(A.27):
\be
m^i=\fr 1{qK}
\Bigl[q^2b^i-(b+gq)v^i\Bigr],
\ee
entailing
\be
K\D{m^i}{y^n}=-m_nl^i
-gm^im_n
-
\fr 1{q}
(b+gq){\cal H}^i_n.
\ee

It is also valid that
\be
K\D{  {\cal H}^{kj}}{y^m}=
 -
gm_m{\cal H}^{jk}
-l^k{\cal H}^j_m  -l^j{\cal H}^k_m
+\fr bq( m^j{\cal H}^k_m
+ m^k {\cal H}^j_m)
\ee
\ses
and
\be
K\D{ \lf( {\cal H}^{kj}\fr{K^2}B\rg)}{y^m}=
\lf[
-l^k{\cal H}^j_m  -l^j{\cal H}^k_m
+\fr bq( m^j{\cal H}^k_m
+ m^k {\cal H}^j_m)
\rg]
\fr{K^2}B.
\ee

{%\pgbrk}

\ses

If we introduce the covariant derivative
$\cS$ operative in the
tangent Riemannian spaces, such that
$$
K\cS_mm^i= K\D{m^i}{y^m} + A^i{}_{mt}m^t,
\qquad
K\cS_mm_i= K\D{m_i}{y^m} - A^t{}_{mi}m_t,
$$
\ses
we obtain
$$
K\cS_mm^i=
-m_ml^i
-gm^im_m
-
\fr 1{q}
(b+gq){\cal H}^i_m
+\fr1Nm^t
 \lf[
A^ih_{tm}  +A_th^i_{m}  +A_mh^i_{t}
-\fr 4{N^2g^2}
A^iA_tA_m
\rg],
$$
\ses
so that
\be
K\cS_mm^i=
-m_ml^i
-
\fr 1{q}
\lf(b+\fr12gq\rg){\cal H}^i_m.
\ee
\ses
Also, with the definition
\be
K\cS_m{  {\cal H}^{kj}}= K\D{{\cal H}^{kj}}{y^m} +
 A^{kt}{}_{m}{\cal H}^{tj}
+ A^{jt}{}_{m}{\cal H}^{tk},
\ee
\ses
we get
\be
K\cS_m{  {\cal H}^{kj}}=
-l^k{\cal H}^j_m  -l^j{\cal H}^k_m
+\fr 1q
\lf(b+\fr12gq\rg)
( m^j{\cal H}^k_m + m^k {\cal H}^j_m).
\ee

{%\pgbrk}

\ses

\ses

The equality
$
\partial {K^2} /\partial g={\bar M}   K^2
$
holds with
\be
{\bar M}=
 - \fr1{h^3}f+
\fr12\fr{G}{hB}   q^2+  \frac1{h^2B} b q.
\ee
In obtaining this formula we have used the derivatives
\be
  \D hg= -\fr14 G, \quad
    \D Gg= \fr1{h^3}, \quad \D{\lf(\fr Gh\rg)} g =\fr1{h^4} \lf(1+\fr{g^2}4 \rg),
\quad
\D fg= -\fr1{2h} +  \fr{b}B\Bigl(\fr14G  q +   \fr1{2h}b\Bigr).
\ee
\ses
Therefore,
\be
\partial^*_nK=
{\bar M}   K^2\D g{x^n}.
\ee
\ses
It follows that
\be
\D {\bar M}{y^h}=\fr{2b^4}{B^2}\D{\fr bq}{y^m}=\fr{4  q^2X}{gBK} A_h
\ee
and
\be
\partial^*_nl_m= \D{ \lf(\partial^*_nK\rg)}{y^m}.
\ee

{%\pgbrk}

The substitution
$$
a^{nt}=
\fr {K^2}B{\cal H}^{nt}
+b^nb^t+\fr1{q^2}v^nv^t
$$
transforms the tensor
 (6.71) to
$$
 M^n{}_{ij}=  \Bigl((1-h)b  +\fr12gq    \Bigr)\fr {K^2}B{\cal H}^{nt}
\fr1h
b_la_t{}^l{}_{ij}
-
a_t{}^n{}_{ij}
 y^t
 $$

\ses

$$
+
\Biggl[
(1-h)b
\lf(b^nb^t+\fr1{q^2}v^nv^t\rg)
+\fr12gq
\lf(b^nb^t+\fr1{q^2}v^nv^t\rg)
+
\lf(
\fr g{2q}v^n
-(1-h)b^n
\rg)  y^t
\Biggr]
\fr1h
b_la_t{}^l{}_{ij}
$$

\ses

\ses

\ses

$$
=  \Bigl((1-h)b  +\fr12gq    \Bigr)\fr {K^2}B{\cal H}^{nt}
\fr1h
b_la_t{}^l{}_{ij}
+
\Biggl[
(1-h)
\lf(\fr b{q^2}v^n-b^n\rg)
+
\fr g{q}v^n
\Biggr]
\fr1h
y^tb_la_t{}^l{}_{ij}
-
a_t{}^n{}_{ij}
 y^t.
$$

{%\pgbrk}

Applying

$$
y_i=(u_i+gqb_i)\fr{K^2}B,
\qquad
b^n=\fr1B\lf[Kqm^n+(b+gq)y^n\rg]
$$
together with
\ses\\
$$
a_t{}^n{}_{ij}  y^t= h^n_pa_t{}^p{}_{ij}y^t
+ l^nl_la_t{}^l{}_{ij}y^t, \qquad
l_la_t{}^l{}_{ij}y^t=gq\fr{K}B
b_la_t{}^l{}_{ij}y^t
$$

{%\pgbrk}

\nin
leads to
$$
 M^n{}_{ij}=  \Bigl((1-h)b  +\fr12gq    \Bigr)\fr {K^2}B{\cal H}^{nt}
\fr1h
b_la_t{}^l{}_{ij}
$$

\ses

$$
+
\Biggl[
(1-h)
\lf(\fr b{q^2}y^n
-\fr {S^2}{q^2}  \fr1B\lf[Kqm^n+(b+gq)y^n\rg]
\rg)
+
\fr g{q}y^n
-
\fr g{q}b      \fr1B\lf[Kqm^n+(b+gq)y^n\rg]
\Biggr]
\fr1h
y^tb_la_t{}^l{}_{ij}
$$

\ses

$$
-
 h^n_pa_t{}^p{}_{ij}y^t
-
 l^n
gq\fr{K}B
b_la_t{}^l{}_{ij}y^t,
$$
\ses
or
$$
 M^n{}_{ij}=  \Bigl((1-h)b  +\fr12gq    \Bigr)\fr {K^2}B{\cal H}^{nt}
\fr1h
b_la_t{}^l{}_{ij}
-
 \lf[
(1-h)
\fr {B}{q}    +hgb\rg]
\fr1B
 Km^n
 \fr1h
y^tb_la_t{}^l{}_{ij}
$$

\ses

\ses

\ses

$$
+
\Biggl[
(1-h)
\lf(\fr b{q^2}y^n
-\fr {B-gbq}{q^2}  \fr1B(b+gq)y^n
\rg)
+
 gqy^n\fr1B
\Biggr]
\fr1h
y^tb_la_t{}^l{}_{ij}
$$

\ses

$$
-
 h^n_pa_t{}^p{}_{ij}y^t
-
 l^n
gq\fr{K}B
b_la_t{}^l{}_{ij}y^t,
$$
which is
\ses\\
$$
 M^n{}_{ij}=  \Bigl((1-h)b  +\fr12gq    \Bigr)\fr {K^2}B{\cal H}^{nt}
\fr1h
b_la_t{}^l{}_{ij}
-
 \lf[
(1-h)
\fr {B}{q}    +hgb\rg]
\fr1B
Km^n
 \fr1h
y^tb_la_t{}^l{}_{ij}
$$

\ses

\ses

\ses

$$
+
\Biggl[
(1-h)
\lf(\fr b{q^2}y^n
-\fr 1{q^2}  (b+gq)y^n
+\fr{g}{q}  \fr1B(B-q^2)y^n
\rg)
+
 gqy^n\fr1B
\Biggr]
\fr1h
y^tb_la_t{}^l{}_{ij}
$$

\ses

$$
-
 h^n_pa_t{}^p{}_{ij}y^t
-
 l^n
gq\fr{K}B
b_la_t{}^l{}_{ij}y^t,
$$
\ses
so that
$$
 M^n{}_{ij}=  \Bigl((1-h)b  +\fr12gq    \Bigr)\fr {K^2}B{\cal H}^{nt}
\fr1h
b_la_t{}^l{}_{ij}
- \lf[
\fr {B}{q} -h\fr{S^2}q\rg]
\fr1B
Km^n
 \fr1h
y^tb_la_t{}^l{}_{ij}
$$

\ses

$$
-
{\cal H}^n_la_t{}^l{}_{ij}y^t
-
m^nm_la_t{}^l{}_{ij}y^t.
 $$
\ses

{%\pgbrk}

The result is
\be
 M^n{}_{ij}=  \Bigl((1-h)b  +\fr12gq    \Bigr)\fr {K^2}B{\cal H}^{nt}
\fr1h
b_la_t{}^l{}_{ij}
-
{\cal H}^n_la_t{}^l{}_{ij}y^t
- \fr {1}{q}
Km^n
 \fr1h
y^tb_la_t{}^l{}_{ij}.
\ee
\ses
Lowering here the index $n$ leads to the representation (6.72).

It is also possible to write
$$
\fr B{K^2} M_{nij}=
  \Bigl((1-h)b  +\fr12gq    \Bigr)
\fr1h
b_la_n{}^l{}_{ij}
-\fr1{q^2}v_ny^t
  \Bigl((1-h)b  +\fr12gq    \Bigr)
\fr1h
b_la_t{}^l{}_{ij}
$$

\ses

$$
-
a_{tnij}y^t
+b_nb_l
a_{tnij}y^t
 -\fr b{q^2}v_nb_l a_t{}^l{}_{ij}y^t
-
\fr B{K^2}
\fr {1}{q}
K
\fr{K}{qB}(q^2b_n-bv_n)
 \fr1h
y^tb_la_t{}^l{}_{ij},
$$
which can be simplified to read
\ses\\
\be
\fr B{K^2} M_{nij}=
\Bigl((1-h)b  +\fr12gq    \Bigr)
\fr1h
b_la_n{}^l{}_{ij}
-
\lf(
\fr g{2q}v_n
+
(1-h)b_n
\rg)
\fr1h
 y^t
 b_la_t{}^l{}_{ij}
-
a_{tnij}
 y^t.
\ee

{%\pgbrk}

Also, considering
the contraction
$$
\fr B{K^2}M^{nij} M_{nij}=
$$

\ses

\ses

$$
\Bigl((1-h)b  +\fr12gq    \Bigr)
\fr1h
b_la^{nlij}
\Biggl[
\Bigl((1-h)b  +\fr12gq    \Bigr)
\fr1h
b_ha_n{}^h{}_{ij}
-
\fr g{2q}u_n
\fr1h
 y^t
 b_ha_t{}^h{}_{ij}
-
a_{tnij}
 y^t
\Biggr]
$$

\ses

\ses

$$
+
\fr g{2q}y^n
 y^t
\fr1h
b_la_t{}^{lij}
\Biggl[
\Bigl((1-h)b  +\fr12gq    \Bigr)
\fr1h
b_la_n{}^l{}_{ij}
-
\lf(
\fr g{2q}v_n
+
(1-h)b_n
\rg)
\fr1h
 y^t
 b_la_t{}^l{}_{ij}
\Biggr]
$$

\ses

\ses

$$
+
\lf(
\fr g{2q}b+1-h\rg)
  y^t
\fr1h
b_la_t{}^{lij}
\Biggl[
(1-h)
\fr1h
 y^t
 b_la_t{}^l{}_{ij}
+
b^na_{tnij}
 y^t
\Biggr]
$$

\ses

\ses

$$
-
a_h{}^{nij}
 y^h
\Biggl[
\Bigl((1-h)b  +\fr12gq    \Bigr)
\fr1h
b_la_n{}^l{}_{ij}
+\lf(\fr{gb}{2q}-(1-h) \rg)\fr1h b_n
  y^t
 b_la_t{}^l{}_{ij}
 -
a_{tnij}
 y^t
\Biggr]
$$

\ses

\ses

\ses

$$
=
\fr1{h^2}
\Bigl((1-h)b  +\fr g2q    \Bigr)^2
b_la^{nlij}
b_ha_n{}^h{}_{ij}
-
\fr 2h
\Bigl((1-h)b  +\fr g2q    \Bigr)
y^ha_h{}^{nij}
b_la_n{}^l{}_{ij}
+
y^ha_h{}^{nij}
a_{tnij}
 y^t
$$
leads to (6.82).

{%\pgbrk}

We can start also with  (A.58), observing that
$$
 M^{nij} M_{nij}=
  \Biggl[
 \Bigl((1-h)b  +\fr12gq    \Bigr)\fr {K^2}B{\cal H}^{nh}
\fr1h
b_pa_h{}^{pij}
-
{\cal H}^n_p a_h{}^{pij}y^h
- \fr {1}{q}
Km^n
 \fr1h
y^hb_pa_h{}^{pij}
\Biggr]
\times
$$

\ses

$$
 \Biggl[
 \Bigl((1-h)b  +\fr12gq    \Bigr)\fr {K^2}B{\cal H}_n^t
\fr1h
b_la_t{}^l{}_{ij}
-
{\cal H}_{nl}a_t{}^l{}_{ij}y^t
- \fr {1}{q}
Km_n
 \fr1h
y^tb_la_t{}^l{}_{ij}
\Biggr]
$$

\ses

 \ses

 $$
=
 \Bigl((1-h)b  +\fr12gq    \Bigr)\fr {K^2}B
\fr1h
b_pa_h{}^{pij}
 \Biggl[
 \Bigl((1-h)b  +\fr12gq    \Bigr)\fr {K^2}B{\cal H}^{th}
\fr1h
b_la_t{}^l{}_{ij}
-
{\cal H}^h_l a_t{}^l{}_{ij}y^t
\Biggr]
$$

\ses

\ses

 \ses

 $$
-
 a_h{}^{pij}y^h
 \Biggl[
 \Bigl((1-h)b  +\fr12gq    \Bigr)\fr {K^2}B{\cal H}_p^t
\fr1h
b_la_t{}^l{}_{ij}
-
{\cal H}_{pl}a_t{}^l{}_{ij}y^t
\Biggr]
+ \fr {1}{q^2}   K^2
 \fr1{h^2}
y^hb_pa_h{}^{pij}
y^tb_la_t{}^l{}_{ij}
$$

\ses

\ses

\ses

 \ses

 $$
=
 \fr1{h^2}
 \Bigl((1-h)b  +\fr12gq    \Bigr)^2\fr {K^2}B
b_pa_h{}^{pij}
\lf( a^{th}-\fr1{q^2}y^ty^h\rg)
b_la_t{}^l{}_{ij}
$$

 \ses

\ses

 $$
-
2\fr1h
 \Bigl((1-h)b  +\fr12gq    \Bigr)\fr {K^2}B
b_pa_h{}^{pij}
\lf(\de^h_l+\fr1{q^2}by^hb_l\rg)
 a_t{}^l{}_{ij}y^t
$$

\ses

\ses

 $$
+
 a_h{}^{pij}y^h
\fr {K^2}B
\lf(a_{pl}-b_pb_l-\fr1{q^2}b^2b_pb_l \rg)
a_t{}^l{}_{ij}y^t
+ \fr {K^2}{h^2q^2}
y^hb_pa_h{}^{pij}
y^tb_la_t{}^l{}_{ij}
\fr{ \Bigl(b  +\fr12gq    \Bigr)^2+h^2q^2}B
$$

 \ses

 \ses

 \ses

$$
 =
\!\Biggl[
\fr1{h^2}
\Bigl((1\!-\!h)b  +\fr g2q    \Bigr)^2
b_la^{nlij}
b_ha_n{}^h{}_{ij}
-
\fr 2h\!
\Bigl((1-h)b  +\fr g2q    \Bigr)
y^ha_h{}^{nij}
b_la_n{}^l{}_{ij}
+
y^ha_h{}^{nij}
a_{tnij}
 y^t
\Biggr]\!
\fr {K^2}B
\!+\!E
$$
\ses
with
$$
E=
-
\lf[
 \fr1{h^2}
 \Bigl((1\!-\!h)b  +\fr12gq    \Bigr)^2\fr {K^2}B
b_pa_h{}^{pij}
\fr1{q^2}
+
\fr2h
 \Bigl((1-h)b  +\fr12gq    \Bigr)\fr {K^2}B
b_pa_h{}^{pij}
\fr b{q^2}
\rg]
y^hb_l
 a_t{}^l{}_{ij}y^t
$$

\ses

\ses

 $$
+
 a_h{}^{pij}y^h
\fr {K^2}B
\lf(-\fr1{q^2}b^2b_pb_l \rg)
a_t{}^l{}_{ij}y^t
+ \fr {1}{q^2}   K^2
 \fr1{h^2}
y^hb_pa_h{}^{pij}
y^tb_la_t{}^l{}_{ij}
\fr{ \Bigl(b  +\fr12gq    \Bigr)^2}B
=0.
$$
Thus we  obtain (6.82) from new standpoint.

 {%\pgbrk}

The contraction
 can be written in the concise form
$$
 M^{nij} M_{nij}
 =
\fr {K^2}{h^2B}
\lf[ hv^l+\lf(b  +\fr g2q    \rg) b^l\rg]
a_l{}^{nij}
\lf[ hv^h+\lf(b  +\fr g2q    \rg)b^h\rg]
a_{hnij}.
$$
With (6.26) and (6.39), we have
$$
\zeta^i=\lf[hv^i  +   ( b+\fr12gq)b^i    \rg] \fr S{\sqrt B},
$$
\ses
obtaining  the simple representation
$$
 M^{nij} M_{nij}=
\ka^2
\zeta^l
a_l{}^{nij}
\zeta^h
a_{hnij}
$$
which is equivalent to (6.83);
$\ka^2= {K^2}/h^2S^2$ in accordance with (3.6).

\ses

\ses

Also, it is possible to  get
$$
E_k{}^n{}_{ij}+\fr1Kl_kM^n{}_{ij}
-\fr1K
l^nM_{kij}
=
\Biggl[
-
\Bigl((1-h)q
-\fr g2
 (b+gq)
  \Bigr)
   b^l
-
h
K
m^l
\Biggr]
m_k
\fr {K}B{\cal H}^{nt}
\fr1h
a_{tlij}
$$

\ses
\ses

$$
+
\lf(\fr g{2}b+q\rg)m^n
K\fr1B
{\cal H}_k^t
\fr1h
b_la_t{}^l{}_{ij}
-
K\fr1Bqm^n
{\cal H}_k^t
b_la_t{}^l{}_{ij}
$$
\ses

\be
-
 \fr g{2q}
{\cal H}^n_k
\fr1h
b_la_t{}^l{}_{ij}y^t
+{\cal H}^n_l  {\cal H}_k^ta_t{}^l{}_{ij}
+m^nm_l  {\cal H}_k^ta_t{}^l{}_{ij}
-
    g  m^n
m_k \fr1{q}
\fr1h
b_la_t{}^l{}_{ij}y^t
\ee
\ses
and
\be
BE_{knij}E^{knij}=
\fr B{K^2} (KE_{knij}+l_kM_{nij}-l_nM_{kij})
(KE^{knij}+l^kM^{nij}-l^nM^{kij})
+
2
\fr B{K^2}
M_{kij}M^{kij},
\ee
\ses
together with
$$
E_{knij}E^{knij}=  a_{knij} a^{knij}
+
g^2
\fr{1}{h^2q^2}
\Bigl(
(N-2)\fr{1}{4}
+
1
\Bigr)
b^t   y^la_{tl}{}^{ij}    y^h    b^s    a_{shij}
$$

\ses

\ses

\be
+
\fr{g^2}{B}
\Biggl[
\Biggl( b - \fr1h\lf(b+\fr12gq\rg)  \Biggr)
b_ha^{nhij}
-
a^n{}_h{}^{ij}
 y^h
\Biggr]
\Biggl[
\Biggl( b - \fr1h\lf(b+\fr12gq\rg)  \Biggr)
b_la_n{}^l_{ij}
-
a_{ntij}
 y^t
\Biggr]
\ee
(see Appendix C).

{%\pgbrk}

Let us verify the formulas (A.45) and (A.47).

\ses

Upon differentiating (A.44) we directly obtain the equality
$$
\D{m_i}{y^n}=\fr1Kl_nm_i
-\fr1{q^2}v_nm_i
-\fr1qb_nl_i
+
\fr b{qK}l_nl_i
-\fr b{qK}h_{in},
$$
in which
the substitutions
\be
b_n=\fr qKm_n+\fr b{K}l_n,
\ee
\ses
\be
bv_n=q^2b_n-\fr1K qBm_n
=
q^2\fr qKm_n+q^2\fr b{K}l_n -\fr1K qBm_n,
\ee
and
\be
v_n
=
q^2\fr 1{K}l_n -\fr1K q(b+gq)m_n
\ee
can conveniently be used.
We obtain
$$
\D{m_i}{y^n}
=
\fr1Kl_nm_i
-\fr1{q^2}
\lf[q^2\fr 1{K}l_n -\fr1K q(b+gq)m_n\rg]
m_i
-\fr1q
\lf[\fr qKm_n+\fr b{K}l_n\rg]
l_i
+
\fr b{qK}l_nl_i
-\fr b{qK}h_{in}
$$

\ses

\ses

$$
=
\fr1q(b+gq)m_n
m_i
-m_nl_i
-\fr b{q}h_{in},
$$
so that the formula (A.45) is valid.

{%\pgbrk}

Also, we differentiate (A.46) and apply
$
b^i=\lf(Kqm^i+(b+gq)y^i\rg)/B
$
together with
$
v^i=q\lf(-Kbm^i+qy^i\rg)/ B.
$
We obtain
$$
K\D{m^i}{y^n}=
-l_nm^i -\fr1{q^2}v_nm^i
+\fr 2{q}v_n
b^i
-\fr 1{q}
\lf(b_n+\fr gqv_n\rg)v^i
-
\fr 1{q}
(b+gq)\lf({\cal H}^i_n+\fr1{q^2}v^iv_n\rg)
$$

\ses

\ses

$$
=
-l_nm^i -\fr K{q^2}v_nm^i
+\fr 2{q}v_n
\fr KB\lf[qm^i+(b+gq)l^i\rg]
-\fr 1{q}
\lf(
\fr qKm_n+\fr b{K}l_n+\fr gqv_n\rg)v^i
$$

\ses

$$
-
\fr 1{q}
(b+gq)\fr1{q^2}v^iv_n
-
\fr 1{q}
(b+gq){\cal H}^i_n
$$

\ses

\ses

\ses

$$
=
-l_nm^i -\fr1{q^2}
\lf[q^2l_n - q(b+gq)m_n\rg]
m^i
$$

\ses

$$
+
\fr 2{q}
q^2\fr 1{K}l_n
\fr KB\lf[qm^i+(b+gq)l^i\rg]
-
\fr 2{q}
\fr1K q(b+gq)m_n
\fr KB\lf[qm^i+(b+gq)l^i\rg]
-\fr 1{q}
\lf(
\fr qKm_n+\fr b{K}l_n\rg)v^i
$$

\ses

$$
-
\fr 1{q}
(b+2gq)\fr1{q^2}v^i
q^2\fr 1{K}l_n
+
\fr 1{q}
(b+2gq)\fr1{q^2}v^i
\fr1K q(b+gq)m_n
-
\fr 1{q}
(b+gq){\cal H}^i_n
$$

\ses

\ses

\ses

$$
=
-2l_nm^i
+\fr1{q^2}
 q(b+gq)m_n
m^i
$$

\ses

$$
+
\fr 2{q}
q^2l_n
\fr 1B\lf[qm^i+(b+gq)l^i\rg]
-
\fr 2{q}
\fr1K q(b+gq)m_n
\fr KB\lf[qm^i+(b+gq)l^i\rg]
-\fr 1K m_nv^i
$$

\ses

$$
-
\fr 2{q}
(b+gq)v^i
\fr 1{K}l_n
+
\fr 1{q}
(b+2gq)\fr1{q^2}v^i
\fr1K q(b+gq)m_n
-
\fr 1{q}
(b+gq){\cal H}^i_n.
$$
\ses
Making here natural reductions
leads to
\ses\\
$$
K\D{m^i}{y^n}=
-2l_nm^i
+\fr1{q}
(b+gq)m_n
m^i
+
\fr 2{B}
q^2l_n
m^i
-
\fr 2B
(b+gq)m_n
\lf[qm^i+(b+gq)l^i\rg]
$$

\ses

$$
-\fr 1K m_nv^i
+
\fr 2{q}
(b+gq)
\fr qBbm^i
l_n
+
\fr 1{q}
(b+2gq)\fr1{q^2}v^i
\fr1K q(b+gq)m_n
-
\fr 1{q}
(b+gq){\cal H}^i_n
$$

\ses

\ses

$$
=-m_nl^i
+
\fr1{q}
(b+gq)m_n
m^i
-
\fr1q(b+2gq)m_n
m^i
-
\fr 1{q}
(b+gq){\cal H}^i_n.
$$
In this way the validity of (A.47) is straightforwardly verified.

{%\pgbrk}

\ses

Finally,
considering the equality
$$
K\D {  {\cal H}^{kj}}{y^m}=
K\D { (h^{kj}-m^km^j)}{y^m}=
-2A^{kj}{}_m-(l^kh^j_m+l_jh^k_m)
$$

\ses

$$
+ m^j \lf[m_ml^k+gm^km_m+\fr 1{q}(b+gq){\cal H}^k_m\rg]
+ m^k \lf[m_ml^j+gm^jm_m+\fr 1{q}(b+gq){\cal H}^j_m\rg],
 $$
which is simplified to read
\ses\\
$$
K\D{  {\cal H}^{kj}}{y^m}=
 -
gm_m{\cal H}^{jk}  -gm^j{\cal H}^k_m  -gm^k{\cal H}^j_m-2gm^km^jm_m
-l^k({\cal H}^j_m+m^jm_m)  -l^j({\cal H}^k_m+m^km_m)
$$

\ses

$$
+ m^j \lf[m_ml^k+gm^km_m+\fr 1{q}(b+gq){\cal H}^k_m\rg]
+ m^k \lf[m_ml^j+gm^jm_m+\fr 1{q}(b+gq){\cal H}^j_m\rg],
 $$
we can readily conclude that the formulas (A.48) and (A.49) are true.

{%\pgbrk}

\ses

\ses

\setcounter{equation}{0}

\nin
{ \bf Appendix B: Representations for connection coefficients}

\ses

\ses

With (6.53) and (6.54)
we evaluate the sum
$$
N^k_n+l^k\D K{x^n}=  N^k_n-l^kN^m_nl_m
=
\fr{gq}By^ky^j   \nabla_nb_j
$$

\ses

$$
+
\Biggl[\lf(b-\fr1h\lf(b+\fr12gq\rg)\rg)
\eta^{kj}
+
\lf(
\fr1{q^2} v^k \lf(b   -  \fr1h(b+gq)\rg)
+\lf(\fr1h-1\rg)b^k
\rg)
y^j
\Biggr]
 \nabla_nb_j
-    h^k_t a^t{}_{nj}y^j
$$

\ses

\ses

$$
=
\Biggl[\lf(b-\fr1h\lf(b+\fr12gq\rg)\rg)
\eta^{kj}
+\fr1{hq^2}[-(b+gq)v^k+q^2b^k]
y^j
 +
\lf(
\fr b{q^2} v^k
-b^k
\rg)
y^j
\Biggr]
 \nabla_nb_j
-    h^k_t a^t{}_{nj}y^j
$$

\ses

$$
+
\fr{gq}By^ky^j   \nabla_nb_j,
$$
\ses
coming to the representation
\be
N^k_n
=-l^k\D K{x^n}
+
\Biggl[\lf(b-\fr1h\lf(b+\fr12gq\rg)\rg)
\eta^{kj}
+\lf(\fr1{hq} -\fr{b^2+q^2}{qB}\rg)Km^k
y^j
\Biggr]
 \nabla_nb_j
-    h^k_t a^t{}_{nj}y^j
\ee
which is obviously equivalent to (6.62).

{%\pgbrk}

Let us differentiate  (6.62) with respect to $y^m$:
\ses
 $$
 N^k_{nm}
=-\fr1K h^k_m\D K{x^n}  -l^k\D {l_m}{x^n}
+
\lf(b_m-\fr1h\lf(b_m+\fr1{2q}gv_m\rg)\rg)
{\cal H}^{kj}\fr{K^2}B
 \nabla_nb_j
 $$

\ses
\ses

$$
+
\lf(b-\fr1h\lf(b+\fr12gq\rg)\rg)
\lf(
-l^k{\cal H}^j_m  -l^j{\cal H}^k_m
+\fr bq( m^j{\cal H}^k_m
+ m^k {\cal H}^j_m)
\rg)
\fr{K}B
 \nabla_nb_j
$$

\ses

$$
+
Z_m
m^k
y^j
 \nabla_nb_j
+
\lf(\fr1{hq}-\fr1q +\fr{gb}{B}\rg)
\lf[
-m_ml^k
-gm^km_m
-
\fr 1{q}
(b+gq){\cal H}^k_m\rg]
y^j
 \nabla_nb_j
$$

\ses

\ses

$$
+
\lf(\fr1{hq}-\fr1q +\fr{gb}{B}\rg)Km^k
 \nabla_nb_m
-    h^k_t a^t{}_{nm}
+\fr1K(l^kh_{tm}+l_th^k_m) a^t{}_{nj}y^j,
$$
\ses
where
$$
Z_m=K
\lf[
\fr{v_m}{q^3}
\fr{h-1}h
+\fr{gb_m}{B}
-\fr{gb}{B^2}
\lf(
2bb_m+gqb_m+gb\fr1qv_m+2v_m
\rg)
\rg]
+l_m\lf(\fr1{hq}-\fr1q +\fr{gb}{B}\rg).
$$

Apply
$$
Kb_m= qm_m+ bl_m,
\qquad
Kv_m
=
q^2l_m - q(b+gq)m_m,
$$
\ses
which yields
$$
Z_m=
-\fr1{hq^3}(q^2l_m - q(b+gq)m_m)+\fr1{q^3}(q^2l_m - q(b+gq)m_m)
 + \fr{g(qm_m+ bl_m)}{B}
$$

\ses

$$
-\fr{gb}{B^2}
\lf(
b(2qm_m+ 2bl_m)+gq(qm_m+ bl_m)+gb(ql_m - (b+gq)m_m)+2(q^2l_m - q(b+gq)m_m)
\rg)
$$

\ses

$$
+l_m\lf(\fr1{hq}-\fr1q +\fr{gb}{B}\rg),
$$
\ses
or
$$
Z_m=
-\fr1{hq^3}[q^2l_m - q(b+gq)m_m]+\fr1{q^3}[q^2l_m - q(b+gq)m_m]
 + g\fr{(q+gb)m_m- bl_m}{B}
$$

\ses

$$
+l_m\lf(\fr1{hq}-\fr1q +\fr{gb}{B}\rg).
$$
\ses
So we   have
\be
Z_m=
\lf[
\fr1{hq^2}(b+gq) -\fr1{q^2}(b+gq)
 + g\fr{q+gb}{B}
 \rg]
 m_m.
\ee

{%\pgbrk}

Thus we can write
\ses
 $$
 N^k_{nm}
=
-\fr1K h^k_m
\lf[
gq\fr{K^2}Bl^j   \nabla_nb_j
+
y_k  a^k{}_{nj}l^j
\rg]
  -l^k\D {l_m}{x^n}
$$

\ses

\ses

$$
+
\lf(
qm_m+ bl_m-\fr1h\lf(qm_m+ bl_m+\fr1{2q}g
(q^2l_m - q(b+gq)m_m)
\rg)\rg)
{\cal H}^{kj}\fr{K}B
 \nabla_nb_j
 $$

\ses
\ses

$$
+
\lf(b-\fr1h\lf(b+\fr12gq\rg)\rg)
\lf(
-l^k{\cal H}^j_m  -l^j{\cal H}^k_m
+\fr bq( m^j{\cal H}^k_m
+ m^k {\cal H}^j_m)
\rg)
\fr{K}B
 \nabla_nb_j
$$

\ses

$$
+
\lf[
\fr1{hq^2}(b+gq) -\fr1{q^2}(b+gq)
 + g\fr{q+gb}{B}
 \rg]
 m_m
 m^k
y^j
 \nabla_nb_j
$$

\ses

\ses

$$
+
\lf(\fr1{hq}-\fr1q +\fr{gb}{B}\rg)
\lf[
-m_ml^k
-gm^km_m
-
\fr 1{q}
(b+gq){\cal H}^k_m\rg]
y^j
 \nabla_nb_j
$$

\ses

\ses

$$
+
\lf(\fr1{hq}-\fr1q +\fr{gb}{B}\rg)Km^k
 \nabla_nb_m
-    h^k_t a^t{}_{nm}
+\fr1K(l^kh_{tm}+l_th^k_m) a^t{}_{nj}y^j.
$$
\ses\\
Cancelling similar terms leads to
 $$
 N^k_{nm}
=-g\fr qK {\cal H}^k_m
\fr{K^2}Bl^j   \nabla_nb_j
  -l^k\lf(\D {l_m}{x^n}  -h_{tm} a^t{}_{nj}l^j\rg)
$$

\ses

\ses

$$
+
\fr1K\lf[q-\fr1hq+\fr1{2h}g(b+gq)\rg]
{\cal H}^{kj}m_m\fr{K^2}B
 \nabla_nb_j
 $$

\ses
\ses

$$
+
\lf(b-\fr1h\lf(b+\fr12gq\rg)\rg)
\fr{K}B
\lf(
{\cal H}^{kj}l_m-l^k{\cal H}^j_m  -l^j{\cal H}^k_m
+\fr bq( m^j{\cal H}^k_m
+ m^k {\cal H}^j_m)
\rg)
 \nabla_nb_j
$$

\ses

$$
+
K\lf[
\fr1{hq^2}b -\fr1{q^2}b
 \rg]
 m_m
 m^k
l^j
 \nabla_nb_j
+
K\lf(\fr1{hq}-\fr1q +\fr{gb}{B}\rg)
\lf[
-m_ml^k
-
\fr 1{q}
(b+gq){\cal H}^k_m\rg]
l^j
 \nabla_nb_j
$$

\ses

\ses

\be
+
\lf(\fr1{hq}-\fr1q +\fr{gb}{B}\rg)Km^k
 \nabla_nb_m
-    h^k_t a^t{}_{nm}.
\ee

{%\pgbrk}

The equality
$$
\D K{x^n}=
gq\fr{K^2}Bl^j   \nabla_nb_j
+
y_t  a^t{}_{nj}l^j
$$
(we assume $g=const$)
ensuing from (6.54) and (6.59)
can readily
be differentiated with respect to $y^m$, yielding
$$
\D {l_m}{x^n} -h_{tm} a^t{}_{nj}l^j=
 g
[ql_m - (b+gq)m_m]
\fr{K}Bl^j   \nabla_nb_j
+
  \fr{gqK}Bh^j_m   \nabla_nb_j
+
g^2\fr qK  m_m\fr{K^2}B
l^j   \nabla_nb_j
+
l_t  a^t{}_{nm},
$$
or
\ses\\
\be
\D {l_m}{x^n} -h_{tm} a^t{}_{nj}l^j=
- gbm_m
\fr{K}Bl^j   \nabla_nb_j
+
g\fr qK  \fr{K^2}B   \nabla_nb_m
+
l_t  a^t{}_{nm}.
\ee
The representation (B.3) takes on the form
\ses\\
 $$
 N^k_{nm}
=-g\fr qK {\cal H}^k_m
\fr{K^2}Bl^j   \nabla_nb_j
-l^k
\lf[
- gbm_m
\fr{K}Bl^j   \nabla_nb_j
+
g\fr qK  \fr{K^2}B   \nabla_nb_m
\rg]
$$

\ses

\ses

$$
+
\fr1K\lf[q-\fr1hq+\fr1{2h}g(b+gq)\rg]
{\cal H}^{kj}m_m\fr{K^2}B
 \nabla_nb_j
 $$

\ses
\ses

$$
+
\lf(b-\fr1h\lf(b+\fr12gq\rg)\rg)
\fr{K}B
\lf(
{\cal H}^{kj}l_m-l^k{\cal H}^j_m  -l^j{\cal H}^k_m
+\fr bq( m^j{\cal H}^k_m
+ m^k {\cal H}^j_m)
\rg)
 \nabla_nb_j
$$

\ses

$$
+
K\lf[
\fr1{hq^2}b -\fr1{q^2}b
 \rg]
 m_m
 m^k
l^j
 \nabla_nb_j
-
K\lf(\fr1{hq}-\fr1q +\fr{gb}{B}\rg)
m_ml^k
l^j
 \nabla_nb_j
$$

\ses

\ses

$$
-
K\lf[\lf(\fr1{hq}-\fr1q\rg)(b+gq) +\fr{g(B-q^2)}{B}\rg]
\fr 1{q}
{\cal H}^k_m
l^j
 \nabla_nb_j
+
\lf(\fr1{hq}-\fr1q +\fr{gb}{B}\rg)Km^k
 \nabla_nb_m
-
a^k{}_{nm},
 $$
\ses\\
or
\ses\\
 $$
 N^k_{nm}
=
-l^k
g\fr qK  \fr{K^2}B   \nabla_nb_m
+
\fr1K\lf[q-\fr1hq+\fr1{2h}g(b+gq)\rg]
{\cal H}^{kj}m_m\fr{K^2}B
 \nabla_nb_j
 $$

\ses

\ses

$$
+
\lf(b-\fr1h\lf(b+\fr12gq\rg)\rg)
\fr{K}B
\lf(
{\cal H}^{kj}l_m-l^k{\cal H}^j_m  -l^j{\cal H}^k_m
+\fr bq( m^j{\cal H}^k_m
+ m^k {\cal H}^j_m)
\rg)
 \nabla_nb_j
$$

\ses

$$
+
K\lf[
\fr1{hq^2}b -\fr1{q^2}b
 \rg]
 m_m
 m^k
l^j
 \nabla_nb_j
-
\fr Kq\lf(\fr1{h}-1 \rg)
m_ml^k
l^j
 \nabla_nb_j
$$

\ses

\ses

\be
-
K
\lf(\fr1{hq}(b+gq)-\fr bq\rg)
\fr 1{q}
{\cal H}^k_m
l^j
 \nabla_nb_j
+
\lf(\fr1{hq}-\fr1q +\fr{gb}{B}\rg)Km^k
 \nabla_nb_m
-
a^k{}_{nm}.
\ee

{%\pgbrk}

By the help of the equalities
$$
B\fr b{q^2}b^i=K\fr bqm^i+\fr{B-q^2}{q^2}y^i,
\qquad
\fr bqm^j-l^j=  \fr BK\fr b{q^2}b^j-\fr B{q^2}l^j
$$
we come to
the representation
 $$
 N^k_{nm}
=
-l^k
g\fr qK  \fr{K^2}B   \nabla_nb_m
+
\fr1K\lf[q-\fr1hq+\fr1{2h}g(b+gq)\rg]
{\cal H}^{kj}m_m\fr{K^2}B
 \nabla_nb_j
 $$

\ses

\ses

$$
+
\lf(b-\fr1h\lf(b+\fr12gq\rg)\rg)
\fr{K}B
\lf(
{\cal H}^{kj}l_m-l^k{\cal H}^j_m
+\fr bq m^k {\cal H}^j_m
\rg)
 \nabla_nb_j
+
K\fr b{q^2}
\fr{1-h}h
 m_m
 m^k
l^j
 \nabla_nb_j
$$

\ses

\ses

$$
-
\fr Kq\lf(\fr1{h}-1 \rg)
m_ml^k
l^j
 \nabla_nb_j
-
\lf(b-\fr1h\lf(b+\fr12gq\rg)\rg)
\fr{K}q
\fr1q{\cal H}^k_m
l^j \nabla_nb_j
$$

\ses

\ses

$$
-
K
\lf(\fr1{hq}(b+gq)-\fr bq\rg)
\fr 1{q}
{\cal H}^k_m
l^j
 \nabla_nb_j
+
\lf(\fr1{hq}-\fr1q +\fr{gb}{B}\rg)Km^k
 \nabla_nb_m
-
a^k{}_{nm}
 $$
\ses
which is reduced to read
\ses\\
 $$
 N^k_{nm}
=
-l^k
g\fr qK  \fr{K^2}B   \nabla_nb_m
+
\fr1K\lf[q-\fr1hq+\fr1{2h}g(b+gq)\rg]
{\cal H}^{kj}m_m\fr{K^2}B
 \nabla_nb_j
 $$

\ses

\ses

$$
+
\lf(b-\fr1h\lf(b+\fr12gq\rg)\rg)
\fr{K}B
\lf(
{\cal H}^{kj}l_m-l^k{\cal H}^j_m
+\fr bq m^k {\cal H}^j_m
\rg)
 \nabla_nb_j
$$

\ses

$$
+
K
\fr b{q^2}
\lf(\fr1{h}-1 \rg)
 m_m
 m^k
l^j
 \nabla_nb_j
-
\fr Kq\lf(\fr1{h}-1 \rg)
m_ml^k
l^j
 \nabla_nb_j
-
g\fr1{2qh}
{\cal H}^k_m
y^j \nabla_nb_j
$$

\ses

\ses

\be
+
\lf(\fr1{hq}-\fr1q +\fr{gb}{B}\rg)Km^k
 \nabla_nb_m
-
a^k{}_{nm}.
\ee

{%\pgbrk}

If we write this expression in the form
 $$
 N^k_{nm}
=
-l^k
g\fr qK  \fr{K^2}B   \nabla_nb_m
+
\fr1K\lf[q-\fr1hq+\fr1{2h}g(b+gq)\rg]
{\cal H}^{kj}m_m\fr{K^2}B
 \nabla_nb_j
 $$

\ses

\ses

$$
+
\lf(b-\fr1h\lf(b+\fr12gq\rg)\rg)
\fr{K}B
\lf(
{\cal H}^{kj}l_m-l^k{\cal H}^j_m
+\fr bq m^k {\cal H}^j_m
\rg)
 \nabla_nb_j
$$

\ses

$$
+
K\lf[
\fr1{hq^2}b -\fr1{q^2}b
 \rg]
 m_m
 m^k
l^j
 \nabla_nb_j
-
\fr Kq\lf(\fr1{h}-1 \rg)
m_ml^k
l^j
 \nabla_nb_j
-
g\fr1{2qh}
{\cal H}^k_m
y^j \nabla_nb_j
$$

\ses

\ses

$$
+
\lf(\fr1{hq}-\fr1q +\fr{gb}{B}\rg)\lf(Km^k-\fr qby^k\rg)
 \nabla_nb_m
+
\lf(\fr1{hq}-\fr1q +\fr{gb}{B}\rg)\fr qby^k
 \nabla_nb_m
 -
a^k{}_{nm}
$$
\ses
and use
$$
{\cal H}_{ij}=
\Bigl( r_{ij}-\fr1{q^2}v_iv_j\Bigr)\fr {K^2}{B}
\quad
\text{and} \quad
{\cal H}_i{}^j=
r_i{}^j-\fr1{q^2}v_iv^j,
$$
we obtain
\ses\\
 $$
 N^k_{nm}
=
\fr1K\lf[q-\fr1hq+\fr1{2h}g(b+gq)\rg]
\Bigl( a^{kj}-\fr1{q^2}v^ky^j\Bigr)
m_m
 \nabla_nb_j
 $$

\ses

\ses

$$
+
\lf(b-\fr1h\lf(b+\fr12gq\rg)\rg)
\fr{1}B
\lf(
\Bigl( a^{kj}-\fr1{q^2}v^ky^j\Bigr)v_m
+\Bigl(a^{kj}-\fr1{q^2}v^ky^j\Bigr)(b+gq)b_m
\rg)
 \nabla_nb_j
$$

\ses

\ses

$$
+
\lf(b-\fr1h\lf(b+\fr12gq\rg)\rg)
\fr{1}B
\lf(
-l^k+\fr bqm^k\rg)
 \nabla_nb_m
$$

\ses

$$
+
\fr1{q^2}\lf(b-\fr1h\lf(b+\fr12gq\rg)\rg)
\fr{1}B
\lf(
       y^k
-\fr bq Km^k
\rg)
v_my^j \nabla_nb_j
$$

\ses

$$
+
K\lf[
\fr1{hq^2}b -\fr1{q^2}b
 \rg]
 m_m
 m^k
l^j
 \nabla_nb_j
-
\fr Kq\lf(\fr1{h}-1 \rg)
m_ml^k
l^j
 \nabla_nb_j
-
g\fr1{2qh}
{\cal H}^k_m
y^j \nabla_nb_j
$$

\ses

\ses

$$
+
\lf(\fr1{hq}-\fr1q +\fr{gb}{B}\rg)\lf(Km^k-\fr qby^k\rg)
 \nabla_nb_m
+
\lf(\fr1{hq}-\fr1q \rg)\fr qby^k
 \nabla_nb_m
 -
a^k{}_{nm}.
$$

\ses

{%\pgbrk}

Taking into account the formula
$$
m_m=\fr K{qB}(q^2b_m-bv_m),
$$
we arrive at
\ses\\
 $$
 N^k_{nm}
=
\lf[q-\fr1hq+\fr1{2h}g(b+gq)\rg]
\fr 1{qB}q^2b_m
  a^{kj}\nabla_nb_j
-
\lf[q-\fr1hq+\fr1{2h}g(b+gq)\rg]
\fr 1{qB}bv_m
  a^{kj}\nabla_nb_j
   $$

\ses

$$
-
\lf[q-\fr1hq+\fr1{2h}g(b+gq)\rg]
\fr1{q^2}v^k
 \fr 1{qB}q^2b_m
y^j \nabla_nb_j
+
\lf[q-\fr1hq+\fr1{2h}g(b+gq)\rg]
\fr1{q^2}v^k
 \fr 1{qB}bv_m
y^j \nabla_nb_j
 $$

\ses

\ses

$$
+
\lf(b-\fr1h\lf(b+\fr12gq\rg)\rg)
\fr{1}B
\lf(
v_m
+ (b+gq)b_m
\rg)
 a^{kj}
 \nabla_nb_j
$$

\ses

$$
-
\lf(b-\fr1h\lf(b+\fr12gq\rg)\rg)
\fr{1}B
\fr1{q^2}v^kv_m
              y^j\nabla_nb_j
-
\lf(b-\fr1h\lf(b+\fr12gq\rg)\rg)
\fr{1}B
\fr1{q^2}v^k(b+gq)b_m
 y^j\nabla_nb_j
$$

\ses

\ses

$$
+\fr1h\fr12gq
\fr{K}B
\lf(
-l^k+\fr bqm^k\rg)
 \nabla_nb_m
$$

\ses

$$
+
\fr1{q^2}\lf(b-\fr1h\lf(b+\fr12gq\rg)\rg)
\fr{1}B
\lf(
       y^k
-\fr bq Km^k
\rg)
v_my^j \nabla_nb_j
$$

\ses

$$
-
\fr qb
\lf[
\fr1{hq^2}b -\fr1{q^2}b
 \rg]
 m_m
 \lf( y^k
-\fr bq Km^k\rg)
l^j
 \nabla_nb_j
$$

\ses

\ses

$$
+
\fr qb
\lf[
\fr1{hq^2}b -\fr1{q^2}b
 \rg]
 m_m
y^k
l^j
 \nabla_nb_j
-
\fr Kq\lf(\fr1{h}-1 \rg)
m_ml^k
l^j
 \nabla_nb_j
-
g\fr1{2qh}
{\cal H}^k_m
y^j \nabla_nb_j
$$

\ses

\ses

$$
+
\lf(\fr1{hq}-\fr b{hqB}(b+gq)-\fr qB\rg)\lf(Km^k-\fr qby^k\rg)
 \nabla_nb_m
+
\lf(\fr1{hq}-\fr1q \rg)\fr qby^k
 \nabla_nb_m
 -
a^k{}_{nm}.
$$
\ses\\
Noting also
$$
y^k-\fr bq Km^k  =        \fr B{q^2}v^k,
$$

{%\pgbrk}

\nin
we observe that
\ses\\
 $$
 N^k_{nm}
=
\lf(1-\fr1h\rg)
b_m  a^{kj}\nabla_nb_j
-\fr g{2hq}
v_m
  a^{kj}\nabla_nb_j
   $$

\ses

\ses

$$
-
\lf(q-\fr1hq\rg)
v^k
 \fr 1{qB}b_m
y^j \nabla_nb_j
+
\lf[q-\fr1hq+\fr1{2h}g(b+gq)\rg]
\fr1{q^2}v^k
 \fr 1{qB}bv_m
y^j \nabla_nb_j
 $$

\ses

\ses

$$
-
\lf(b-\fr1h\lf(b+\fr12gq\rg)\rg)
\fr{1}B
\fr1{q^2}v^kv_m
 y^j\nabla_nb_j
-
\lf(b-\fr1hb\rg)
\fr{K}B
\fr1{q^2}v^k(b+gq)b_m
 y^j\nabla_nb_j
$$

\ses

\ses

$$
-\fr1h\fr g{2q} v^k
 \nabla_nb_m
+
\fr1{q^2}\lf(b-\fr1h\lf(b+\fr12gq\rg)\rg)
\fr{1}{q^2}v^k
v_my^j \nabla_nb_j
$$

\ses

\ses

$$
-
\fr qb
\lf[
\fr1{hq^2}b -\fr1{q^2}b
 \rg]
 m_m
\fr{B}{q^2}v^k
l^j
 \nabla_nb_j
-
g\fr1{2qh}
{\cal H}^k_m
y^j \nabla_nb_j
+
\lf(\fr1h-1\rg)b^k
 \nabla_nb_m
 -
a^k{}_{nm}.
$$
\ses
These coefficients fulfill the equality
 $ N^k_{nm}=- D^k{}_{nm}$ with $D^k{}_{nm}$ given by (6.49).

{%\pgbrk}

                                  \ses

\ses

\ses

\setcounter{equation}{0}

\nin { \bf Appendix C: Evaluation of curvature tensor}

\ses

\ses

The substitution
$$
 a^{nt}= \fr {K^2}B{\cal H}^{nt}
+b^nb^t+\fr1{q^2}v^nv^t
$$
changes the representation  (6.75) to
the form
$$
E_k{}^n{}_{ij}=
-
\Bigl((1-h)b_k  +\fr12\fr gqv_k   \Bigr)
\fr {K^2}B{\cal H}^{nt}
\fr1h
b_ma_t{}^m{}_{ij}
$$

\ses

$$
-
\Biggl[
\Bigl((1-h)b_k  +\fr12\fr gqv_k   \Bigr)
\fr1{q^2}v^ny^t
+   \fr g{2q} \eta^n_k  y^t
+
\lf(
\fr g{2q}v^n
-(1-h)b^n
\rg)
 \de^t_k
\Biggr]
\fr1h
b_ma_t{}^m{}_{ij}
+  a_k{}^n{}_{ij}.
$$
\ses
Noting also
the expansion
$$
\de_k^t=
{\cal H}_k^t
+b_kb^t+\fr1{q^2}v_kv^t,
$$
we obtain
$$
E_k{}^n{}_{ij}=
-
\Bigl((1-h)b_k  +\fr12\fr gqv_k   \Bigr)
\fr {K^2}B{\cal H}^{nt}
\fr1h
b_ma_t{}^m{}_{ij}
-
\lf(
\fr g{2q}v^n
-(1-h)b^n
\rg)
{\cal H}_k^t
\fr1h
b_ma_t{}^m{}_{ij}
$$

\ses

$$
-
\Biggl[
\Bigl((1-h)b_k  +\fr12\fr gqv_k   \Bigr)
\fr1{q^2}v^ny^t
+   \fr g{2q} \eta^n_k  y^t
+
\lf(
\fr g{2q}v^n
-(1-h)b^n
\rg)
\fr1{q^2}v_ky^t
\Biggr]
\fr1h
b_ma_t{}^m{}_{ij}
$$

\ses

$$
+h^n_l  a_k{}^l{}_{ij}+  l^nl_la_k{}^l{}_{ij}
$$

\ses

\ses

\ses

$$
=
-
\Bigl((1-h)b_k  +\fr12\fr gqv_k   \Bigr)
\fr {K^2}B{\cal H}^{nt}
\fr1h
b_ma_t{}^m{}_{ij}
-
\lf(
\fr g{2q}v^n
-(1-h)b^n
\rg)
{\cal H}_k^t
\fr1h
b_ma_t{}^m{}_{ij}
$$

\ses

\ses

$$
-
 \fr g{2q}
{\cal H}^n_k  y^t
\fr1h
b_ma_t{}^m{}_{ij}
-
\Biggl[
(1-h)(b_k v^n-b^nv_k)
+
\fr g{q}v^n
v_k
\Biggr]
\fr1{q^2}
\fr1h
b_ma_t{}^m{}_{ij}y^t
$$

\ses

\ses

$$
+h^n_l  h_k^pa_p{}^l{}_{ij}
+h^n_l  l_kl^pa_p{}^l{}_{ij}
+  l^nl_lh_k^pa_p{}^l{}_{ij}
+  l^nl_ll_kl^pa_p{}^l{}_{ij}.
$$

{%\pgbrk}

Let us consider the term
$$
b_k v^n-b^nv_k
+
\fr g{q}v^n
v_k
=
y^n\lf(b_k+\fr g{q}v_k\rg)
-
\lf(
bb_k +v_k
+
\fr g{q}bv_k
\rg)
b^n
$$
\ses
and apply
(A.63)--(A.65).
We obtain
$$
b_k v^n-b^nv_k
+
\fr g{q}v^n
v_k
=
\fr1bl^n
\lf(bqm_k+ (B-q^2)l_k - g(B-q^2)m_k\rg)
+
B\fr1K
( gm_k-l_k)
b^n.
$$
\ses
Using here
$$
b^n=K\fr1B\lf[qm^n+(b+gq)l^n\rg]
$$
leads to
$$
b_k v^n-b^nv_k
+
\fr g{q}v^n
v_k
=
\fr1bl^n
\lf(bqm_k+ (B-q^2)l_k - g(B-q^2)m_k\rg)
$$

\ses

$$
+
 gm_k
\lf[qm^n+(b+gq)l^n\rg]
-
l_k
\lf[qm^n+(b+gq)l^n\rg]
=
 q[  l^n  m_k    +    gm_k   m^n-  l_k   m^n   ].
$$
\ses
Therefore, the term
$$
b_k v^n-b^nv_k
=
l^n(qm_k+ bl_k)
-
\fr1K\lf(b(qm_k+ bl_k)
+     q               \lf(q l_k -(b+gq)m_k\rg)\rg)
b^n
$$
\ses
can be traced to be
$$
b_k v^n-b^nv_k
=
l^n(qm_k+ bl_k)
-
\fr1K\lf((B-gbq)l_k-gq^2m_k
\rg)
b^n
$$

\ses

$$
=
l^n(qm_k+ bl_k)
-
\lf((B-gbq)l_k-gq^2m_k
\rg)
\fr1B
qm^n
-
\lf((B-gbq)l_k-gq^2m_k
\rg)
\fr1B
(b+gq)l^n.
$$
In this way we come to
$$
b_k v^n-b^nv_k
=
l^n(qm_k+ bl_k)
-
\lf((B-gbq)l_k-gq^2m_k
\rg)
\fr1B
qm^n
$$

\ses

$$
-
l_k
(b+gq)l^n
+
gbql_k
\fr1B
(b+gq)l^n
+
gq^2m_k
\fr1B
(b+gq)l^n,
$$
\ses
or
\be
b_k v^n-b^nv_k
=
ql^nm_k
-
\lf((B-gbq)l_k-gq^2m_k
\rg)
\fr1B
qm^n
-\fr{gq^3}Bl_kl^n
+
gq^2m_k
\fr1B
(b+gq)l^n.
\ee

It will be noted also that
\be
l_la_t{}^l{}_{ij}y^t=gq\fr{K}B
b_la_t{}^l{}_{ij}y^t.
\ee

{%\pgbrk}

Thus we can write
\ses\\
$$
E_k{}^n{}_{ij}=
-
\Bigl((1-h)b_k  +\fr12\fr gqv_k   \Bigr)
\fr {K^2}B{\cal H}^{nt}
\fr1h
b_ma_t{}^m{}_{ij}
-
\lf(
\fr g{2q}v^n
-(1-h)b^n
\rg)
{\cal H}_k^t
\fr1h
b_ma_t{}^m{}_{ij}
$$

\ses

\ses

$$
-
 \fr g{2q}
{\cal H}^n_k  y^t
\fr1h
b_ma_t{}^m{}_{ij}
-
[  l^n  m_k    +    gm_k   m^n-  l_k   m^n   ]\fr1{q}
\fr1h
b_ma_t{}^m{}_{ij}y^t
$$

\ses

\ses

$$
+
\Biggl[
ql^nm_k
-
\lf((B-gbq)l_k-gq^2m_k
\rg)
\fr1B
qm^n
-\fr{gq^3}Bl_kl^n
+
gq^2m_k
\fr1B
(b+gq)l^n
\Biggr]
\fr1{q^2}
b_ma_t{}^m{}_{ij}y^t
$$

\ses

\ses

$$
+h^n_l  h_k^pa_p{}^l{}_{ij}
+h^n_l  l_kl^pa_p{}^l{}_{ij}
+  l^nl_lh_k^pa_p{}^l{}_{ij}
+  l^nl_k
gq\fr{1}B
b_la_t{}^l{}_{ij}y^t.
$$

Applying  the equalities
$$
y_i=(u_i+gqb_i)\fr{K^2}B,
\qquad
b^n=\fr1B\lf[Kqm^n+(b+gq)y^n\rg]
$$
to
$$
l_lKqm^ta_t{}^l{}_{ij}=
l_l\lf[Bb^t-(b+gq)y^t \rg]a_t{}^l{}_{ij}
$$
yields
$$
l_lBqm^ta_t{}^l{}_{ij}=
(u_l+gqb_l)
\lf[Bb^t-(b+gq)y^t \rg]a_t{}^l{}_{ij}
=
\lf[Bb^tu_l-gqb_l(b+gq)y^t \rg]a_t{}^l{}_{ij},
$$
so that
\be
l_lm^ta_t{}^l{}_{ij}=
-\fr1{Bq}\lf[B+gq(b+gq) \rg]b_la_t{}^l{}_{ij}y^t.
\ee

{%\pgbrk}

So we have
$$
E_k{}^n{}_{ij}=
-
\Bigl((1-h)b_k  +\fr12\fr gqv_k   \Bigr)
\fr {K^2}B{\cal H}^{nt}
\fr1h
b_ma_t{}^m{}_{ij}
-
\lf(
\fr g{2q}v^n
-(1-h)b^n
\rg)
{\cal H}_k^t
\fr1h
b_ma_t{}^m{}_{ij}
$$

\ses

$$
-
 \fr g{2q}
{\cal H}^n_k  y^t
\fr1h
b_ma_t{}^m{}_{ij}
-
[  l^n  m_k    +    gm_k   m^n-  l_k   m^n   ]\fr1{q}
\fr1h
b_ma_t{}^m{}_{ij}y^t
$$

\ses

\ses

$$
+
\Biggl[
ql^nm_k
-
\lf((B-gbq)l_k-gq^2m_k
\rg)
\fr1B
qm^n
+
gq^2m_k
\fr1B
(b+gq)l^n
\Biggr]
\fr1{q^2}
b_ma_t{}^m{}_{ij}y^t
$$

\ses

\ses

$$
+h^n_l  h_k^ta_t{}^l{}_{ij}
+h^n_l  l_kl^ta_t{}^l{}_{ij}
+  l^nl_l {\cal H}_k^ta_t{}^l{}_{ij}
-\fr1{Bq}\lf[B+gq(b+gq) \rg]l^nm_k
b_la_t{}^l{}_{ij}y^t
$$

\ses

\ses

\ses

$$
=
-
\Bigl((1-h)b_k  +\fr12\fr gqv_k   \Bigr)
\fr {K^2}B{\cal H}^{nt}
\fr1h
b_ma_t{}^m{}_{ij}
-
\lf(
\fr g{2q}v^n
-(1-h)b^n
\rg)
{\cal H}_k^t
\fr1h
b_ma_t{}^m{}_{ij}
$$

\ses

$$
-
 \fr g{2q}
{\cal H}^n_k  y^t
\fr1h
b_ma_t{}^m{}_{ij}
-
[  l^n  m_k    +    gm_k   m^n-  l_k   m^n   ]\fr1{q}
\fr1h
b_ma_t{}^m{}_{ij}y^t
+gq^2m_k
\fr1B
m^n
\fr1{q}
b_ma_t{}^m{}_{ij}y^t
$$

\ses

\ses

$$
+{\cal H}^n_l  l_kl^ta_t{}^l{}_{ij}
+  l^nl_l {\cal H}_k^ta_t{}^l{}_{ij}
b_la_t{}^l{}_{ij}y^t
$$

\ses

$$
+{\cal H}^n_l  {\cal H}_k^ta_t{}^l{}_{ij}
+{\cal H}^n_l  m_km^ta_t{}^l{}_{ij}
+m^nm_l  {\cal H}_k^ta_t{}^l{}_{ij}
+m^n m_k
m_l m^ta_t{}^l{}_{ij}.
$$

{%\pgbrk}

The next step is to consider the term
$$
m_lKqm^ta_t{}^l{}_{ij}=
m_l\lf[Bb^t-(b+gq)y^t \rg]a_t{}^l{}_{ij}
$$
and use
$
(1/K)
qBm_l=q^2b_l-bv_l=S^2b_l-bu_l,
$
obtaining
$$
qBm_l
qm^ta_t{}^l{}_{ij}=
(S^2b_l-bu_l)
\lf[Bb^t-(b+gq)y^t \rg]a_t{}^l{}_{ij}
=
-
\lf[bu_lBb^t+S^2b_l(b+gq)y^t \rg]a_t{}^l{}_{ij},
$$
so that
\be
m_lm^ta_t{}^l{}_{ij}
=
-\fr{gq}B
b_la_t{}^l{}_{ij}y^t.
\ee

The studied tensor takes now the form
$$
E_k{}^n{}_{ij}=
-
\Bigl((1-h)b_k  +\fr12\fr gqv_k   \Bigr)
\fr {K^2}B{\cal H}^{nt}
\fr1h
b_ma_t{}^m{}_{ij}
-
\lf(
\fr g{2q}v^n
-(1-h)b^n
\rg)
{\cal H}_k^t
\fr1h
b_ma_t{}^m{}_{ij}
$$

\ses

\ses

$$
-
 \fr g{2q}
{\cal H}^n_k
\fr1h
b_la_t{}^l{}_{ij}y^t
+{\cal H}^n_l  l_kl^ta_t{}^l{}_{ij}
+  l^nl_l {\cal H}_k^ta_t{}^l{}_{ij}
$$

\ses

\ses

$$
+{\cal H}^n_l  {\cal H}_k^ta_t{}^l{}_{ij}
+{\cal H}^n_l  m_km^ta_t{}^l{}_{ij}
+m^nm_l  {\cal H}_k^ta_t{}^l{}_{ij}
-
[  l^n  m_k    +    gm_k   m^n-  l_k   m^n   ]\fr1{q}
\fr1h
b_la_t{}^l{}_{ij}y^t
$$

\ses

\ses

 \ses

$$
=
-
\Bigl((1-h)b_k  +\fr12\fr gqv_k   \Bigr)
\fr {K^2}B{\cal H}^{nt}
\fr1h
b_la_t{}^l{}_{ij}
-
[  l^n  m_k    +    gm_k   m^n-  l_k   m^n   ]\fr1{q}
\fr1h
b_la_t{}^l{}_{ij}y^t
$$

\ses

$$
-
\lf[\fr g{2q}y^n
-\lf(\fr g{2q}b+1\rg)
K\fr1B\lf[qm^n+(b+gq)l^n\rg]
+
h
K\fr1B\lf[qm^n+(b+gq)l^n\rg]
\rg]
{\cal H}_k^t
\fr1h
b_ma_t{}^m{}_{ij}
$$

\ses

\ses

$$
-
 \fr g{2q}
{\cal H}^n_k
\fr1h
b_la_t{}^l{}_{ij}y^t
+{\cal H}^n_l  l_kl^ta_t{}^l{}_{ij}
+  l^nl_l {\cal H}_k^ta_t{}^l{}_{ij}
$$

\ses

\ses

$$
+{\cal H}^n_l  {\cal H}_k^ta_t{}^l{}_{ij}
+{\cal H}^n_l  m_km^ta_t{}^l{}_{ij}
+m^nm_l  {\cal H}_k^ta_t{}^l{}_{ij}.
$$

Recollecting
$$
b_k=\fr1K(qm_k+ bl_k), \qquad
bv_k=
\fr1Kq               \lf(q bl_k -(B-q^2)m_k\rg),
$$

 {%\pgbrk}

\ses

\nin
we can eventually write
$$
E_k{}^n{}_{ij}=
-
\Bigl((1-h)(qm_k+ bl_k)  +\fr g2
             \lf(q l_k -(b+gq)m_k\rg)   \Bigr)
\fr {K}B{\cal H}^{nt}
\fr1h
b_la_t{}^l{}_{ij}
$$

\ses

$$
-
\lf[\fr g{2q}y^n
-\fr g{2q}b
K\fr1B\lf[qm^n+(b+gq)l^n\rg]
-
K\fr1B\lf[qm^n+(b+gq)l^n\rg]
\rg]
{\cal H}_k^t
\fr1h
b_ma_t{}^m{}_{ij}
$$

\ses

\ses

$$
-
K\fr1B\lf[qm^n+(b+gq)l^n\rg]
{\cal H}_k^t
b_ma_t{}^m{}_{ij}
-
 \fr g{2q}
{\cal H}^n_k
\fr1h
b_la_t{}^l{}_{ij}y^t
$$

\ses

\ses

\ses

$$
+{\cal H}^n_l  l_kl^ta_t{}^l{}_{ij}
+  l^nl_l {\cal H}_k^ta_t{}^l{}_{ij}
+{\cal H}^n_l  {\cal H}_k^ta_t{}^l{}_{ij}
+{\cal H}^n_l  m_km^ta_t{}^l{}_{ij}
+m^nm_l  {\cal H}_k^ta_t{}^l{}_{ij}
$$

\ses

\ses

\be
-
[  l^n  m_k    +    gm_k   m^n-  l_k   m^n   ]\fr1{q}
\fr1h
b_la_t{}^l{}_{ij}y^t.
\ee
Here,
$$
{\cal H}^n_l a_t{}^l{}_{ij}=-\fr{K^2}B{\cal H}^{nl}a_{tlij}.
$$

\ses

{%\pgbrk}

 The tensor
$$
E_k{}^n{}_{ij}+\fr1Kl_kM^n{}_{ij}
-\fr1K
l^nM_{kij}
=
-
\Bigl((1-h)q
-\fr g2
 (b+gq)  \Bigr)
m_k
\fr {K}B{\cal H}^{nt}
\fr1h
b_la_t{}^l{}_{ij}
$$

\ses

\ses

$$
-
\lf[\fr {gq}{2}l^n
-\fr g{2}b
m^n
-
qm^n-(b+gq)l^n
\rg]
K\fr1B
{\cal H}_k^t
\fr1h
b_la_t{}^l{}_{ij}
$$

\ses

\ses

$$
-
K\fr1B\lf[qm^n+(b+gq)l^n\rg]
{\cal H}_k^t
b_la_t{}^l{}_{ij}
-
 \fr g{2q}
{\cal H}^n_k
\fr1h
b_la_t{}^l{}_{ij}y^t
+  l^n
(u_l+gqb_l)\fr{K}B
 {\cal H}_k^ta_t{}^l{}_{ij}
$$

\ses

$$
+{\cal H}^n_l  {\cal H}_k^ta_t{}^l{}_{ij}
                                          +{\cal H}^n_l  m_km^ta_t{}^l{}_{ij}
+m^nm_l  {\cal H}_k^ta_t{}^l{}_{ij}
-
(  l^n   +    g  m^n )
m_k \fr1{q}
\fr1h
b_la_t{}^l{}_{ij}y^t
$$

\ses

\ses

\ses

$$
=
\Biggl[
-
\Bigl((1-h)q
-\fr g2
 (b+gq)
  \Bigr)
   b^l
-
h
K
m^l
\Biggr]
m_k
\fr {K}B{\cal H}^{nt}
\fr1h
a_{tlij}
$$

\ses
\ses

$$
+
\lf(\fr g{2}b+q\rg)m^n
K\fr1B
{\cal H}_k^t
\fr1h
b_la_t{}^l{}_{ij}
-
K\fr1Bqm^n
{\cal H}_k^t
b_la_t{}^l{}_{ij}
$$
\ses

$$
-
 \fr g{2q}
{\cal H}^n_k
\fr1h
b_la_t{}^l{}_{ij}y^t
+{\cal H}^n_l  {\cal H}_k^ta_t{}^l{}_{ij}
+m^nm_l  {\cal H}_k^ta_t{}^l{}_{ij}
-
    g  m^n
m_k \fr1{q}
\fr1h
b_la_t{}^l{}_{ij}y^t
$$

\ses

\ses
        \ses

$$
=
\Biggl[
-
\Bigl(q
-\fr g2
 (b+gq)
  \Bigr)
   b^l
-
h
  \fr1q[(B-q^2)b^l-(b+gq)y^l]
\Biggr]
m_k
\fr {K}B{\cal H}^{nt}
\fr1h
a_{tlij}
$$

\ses
\ses

$$
+
\lf(\fr g{2}b+q\rg)m^n
K\fr1B
{\cal H}_k^t
\fr1h
b_la_t{}^l{}_{ij}
-
K\fr1Bqm^n
{\cal H}_k^t
b_la_t{}^l{}_{ij}
$$

\ses

\ses

$$
-
 \fr g{2q}
{\cal H}^n_k
\fr1h
b_la_t{}^l{}_{ij}y^t
+{\cal H}^n_l  {\cal H}_k^ta_t{}^l{}_{ij}
+m^n
     \fr K{qB}    [(b^2+q^2)b_l-bu_l]
  {\cal H}_k^ta_t{}^l{}_{ij}
-
    g  m^n
m_k \fr1{q}
\fr1h
b_la_t{}^l{}_{ij}y^t
$$
\ses
is found in the simple expansion form
$$
E_k{}^n{}_{ij}+\fr1Kl_kM^n{}_{ij}
-\fr1K
l^nM_{kij}
=
\!-\!
\Biggl[
\fr1h
\Bigl(q
-\fr g2
 (b+gq)
  \Bigr)
   b^l
+
  \fr1q[(B-q^2)b^l-(b+gq)y^l]
\Biggr]
m_k
\fr {K}B{\cal H}^{nt}
a_{tlij}
$$

\ses

\ses

$$
+
\Biggl[
\lf(\fr g{2}bq+q^2\rg)
b_l
  + h   (b^2b_l-bu_l)
\Biggr]
m^n
   K\fr1{Bq}
{\cal H}_k^t
\fr1h
a_t{}^l{}_{ij}
$$

\ses

\ses

\be
-
 \fr g{2q}
{\cal H}^n_k
\fr1h
b_la_t{}^l{}_{ij}y^t
+{\cal H}^n_l  {\cal H}_k^ta_t{}^l{}_{ij}
-
    g  m^n
m_k \fr1{q}
\fr1h
b_la_t{}^l{}_{ij}y^t.
\ee

\ses

{%\pgbrk}

The contraction
$$
\fr B{K^2}
h^2
(KE_{knij}+l_kM_{nij}-l_nM_{kij})
(KE^{knij}+l^kM^{nij}-l^nM^{kij})
$$

\ses

     \ses

$$
=
\Biggl[
\Bigl(q
-\fr g2
 (b+gq)
+
h
  \fr1q(B-q^2)
    \Bigr)
   b^l
-
\fr1q
h
(b+gq)y^l
\Biggr]
a_{tl}{}^{ij}
\times
$$

\ses

\ses

$$
\Biggl[
\Bigl(q
-\fr g2
 (b+gq)
+
h
  \fr1q(B-q^2)
    \Bigr)
   b^h
-
\fr1qh
(b+gq)y^h
\Biggr]a^t{}_{h}{}_{ij}
$$

\ses

\ses

\ses

\ses

$$
-
\fr{S^2}{q^2}
\fr1{q^2}h^2(b+gq)^2
b^t   y^la_{tl}{}^{ij}    y^h    b^s    a_{shij}
$$

\ses

\ses

$$
+
2\fr{b}{q^2}
\Biggl[
\Bigl(q
-\fr g2
 (b+gq)
  \Bigr)
+
h
  \fr1q(B-q^2)
\Biggr]
\fr1q
h
(b+gq)
b^t   y^la_{tl}{}^{ij}    y^h    b^s    a_{shij}
$$

\ses

\ses

\ses

$$
-
\fr{1}{q^2}
y^ty^s
\Biggl[
\Bigl(q
-\fr g2
 (b+gq)
  \Bigr)
   b^l
+
h
  \fr1q(B-q^2)b^l
\Biggr]
a_{tl}{}^{ij}
\Biggl[
\Bigl(q
-\fr g2
 (b+gq)
  \Bigr)
   b^h
+
h
  \fr1q(B-q^2)b^h
\Biggr]
a_{shij}
 $$

\ses

\ses

     \ses

$$
+
\fr1{q^2}
\Biggl[
\lf(\fr g{2}bq+q^2  +  h   b^2\rg)
b^l
- h   by^l
\Biggr]
a_{tl}{}^{ij}
\Biggl[
\lf(\fr g{2}bq+q^2  +  h   b^2\rg)
b^h
- h   by^h
\Biggr]
a^t{}_{h}{}_{ij}
$$

\ses

\ses

\ses

\ses

$$
-
\fr1{q^2}
\fr{S^2}{q^2}
h^2b^2
b^t
y^la_{tl}{}^{ij}
y^h
b^s
a_{shij}
+2h
\fr1{q^2}
\fr{b^2}{q^2}
\lf(\fr g{2}bq+q^2  +  h   b^2\rg)
b^t   y^la_{tl}{}^{ij}    y^h    b^s    a_{shij}
$$

\ses

\ses

\ses

$$
-\fr1{q^2}
\fr{1}{q^2}
y^ty^s
\lf(\fr g{2}bq+q^2  +  h   b^2\rg)
b^l
a_{tl}{}^{ij}
\lf(\fr g{2}bq+q^2  +  h   b^2\rg)
b^h
a_{shij}
 $$

\ses

\ses

$$
+
(N-2)\fr{g^2}{4q^2}
b_la_t{}^{lij}y^t
b_ha_s{}^h{}_{ij}y^s
B
+
g^2
 \fr1{q^2}
b_la_t{}^{lij}y^t
b_ha_s{}^h{}_{ij}y^s
B
+Bh^2{\cal H}_{lh}  {\cal H}^{ts}
a_t{}^{lij}
a_s{}^h{}_{ij}
$$

{%\pgbrk}

\ses

\nin
can be written in the simple form
$$
(KE_{knij}+l_kM_{nij}-l_nM_{kij})
(KE^{knij}+l^kM^{nij}-l^nM^{kij})    =   KE_{knij}
(KE^{knij}+l^kM^{nij}-l^nM^{kij})
$$

\ses

\be
=
K^2E_{knij}E^{knij}
-2
M_{kij}M^{kij}.
\ee

\ses

\ses

Therefore,
\be
BE_{knij}E^{knij}=
\fr B{K^2} (KA_{knij}+l_kM_{nij}-l_nM_{kij})
(KA^{knij}+l^kM^{nij}-l^nM^{kij})
+
2
\fr B{K^2}
M_{kij}M^{kij}.
\ee
\ses
The term
$M_{kij}M^{kij}$
can be taken from
(6.82).

{%\pgbrk}

The contraction (C.8) can be written in the alternative form
$$
Bh^2E_{knij}E^{knij}=
Bh^2
a_{lh}
a^{ts}
a_t{}^{lij}
a_s{}^h{}_{ij}
 $$

\ses

\ses

$$
+
\fr1{q^2}
T
b^l
a_{tl}{}^{ij}
b^h
a^t{}_{h}{}_{ij}
-
2h\fr1{q^2}
(Z-2hbB)
b^l
a_{tl}{}^{ij}
   y^h
a^t{}_{h}{}_{ij}
+
\fr{1}{q^2}
 h^2
 [g^2q^2-2q^2]
    y^l
a_{tl}{}^{ij}
y^h
a^t{}_{h}{}_{ij}
$$

\ses

\ses

\ses

         \ses

$$
+
g^2\fr{B}{q^2}
\Bigl(
(N-2)\fr{1}{4}
+
1
\Bigr)
b^t   y^la_{tl}{}^{ij}    y^h    b^s    a_{shij}
$$

\ses

\be
+2
\Biggl(
\Bigl((1-h)b  +\fr g2q    \Bigr)
b_ha^{nhij}
-
h
a_h{}^{nij}
 y^h
\Biggr)
\Biggl(
\Bigl((1-h)b  +\fr g2q    \Bigr)
b_la_n{}^l{}_{ij}
-
h
a_{tnij}
 y^t
\Biggr).
\ee

{%\pgbrk}

\nin
Here,
$$
T=
\lf(\fr g{2}bq+q^2  +  h   b^2\rg)^2
+
 \Bigl(q^2
-\fr g2q
 (b+gq)
+
h
(B-q^2)
    \Bigr)^2
-2h^2b^2B    -2h^2q^2B
$$

\ses

$$
=
\lf(\fr g{2}bq+q^2 \rg)^2
+
2\lf(\fr g{2}bq+q^2\rg)h   b^2+h^2b^4
$$

\ses

$$
+
 \Bigl(q^2
-\fr g2q
 (b+gq)
    \Bigr)^2
+
2
 \Bigl(q^2
-\fr g2q
 (b+gq)
    \Bigr)
hb(b+gq)
+
h^2b^2(b^2+2gbq+g^2q^2)
$$

\ses

$$
-2h^2b^2(b^2+gbq+q^2)    -2h^2q^2(b^2+gbq+q^2),
$$
\ses
so that
$$
T+2\Bigl((1-h)b  +\fr g2q    \Bigr)^2q^2
=
2(1-2h+h^2)b^2q^2  + 2gqb(1-h)q^2
 +2(1-h^2)q^4
 $$

\ses

$$
+
\lf(\fr g{2}bq+q^2 \rg)^2
+
2q^2h   b^2
$$

\ses

$$
+
 \Bigl(q^2
-\fr g2q
 (b+gq)
    \Bigr)^2
+  2   \Bigl(q^2-\fr g2q gq    \Bigr)hb^2
+  2   \Bigl(q^2-\fr g2q (b+gq)    \Bigr)hbgq
+h^2g^2b^2q^2
-
2
h^2b^2q^2
$$

\ses

$$
 -2h^2q^2(b^2+gbq+q^2)
 $$

\ses

\ses

\ses

$$
=
4b^2q^2  + 2gqbq^2
 +4(1-h^2)q^4
-g^2q^4+\fr{g^2}4q^2(2gbq+g^2q^2)
$$

\ses

$$
- g^2q^2hb(2b+gq)
+h^2g^2b^2q^2
-
2
h^2b^2q^2
 -2h^2q^2(b^2+gbq).
 $$
\ses
Thus we have
simply
\be
\fr1{q^2}
\Biggl[T+2\Bigl((1-h)b  +\fr g2q    \Bigr)^2q^2  \Biggr]
=
g^2
\lf(
-\fr {1}2 (2b+gq)+hb
\rg)^2.
\ee

{%\pgbrk}

Another coefficient is
$$
Z=
\Bigl(q^2
-\fr g2q
 (b+gq)
+
h
(B-q^2)
    \Bigr)
b
+\Bigl(q^2
-\fr g2q
 (b+gq)
+
h
(B-q^2)
    \Bigr)
gq
+b
\lf(\fr g{2}bq+q^2  +  h   b^2\rg).
$$
\ses
We can follow the steps
$$
Z-2q^2\lf((1-h)b+\fr g2q\rg)
=
-2q^2(1-h)b
+
\Bigl(q^2
-\fr g2q
 (b+gq)
+
h
(B-q^2)
    \Bigr)
b
$$

\ses

$$
+\Bigl(
-\fr g2q
 (b+gq)
+
h
(B-q^2)
    \Bigr)
gq
+b
\lf(\fr g{2}bq+q^2  +  h   b^2\rg)
$$

\ses

\ses

        \ses

$$
=
2q^2hb
-\fr g2qb(b+gq)
+hb^2
(b+gq)
+\Bigl(
-\fr g2q
+
h
b
    \Bigr)
gq (b+gq)
+b
\lf(\fr g{2}bq +  h   b^2\rg)
$$

\ses

\ses

        \ses

$$
=
-\fr{ g^2}2bq^2
-\fr {g^2}2 (b+gq)q^2+g^2hbq^2
+2hbB,
$$
\ses
obtaining
\be
\fr1{q^2}\Biggl[Z-2q^2\lf((1-h)b+\fr g2q\rg)
-2hbB\Biggr]
=
-\fr {g^2}2 (2b+gq)+g^2hb.
\ee

{%\pgbrk}

The contraction becomes eventually
\ses\\
$$
E_{knij}E^{knij}=  a_{knij} a^{knij}
+
g^2
\fr{1}{h^2q^2}
\Bigl(
(N-2)\fr{1}{4}
+
1
\Bigr)
b^t   y^la_{tl}{}^{ij}    y^h    b^s    a_{shij}
$$

\ses

\ses

\be
+
\fr{g^2}{B}
\Biggl[
\Biggl( b - \fr1h\lf(b+\fr12gq\rg)  \Biggr)
b_ha^{nhij}
-
a^n{}_h{}^{ij}
 y^h
\Biggr]
\Biggl[
\Biggl( b - \fr1h\lf(b+\fr12gq\rg)  \Biggr)
b_la_n{}^l_{ij}
-
a_{ntij}
 y^t
\Biggr].
\ee

{%\pgbrk}

Now we are to verify  the representation (6.77)--(6.78).

\ses

Using  (A.60) together with (6.73) yields
\ses\\
$$
\Rho_k{}^n{}_{ij}+\fr1Kl_kM^n{}_{ij} -\fr1K l^nM_{kij}
 = \lf[ -
\Bigl((1-h)q -\fr g2
 (b+gq)
  \Bigr)
   b^l
- h K m^l
\rg] m_k \fr {K}B{\cal H}^{nt} \fr1h a_{tlij}
$$

\ses \ses

$$
+
\lf(\fr g{2}b+q\rg)m^n
K\fr1B
{\cal H}_k^t
\fr1h
b_la_t{}^l{}_{ij}
-
K\fr1Bqm^n
{\cal H}_k^t
b_la_t{}^l{}_{ij}
$$
\ses

$$ -
 \fr g{2q}
{\cal H}^n_k \fr1h b_la_t{}^l{}_{ij}y^t +{\cal H}^n_l  {\cal
H}_k^ta_t{}^l{}_{ij} +m^nm_l  {\cal H}_k^ta_t{}^l{}_{ij} -
    g  m^n
m_k \fr1{q} \fr1h b_la_t{}^l{}_{ij}y^t
$$

\ses

$$
+
 {\cal H}^n_k
\fr g{2qh}
b_m a_t{}^m{}_{ij}y^t
-\fr1K
\fr1N (A_kM^n{}_{ij}+A^nM_{kij})
$$

\ses

\ses

\ses

$$
 = \Biggl[ -
\Bigl((1-h)q -\fr g2
 (b+gq)
  \Bigr)
   b^l
- h K m^l \Biggr] m_k \fr {K}B{\cal H}^{nt} \fr1h a_{tlij}
+
\lf(\fr g{2}b+q\rg)m^n
K\fr1B
{\cal H}_k^t
\fr1h
b_la_t{}^l{}_{ij}
$$

\ses

\ses

$$
-
K\fr1Bqm^n
{\cal H}_k^t
b_la_t{}^l{}_{ij}
 +
 {\cal H}^n_l  {\cal H}_k^ta_t{}^l{}_{ij} +m^nm_l  {\cal H}_k^ta_t{}^l{}_{ij} -
    g  m^n
m_k \fr1{q} \fr1h b_la_t{}^l{}_{ij}y^t
$$

\ses

\ses

$$
-\fr1K
\fr g2 m_k
\Biggl[
 \Bigl((1-h)b  +\fr12gq    \Bigr)\fr {K^2}B{\cal H}^{nt}
\fr1h
b_la_t{}^l{}_{ij}
-
{\cal H}^n_la_t{}^l{}_{ij}y^t
- \fr {1}{q}
Km^n
 \fr1h
y^tb_la_t{}^l{}_{ij}
\Biggr]
$$

\ses

\ses

$$
-\fr1K
\fr g2 m^n
\Biggl[
  \Bigl((1-h)b  +\fr12gq    \Bigr)\fr {K^2}B{\cal H}_k^t
\fr1h
b_la_t{}^l{}_{ij}
-
{\cal H}_{kl}a_t{}^l{}_{ij}y^t
- \fr {1}{q}
Km_k
 \fr1h
y^tb_la_t{}^l{}_{ij}
\Biggr],
$$
where (6.72) has been applied.

{%\pgbrk}

\nin
Reducing similar terms leaves us with
\ses\\
$$
\Rho_k{}^n{}_{ij}+\fr1Kl_kM^n{}_{ij}
-\fr1K
l^nM_{kij}
 = \Bigl[
\fr g2
 gq
   b^l
- h K m^l \Bigr] m_k \fr {K}B{\cal H}^{nt} \fr1h a_{tlij}
$$

\ses \ses

$$
+
\fr1h(1-h)q
K\fr1B
({\cal H}_k^tm^n-{\cal H}^{nt}m_k)
b_la_t{}^l{}_{ij}
 +
 {\cal H}^n_l  {\cal H}_k^ta_t{}^l{}_{ij} +m^nm_l  {\cal H}_k^ta_t{}^l{}_{ij}
$$

\ses

$$
-\fr1K
\fr g2 m_k
\Biggl[
 \Bigl(-hb  +\fr12gq    \Bigr)\fr {K^2}B{\cal H}^{nt}
\fr1h
b_la_t{}^l{}_{ij}
-
{\cal H}^n_la_t{}^l{}_{ij}y^t
\Biggr]
$$

\ses

\ses

$$
-\fr1K
\fr g2 m^n
\Biggl[
  \Bigl(-hb  +\fr12gq    \Bigr)\fr {K^2}B{\cal H}_k^t
\fr1h
b_la_t{}^l{}_{ij}
-
{\cal H}_{kl}a_t{}^l{}_{ij}y^t
\Biggr].
$$
\ses
Here it is convenient to apply
the relation
\be
m_ia^{il}
=
\fr {K}{B}\Bigl[
\fr{q^2+b^2}BKm^l
+\fr{gq^2}By^l
\Bigl],
\ee
which comes from the chain
$$
m_ia^{il}=\fr {K}{qB}(q^2b^l-bv^l)
=
\fr {K}{qB}\Bigl[
q^2
\fr1B\lf[Kqm^l+(b+gq)y^l\rg]
-b
\fr qB\lf[-Kbm^l+qy^l\rg]
\Bigl].
$$

By lowering the index, we obtain
$$
\Rho_{knij}+\fr1Kl_kM_{nij} -\fr1K l_nM_{kij}
 =
\lf[\fr g2 gq   b^l- h K m^l \rg]
 m_k \fr {K}B{\cal H}_n^{t} \fr1h a_{tlij}
$$

\ses

$$
+
\fr1h(1-h)q
K\fr1B
({\cal H}_k^tm_n-{\cal H}_n^tm_k)
b_la_t{}^l{}_{ij}
$$

\ses

$$
 +
 {\cal H}_n^l  {\cal H}_k^ta_{tlij}
 \fr {K^2}B
  +m_n\fr {K}{B}\Bigl[
\fr{q^2+b^2}BKm^l
+\fr{gq^2}By^l
\Bigl]  {\cal H}_k^ta_{tlij}
$$

\ses

$$
-\fr1K
\fr g2 m_k
\Biggl[
 \Bigl(-hb  +\fr12gq    \Bigr){\cal H}_n^t
\fr1h
b_la_t{}^l{}_{ij}
-
{\cal H}_n^la_{tlij}y^t
\Biggr]
\fr {K^2}B
$$

\ses

$$
-\fr1K
\fr g2 m_n
\Biggl[
  \Bigl(-hb  +\fr12gq    \Bigr){\cal H}_k^t
\fr1h
b_la_t{}^l{}_{ij}
-
{\cal H}_k^la_{tlij}y^t
\Biggr]
\fr {K^2}B.
$$

{%\pgbrk}

\nin
Using  the equality
$
Bb^l=\lf[Kqm^l+(b+gq)y^l\rg]
$
leads to the following result after a short  simplification:
\ses\\
\be
\Rho_{knij}=-\fr1K(l_kM_{nij} - l_nM_{kij})
+(m_k{\cal H}_n^t-m_n{\cal H}_k^t)
P_{tij}
 +
{\cal H}_k^t {\cal H}_n^l  a_{tlij}
 \fr {K^2}B
 \ee
with
\ses\\
$$
P_{tij}=
-\fr1h(1-h)q
\fr KB
b^la_{tlij}
-\fr {K^2}{B^2}
                    \lf(B-\fr12gbq\rg)m^l
a_{tlij}
$$

\ses

\ses

\be
-
\fr g2
q^2y^la_{tlij}
\fr {K}{B^2}
+
\fr1h
\fr {g^2q}4
(b+gq)y^l
a_{tlij}
\fr {K}{B^2}
+
\fr1h
\fr {g^2q^2}4
m^l
a_{tlij}
\fr {K^2}{B^2}.
\ee
\ses
Inserting the vector
$$
m^l=\fr1{Kq}[Bb^l-(b+gq)y^l]
$$
leads to
$$
P_{tij}=
-\fr1h(1-h)q
\fr KB
b^la_{tlij}
-\fr {K}{B^2}
                    \lf(B-\fr12gbq\rg)\fr1{q}[Bb^l-(b+gq)y^l]
a_{tlij}
$$

\ses

\ses

\be
-
\fr g2
q^2y^la_{tlij}
\fr {K}{B^2}
+
\fr1h
\fr {g^2q}4
(b+gq)y^l
a_{tlij}
\fr {K}{B^2}
+
\fr1h
\fr {g^2q^2}4
\fr1{q}[Bb^l-(b+gq)y^l]
a_{tlij}
\fr {K}{B^2},
\ee
\ses
so that
$$
P_{tij}=
-\lf[hq^2+b\lf(b+\fr12gq\rg)\rg]
\fr K{qB}
b^la_{tlij}
+\fr {K}{B^2}
\lf(B-\fr12gbq\rg)\fr1{q}(b+gq)y^l
a_{tlij}
-
\fr g2
q^2y^la_{tlij}
\fr {K}{B^2},
$$
which can be simplified to read
\be
P_{tij}=
\Biggl[
-\lf[hq^2+b\lf(b+\fr12gq\rg)\rg]
b^la_{tlij}
+
\lf(b+\fr12gq\rg)
y^l a_{tlij}
\Biggr]
\fr {K}{qB}.
\ee
The representation (6.77)--(6.78) is valid.

 {%\pgbrk}

Let us find
\ses\\
$$
\fr{qB}K
{\eta}^{nt}  P_n{}^{ij}
=
\Biggl[
-\lf[hq^2+b\lf(b+\fr12gq\rg)\rg]
b^l
+
\lf(b+\fr12gq\rg)
y^l
\Biggr]
a^t{}_{lij}
$$

\ses

$$
-
b^t
\lf(b+\fr12gq\rg)
y^l
b^ha_{hlij}
-\fr1{q^2}v^t
\Biggl[
-\lf[hq^2+b\lf(b+\fr12gq\rg)\rg]
b^l
+
\lf(b+\fr12gq\rg)
y^l
\Biggr]
v^ha_{hlij}
$$

\ses

\ses

\ses

$$
=
\Biggl[
-\lf[hq^2+b\lf(b+\fr12gq\rg)\rg]
b^l
+
\lf(b+\fr12gq\rg)
y^l
\Biggr]
a^t{}_{lij}
-
\lf[
hv^t+
\lf(b+\fr12gq\rg)
b^t
\rg]
y^l
b^ha_{hlij}
$$
\ses
and
$$
\lf(\fr{qB}K\rg)^2
 {\eta}^{nt}  P_n{}^{ij}      P_{tij}
=
hq^2b^l a_{tlij}
\lf[hq^2+b\lf(b+\fr12gq\rg)\rg]
b^h
a^t{}_{hij}
$$

\ses

\ses

$$
-
\lf(b+\fr12gq\rg)
v^l a_{tlij}
\lf[hq^2+b\lf(b+\fr12gq\rg)\rg]
b^h
a^t{}_{hij}
$$

\ses

\ses

$$
+
\Biggl[
-hq^2b^l a_{tlij}
+
\lf(b+\fr12gq\rg)
v^l a_{tlij}
\Biggr]
\lf(b+\fr12gq\rg)
v^h
a^t{}_{hij}
$$

\ses

\ses

$$
+
\Biggl[
-hq^2b^l a_{tlij}
+
\lf(b+\fr12gq\rg)
v^l a_{tlij}
\Biggr]
\lf(b+\fr12gq\rg)
bb^h
a^t{}_{hij}
$$

\ses

     \ses

$$
-
\Biggl[
-hq^2b^l a_{tlij}
+
\lf(b+\fr12gq\rg)
v^l a_{tlij}
\Biggr]
\lf[
hv^t+
\lf(b+\fr12gq\rg)
b^t
\rg]
v^h
b^ua_{uhij}.
$$
So we have
\ses
$$
\lf(\fr{qB}K\rg)^2
 {\eta}^{nt}  P_n{}^{ij}      P_{tij}
=
-
Bb^tb^u
  a_{uv}{}^{ij}
v^vv^l
  a_{tlij}
$$

\ses

\be
+
h^2q^4b^l a_{tlij}
b^h
a^t{}_{hij}
+
\Biggl[
-2
hq^2b^l a_{tlij}
+
\lf(b+\fr12gq\rg)
v^l a_{tlij}
\Biggr]
\lf(b+\fr12gq\rg)
v^h
a^t{}_{hij}.
\ee

\ses

{%\pgbrk}

From (C.14) it follows that
\ses\\
$$
\Rho^{knij}\Rho_{knij}=
2\fr1{K^2}M^{nij}M_{nij}
+2{\cal H}^{nu}  P_u{}^{ij}
{\cal H}_n^tP_{tij}
 +
{\cal H}^{ku} {\cal H}^{nv}  a_{uv}{}^{ij}
 \fr {K^2}B
{\cal H}_k^t {\cal H}_n^l  a_{tlij}
 \fr {K^2}B
$$

\ses

\ses

\ses

$$
=
2
\fr1{K^2}M^{nij}M_{nij}
+
2
{\eta}^{nt}  P_n{}^{ij}   P_{tij}
\fr B{K^2}
 +
{\eta }^{tu} {\eta}^{lv}  a_{uv}{}^{ij}
           a_{tlij}
$$

\ses

\ses

\ses

$$
=
2
\fr1{K^2}M^{nij}M_{nij}
-
2\fr1{q^2}
b^tb^u
  a_{uv}{}^{ij}
v^vv^l
              a_{tlij}
+
a^{tu}
  a_{uv}{}^{ij}
\lf( a^{lv} -b^lb^v-\fr1{q^2}v^vv^l\rg)
  a_{tlij}
 $$

\ses

\ses

$$
+
\fr 2{B}
h^2q^2b^l a_{tlij}
b^h
a^t{}_{hij}
+
\fr 2{Bq^2}
\Biggl[
-2
hq^2b^l a_{tlij}
\lf(b+\fr12gq\rg)
+
\lf(b+\fr12gq\rg)^2
v^l a_{tlij}
\Biggr]
v^h
a^t{}_{hij}
$$

\ses

\ses

$$
-b^tb^u
  a_{uv}{}^{ij}
\lf( a^{lv} -\fr1{q^2}v^vv^l\rg)
  a_{tlij}
-\fr1{q^2}v^uv^t
  a_{uv}{}^{ij}
\lf( a^{lv} -b^lb^v\rg)
  a_{tlij},
 $$
or
\ses\\
$$
\Rho^{knij}\Rho_{knij}=
a^{knij}a_{knij}
+
2
\fr1{K^2}M^{nij}M_{nij}
-
\fr 2{B}
\lf(b+\fr12gq\rg)^2
b^l a_{tlij}
b^h
a^t{}_{hij}
$$

\ses

\ses

$$
-
h\fr 2{B}
b^l a_{tlij}
\lf(b+\fr12gq\rg)
v^h
a^t{}_{hij}
-
h\fr 2{B}
\Biggl[
b^l a_{tlij}
\lf(b+\fr12gq\rg)
+h
v^l a_{tlij}
\Biggr]
v^h
a^t{}_{hij}
$$

\ses

\ses

\ses

$$
=
a^{knij}a_{knij}
+
2
 \fr1{K^2}M^{nij}M_{nij}
-\fr 2{B}
\lf[ hv^l+\lf(b  +\fr g2q    \rg) b^l\rg]
a_l{}^{nij}
\lf[ hv^h+\lf(b  +\fr g2q    \rg)b^h\rg]
a_{hnij}.
$$

{%\pgbrk}

Inserting here
(6.82), we find
\be
\Rho^{knij}\Rho_{knij}=
a^{knij}a_{knij}
+2
\lf(\fr1{h^2}-1\rg){B}
\lf[ hv^l+\lf(b  +\fr g2q    \rg) b^l\rg]
a_l{}^{nij}
\lf[ hv^h+\lf(b  +\fr g2q    \rg)b^h\rg]
a_{hnij}.
\ee

{%\pgbrk}

Using
$$
\zeta^i=\lf[hv^i  +   ( b+\fr12gq)b^i    \rg] \fr S{\sqrt B}
$$
\ses
(see (6.26) and (6.39)),
we arrive at the representation
\be
\Rho^{knij}\Rho_{knij}=
a^{knij}a_{knij}
+\fr 2{S^2}
\lf(\fr1{h^2}-1\rg)
\zeta^l
a_l{}^{nij}
\zeta^h
a_{hnij}
\ee
which is equivalent to (6.87).

{%\pgbrk}

\ses

\ses

\setcounter{equation}{0}

\nin
{ \bf Appendix D: Important coefficients}

\ses

\ses

In processing the involved calculations it is useful to take into account
the equalities
\be
y^n_ma_h{}^m{}_{ij}\zeta^h_k=
-\lf(
\fr1N C_k
-\fr1{K^2}(1-h)y_k
\rg)
M^n{}_{ij}
+
y^n_ma_h{}^m{}_{ij}E^h_k
\ee
($\zeta^h_k$ and $E^h_k$ are  indicated in (6.36)),
\ses
$$
\zeta^h_{km}=\D{\zeta^h_{k}}{y^m}
=
 \fr1{2q}g \eta_{km}  b^h
 J\fr1{\varkappa h}
+\fr1N E^h_kC_m
- (1-h)y_m\fr1{K^2}  E^h_k
$$

\ses

\be
+\fr1N \zeta^h_mC_k  - (1-h)y_k\fr1{K^2}  \zeta^h_m
+\fr1N \zeta^h\D{C_k}{y^m}  -(1-h)g_{km}\fr1{K^2}\zeta^h  +2(1-h)y_ky_m\fr1{K^4}\zeta^h,
\ee
\ses
and
$$
y_h^n\zeta^h_{km}=   \fr1N \de^n_mC_k  - (1-h)y_k\fr1{K^2}  \de^n_m
+\fr1NC^n J^2\eta_{km}
+\fr1N
h^n_k
C_m
- (1-h)y_m\fr1{K^2} \Biggl(\de^n_k-\fr1h \fr1N C_ky^n  \Biggr)
$$

\ses

\ses

\be
-y_m\fr1{K^2} \fr1{h}(1-h)^2y_k\fr1{K^2} y^n
-\fr1h y^n\fr1N
\fr1{K^2}y_mC_k
+(1-h) y^ng_{km}\fr1{K^2}
 +2\fr1h y^n(1-h)y_ky_m\fr1{K^4},
\ee
\ses
together with
$$
y_h^n\zeta^h_{km}M^m{}_{ij}=   \lf(
\fr{ C_k}N  - (1-h)\fr{y_k}{K^2}
\rg)
M^n{}_{ij}
+\lf(
\fr {C^n}N J^2\eta_{km}
+ h^n_k\fr{C_m }N
+
\fr{ y^n}{K^2}(1-h)g_{km}
 \rg)
 M^m{}_{ij}
$$

\ses

\ses

$$
=   \fr1N C_k M^n{}_{ij} - (1-h)y_k\fr1{K^2}M^n{}_{ij}
+\fr1N \eta^n_kC_m  M^m{}_{ij}
+\fr1NC^n g_{km}M^m{}_{ij}
+ y^n(1-h)g_{km}\fr1{K^2}M^m{}_{ij}.
$$
\ses

{%\pgbrk}

With the help of the formula (A.8) of Appendix A, the last representation
can be written merely as
\be
y_h^n\zeta^h_{km}M^m{}_{ij}=
C^n{}_{km}M^m{}_{ij}+(1-h)
\fr1{K^2}\lf(
 y^ng_{km}M^m{}_{ij}
-y_k
 M^n{}_{ij}
\rg).
\ee
\ses
Also,
\be
\fr1{J^2}g_{kl}M^l{}_{ij}=
a_k{}^m{}_{ij}v_m
+\fr1h
\lf[b + \fr12gq\rg]
b_ma_k{}^m{}_{ij}
-\fr1h
\lf[
\fr g{2q}v_k
+
(1-h)
b_k
\rg]
b_my^ta_t{}^m{}_{ij}.
\ee

{%\pgbrk}

%Now we want to find the coefficients
%$y^k_j\zeta^m_{kn}y^n_h$ which enter the transformation (6.78).
If we use the representation (6.36) and apply
 the formulas
(A.24),    (A.25),   (A.29),   (A.30), and    (A.33),
 we obtain
$$
\zeta^m_{nj}
=\lf(\fr1N C_n  -(1-h)\fr1{K^2}y_n\rg)    \zeta^m_j
+
\lf(\fr1N C_j  -(1-h)\fr1{K^2}y_j\rg)    \zeta^m_n
+ \fr g{2q} \eta_{nj}b^m  J\fr1{\varkappa h}
$$

\ses

$$
-\lf(\fr1N C_n -(1-h)\fr1{K^2}y_n    \rg)         \fr1NC_j\zeta^m
+
\fr1{K^2}(1-h)
\lf(\fr1N C_n -(1-h)\fr1{K^2}y_n   \rg) y_j
\zeta^m
$$

\ses

$$
+\lf(
- \fr1 Kl_n \fr1N C_j-\fr1 Kl_j\fr1N C_n
-\fr1 B\fr 12\fr{gb}{q}     \eta_{nj}+\fr{2}N\fr1N C_nC_j\rg)
\zeta^m
-(1-h)  ( g_{nj}-2l_nl_j)\fr1{K^2}\zeta^m.
  $$
\ses
Simplifying yields
the representation
$$
\zeta^m_{nj}
=\lf(\fr1N C_n  -(1-h)\fr1{K^2}y_n\rg)    \zeta^m_j
+
\lf(\fr1N C_j  -(1-h)\fr1{K^2}y_j\rg)    \zeta^m_n
+ \fr g{2q} \eta_{nj}b^m  J\fr1{\varkappa h}
$$

\ses

$$
-h\fr1{K^2}      \fr1N(y_nC_j+y_jC_n)
\zeta^m
+
h(1-h)
\fr1{K^4}y_n   y_j
\zeta^m
$$

\ses

\be
+\lf(
-\fr1 B\fr 12\fr{gb}{q}     \eta_{nj}+\fr{1}N\fr1N C_nC_j\rg)
\zeta^m
-(1-h)   h_{nj}\fr1{K^2}\zeta^m,
\ee
\ses
from which it follows that
$$
y^k_j\zeta^m_{kn}y^n_h=
\lf(\fr1N C_ny^n_h  -\fr{1-h}{hS^2}  \zeta_h\rg)    \de^m_j
+\lf(\fr1N C_ny^n_j  -\fr{1-h}{hS^2}  \zeta_j\rg)    \de^m_h
+ \fr g{2q} b^m  \fr J{\varkappa h}
\eta_{nk}y^n_hy^k_j
$$

\ses

$$
-      \fr1NC_ky^k_j\fr{1}{S^2}  \zeta_h\zeta^m
-      \fr1NC_ky^k_h\fr{1}{S^2}  \zeta_j\zeta^m
+
h(1-h)
\fr1{h^2S^4}  \zeta_h
  \zeta_j
\zeta^m
$$

\ses

$$
-\fr1 B\fr 12\fr{gb}{q}     \eta_{nk}
y^n_hy^k_j\zeta^m
+\fr{1}N\fr1N C_nC_ky^n_hy^k_j
\zeta^m
 -(1-h)   \eta_{nk}\fr1{B}y^n_hy^k_j\zeta^m
  -(1-h) \fr2g\fr2g\fr1N\fr1N  C_nC_k    y^n_hy^k_j\zeta^m.
     $$

{%\pgbrk}

\nin
Eventually, we arrive at
$$
y^k_j\zeta^m_{kn}y^n_h=
Q_h  \de^m_j+Q_j  \de^m_h
+ \fr g{2q} b^m  \fr {\varkappa }{Jh}
\eta_{hj}
+
h(1-h)
\fr1{h^2S^4}  \zeta_h
  \zeta_j
\zeta^m
$$

\ses

$$
- M_j
     \fr{1}{S^2}  \zeta_h\zeta^m
- M_h
\fr{1}{S^2}  \zeta_j\zeta^m
-\fr1 B\fr 12\fr{gb}{q}     \eta_{hj}
\fr {\varkappa }{J} \fr {\varkappa }{J}
  \zeta^m
+
M_hM_j
\zeta^m
$$

\ses

\be
  -(1-h)   \eta_{nk}\fr1{B}y^n_hy^k_j\zeta^m
  -(1-h) \fr2g\fr2g\fr1N\fr1N  C_nC_k    y^n_hy^k_j\zeta^m,
\ee
\ses
where
$$
Q_h=M_h
          -(1-h)\fr{1}{hS^2}  \zeta_h
$$
with
$$
M_h=\Biggl[
\fr1NC_h+\fr{gq}{2B}
 \lf((h-1) b_h  -          \fr {gT}{2} [\zeta_h-\zeta^lb_lb_h]  \rg)
\Biggr]
\fr {\varkappa }J.
 $$
Using the equalities
$$
\fr1NC_h=\fr 1hg\fr1{2qB}(hq^2b_h- bhv_h)
=\fr 1hg\fr1{2qB}[hq^2+ b( b+\fr12gq]b_h
-g\fr1{2qB}
 b\zeta_h\fr {\varkappa }J
$$
yields
$$
M_h=
\Biggl[
-g\fr1{2qB}
 b\zeta_h\fr {\varkappa }J
 +\fr{gq}{2B}
 \lf(h b_h  -          \fr {gT}{2} \zeta_h +\fr1{hq^2}(b+\fr12gq)^2 \rg)
\Biggr]
\fr {\varkappa }J
$$

\ses

$$
=
\Biggl[
-\fr g{2qB}
 b\zeta_h\fr {\varkappa }J
 +\fr{gq}{2B}
 \lf(\fr B{hq^2} b_h  -        \fr {\varkappa }J  \fr {g}{2q} \zeta_h  \rg)
\Biggr]
\fr {\varkappa }J,
 $$
\ses
or
$$
M_h=
-\fr g{2qB}
( b+\fr12gq)\zeta_h\fr {\varkappa }J  \fr {\varkappa }J
 +\fr{g}{2qh}
  b_h
\fr {\varkappa }J.
 $$
Here we can apply (6.39),
which yields
\be
M_h=
-\fr g{2qhS^2}
( b+\fr12gq)\zeta_h
 +\fr{g}{2qh}
  b_h
\fr {\varkappa }J.
\ee

By means of the transition rule (4.6) the tensor
$E_k{}^n{}_{ij}$ can be transformed into the tensor
$L_k{}^n{}_{ij}~:= y^t_k\zeta^n_lE_t{}^l{}_{ij}$
of the Riemannian space $\cR^N$, which yields
\be
L_k{}^n{}_{ij}= a_k{}^n{}_{ij}-\ga^n_{kt} \zeta^h a_h{}^t{}_{ij}
\ee
with the coefficients $ \ga^n_{kt}=y^r_ky^s_t\zeta^n_{rs}$
given by (D.7).

   {%\pgbrk}

The coefficients (6.45) can be transformed as follows:
$$
y^k_i=
\Biggl[
\Bigl(b_i  -        \fr g{2qh}v_i
\Bigr)
b^k+
\fr1h
\Bigl(\de^k_i- b_ib^k\Bigr)
\Biggr]
 \fr {h\varkappa }J
$$

\ses

$$
+
\fr1B
\Biggl\{
\fr{1-h}{h}  \lf  [hv_i+(b+\fr12gq)b_i\rg]
+
   \fr g{2h^2q}
(b+\fr12gq)
\lf   [hv_i+(b+\fr12gq)b_i\rg]
-
   \fr g{2h^2q}
b_i
\Biggr\}
 \fr {h\varkappa }J
 y^k,
 $$
\ses
getting

$$
y^k_i=
\Biggl[
\Bigl(b_i  -        \fr g{2qh}v_i
\Bigr)
b^k+
\fr1h
\Bigl(\de^k_i- b_ib^k\Bigr)
\Biggr]
 \fr {h\varkappa }J
$$

\ses

\be
+
\fr1{Bhq}
\Biggl[
(h-h^2)q
+
   \fr g{2}
(b+\fr12gq)
 \Biggr]
 \fr {h\varkappa }J
v_i
 y^k
+
\fr1{B}
\Biggl(
\fr1h   (b+\fr12gq)
-b-gq
\Biggr)
 \fr {h\varkappa }J
b_i
 y^k,
\ee
\ses
or
\be
y^n_m=
\de^n_m
 \fr {\varkappa }J
-
v_m
C^n\fr{J^2}{N}
 \fr {\varkappa }J
-
(1-h)
b_m
 C^n\fr{2qJ^2}{gN}
 \fr {\varkappa }J
-\fr1{B}
\lf((1-h)v_m+\fr{gq}2b_m\rg)
 \fr {\varkappa }J
 y^n,
\ee
\ses
which
 entails
\be
y^k_i b^i=
\Biggl[
b^k
+
\fr1{B}\Biggl(
\fr{1}{h}
\lf(b+\fr12gq\rg)
-b -gq
\Biggr)
y^k
\Biggr]
 \fr {h\varkappa }J.
\ee

{%\pgbrk}

\ses

\setcounter{equation}{0}

\nin
{ \bf Appendix E:
  Fixed  tangent space
  of the ${\mathbf\cF\cF^{PD}_{g} } $-space}

\ses

\ses

Let us introduce the {\it orthonormal frame}
  $h^p_i(x)$ of the input  Riemannian  metric tensor
$a_{ij}(x)$:
\be
a_{ij}=e_{pq}h^p_ih^q_j,
\ee
where $\{e_{pq}\}$ is the  Euclidean  diagonal:
$
e_{pq}=\text{diagonal}(+...+);
$
the indices
$p,q,...$ will be specified on the range $1,...,N$;
and
the indices
$a,b,...$ on the range $1,...,N-1$.
 Denote by $h_p^i$ the reciprocal frame,
 so that
 $
 h_p^jh^p_i=\de^j_i.
$
At any fixed point $x$, we can represent the tangent vectors $y$ by their frame-components:
\be
R^p=h^p_iy^i,
\ee
and  use the components
\be
g_{pq}=h_p^ih_q^jg_{ij}
\ee
of the Finslerian metric tensor $g_{ij}$.

In the ${\mathbf\cF\cF^{PD}_{g} } $-space,
it is convenient to specify the frame such that the $N$-th component $h^{N}_i(x)$
becomes collinear to the input vector field $b_i(x)$.
Under these conditions,  the  1-form $b$ reads merely
$b=z$
and we have
$  b^p=\{0,0,...,1\} $
and
$   b_p=\{0,0,...,1\}.   $
 We obtain the decomposition
$
R^p=\{R^a,R^{N}\},
$
together with
$    q^2=e_{ad}R^aR^d.    $
Also the notation
\be
R^{N}=z
\ee
will be used.

In any fixed  tangent space $T_xM$
we can obtain  the covariant components $R_p=h_p^iy_i$ through the definition
$$
R_p~:=\fr12\D{K^2(g;R)}{R^p}.
$$
With the help of  (6.1)-(6.5)
we find
\be
R_a=
e_{ab}R^b
J^2,
\qquad
R_N=(z+gq)J^2.
\ee
For the respective Finsleroid metric tensor components
$$
g_{pq}(g;R)
~:=\fr12\,
\fr{\prtl^2K^2(g;R)}{\prtl R^p\prtl R^q}
=\fr{\prtl R_p(g;R)}{\prtl R^q}
$$
 we obtain
\be
 g_{NN}(g;R)=[(z+gq)^2+q^2]
J^2,
\qquad
g_{Na}(g;R)=gq
e_{ab}R^b
J^2,
\ee

\ses

\be
g_{ab}(g;R)=
J^2
e_{ab}-g\fr{
e_{ad}R^d
e_{be}R^e
z
}{q}
J^2.
\ee

{%\pgbrk}

The components of the inverted metric tensor read
\ses
\be
g^{NN}(g;R)=(z^2+q^2)
\fr1{K^2},
\qquad
g^{Na}(g;R)=-gqR^a
\fr1{K^2},
\ee
\ses
\be
g^{ab}(g;R)=
\fr{B}{K^2}
e^{ab}+g(z+gq)\fr{R^aR^b}{q}
\fr1{K^2}.
\ee

It can readily be verified that
\be
\det(g_{pq})=J^{2N}>0.
\ee
The above formulas are valid at  an arbitrary dimension $N\ge2$.

{%\pgbrk}

\vskip 1cm

\def\bibit[#1]#2\par{\rm\noindent\parskip1pt
                     \parbox[t]{.05\textwidth}{\mbox{}\hfill[#1]}\hfill
                     \parbox[t]{.925\textwidth}{\baselineskip11pt#2}\par}

\bc {\bf  References}
\ec

\ses

\bibit[1] H. Rund, \it The Differential Geometry of Finsler
 Spaces, \rm Springer, Berlin 1959.

\bibit[2] D. Bao, S. S. Chern, and Z. Shen, {\it  An
Introduction to Riemann-Finsler Geometry,}  Springer, N.Y., Berlin 2000.

\bibit[3] L. Kozma and L. Tam{\' a}ssy,
Finsler geometry without line elements faced to applications,
{   \it Rep. Math. Phys.} {\bf 51} (2003), 233--250.

\bibit[4] L. Tam{\' a}ssy,
 Metrical almost linear connections in $TM$ for Randers spaces,
{\it Bull. Soc. Sci. Lett. Lodz Ser. Rech. Deform } {\bf 51} (2006), 147-152.

\bibit[5] Z. L. Szab{\' o},  All regular Landsberg metrics are Berwald,
{\it Ann Glob Anal Geom  } {\bf 34} (2008), 381-386.

\bibit[6] L. Tam{\' a}ssy,  Angle in Minkowski and Finsler spaces,
{\it Bull. Soc. Sci. Lett. Lodz Ser. Rech. Deform } {\bf 49} (2006), 7-14.

\bibit[7] G. S. Asanov,  Finsleroid-space supplemented by angle,
  {\it  arXiv}: 0310019 [math-ph] (2003).

\bibit[8] G. S. Asanov,
{Finsleroid space with angle and scalar product,
   \it  Publ. Math. Debrecen} {\bf 67} (2005),  209-252.

\bibit[9] G. S. Asanov,
Finslerian angle-preserving  connection   in two-dimensional  case.
Regular realization,
{\it  arXiv:} 0909.1641v1 [math.DG],  (2009).

\bibit[10] G. S. Asanov,
{   Finsler cases of GF-space,
   \it  Aeq. Math.} {\bf 49} (1995),  234-251.

 \end{document}